\documentclass[11pt,leqno]{amsart}
\usepackage{amsmath,amsfonts,amssymb,amscd,amsthm,amsbsy,upref}
\theoremstyle{plain}
\usepackage[dvipsnames]{xcolor}   % to write in color

\newtheorem{thm}{Theorem}[section]
\newtheorem{lem}[thm]{Lemma}
\newtheorem{cor}[thm]{Corollary}
\newtheorem{prop}[thm]{Proposition}

\newtheorem{probs}[thm]{Problems}
\theoremstyle{definition}
\newtheorem{defin}[thm]{Definition}
\newtheorem{ex}[thm]{Example}
\newtheorem{exs}[thm]{Examples}
\newtheorem*{rem}{Remark}

\newtheorem{remark}[thm]{Remark}

\newtheorem*{thmA}{Theorem A}
\newtheorem*{thmB}{Theorem B}
\newtheorem*{thmC}{Theorem C}

\usepackage{amssymb}

\newcommand{\N}{\mathbb{N}}

\newcommand{\E}{\mathbb{E}}

\newcommand{\cA}{\mathcal A}

\newcommand{\cC}{\mathcal C}

\newcommand{\cF}{\mathcal F}
\newcommand{\cE}{\mathcal E}
\newcommand{\cG}{\mathcal G}

\newcommand{\cM}{\mathcal M}

\newcommand{\cT}{\mathcal T}
\newcommand{\cS}{\mathcal S}

\newcommand{\cU}{\mathcal U}

\newcommand{\cSB}{\mathcal{SB}}

\newcommand{\xb}{\overline{x}}
\newcommand{\yb}{\overline{y}}

\newcommand{\Ab}{\overline{A}}

\newcommand{\mt}{\tilde m}

\newcommand{\xt}{\tilde x}
 
\newcommand{\zt}{\tilde z}
\newcommand{\lt}{\tilde l}

\newcommand{\At}{\tilde A}
\newcommand{\Bt}{\tilde B}
\newcommand{\Ct}{\tilde C}
\newcommand{\Dt}{\tilde D}
\newcommand{\at}{\tilde a}

\newcommand{\keq}{\!=\!}
 \newcommand{\kneq}{\!\neq\!}
\newcommand{\kleq}{\!\leq\!}
\newcommand{\kge}{\!\ge\!}
\newcommand{\kgeq}{\!\ge\!}
\newcommand{\kle}{\!<\!}
\newcommand{\kgr}{\!>\!}
\newcommand\kin{\!\in\!}
\newcommand{\ksubset}{\!\subset\!}
\newcommand{\ksupset}{\!\supset\!}
\newcommand{\ksetminus}{\!\setminus\!}
\newcommand{\kplus}{\!+\!}
\newcommand{\kminus}{\!-\!}

\newcommand{\supp}{\text{\rm supp}}
\newcommand{\ran}{\text{\rm ran}}

\newcommand{\spa}{\text{\rm span}}

\newcommand{\MAX}{\text{\rm MAX}}
\newcommand{\CB}{\text{{\rm CB}}}

\newcommand{\Sz}{\text{\rm Sz}}

\newcommand{\vp}{\varepsilon}

\newcommand{\Cp}{\text{\rm Cp}}

\newcommand{\ie}{\textit{i.e.,}\ }

\begin{document}
\allowdisplaybreaks

\title{A metric interpretation of reflexivity for Banach spaces}

\author{P. Motakis}
\address{Department of Mathematics, Texas A\&M University, College
  Station, TX 77843, USA }
\email{Pavlos@math.tamu.edu}
\author{Th. Schlumprecht}
\address{Department of Mathematics, Texas A\&M University, College
  Station, TX 77843, USA and Faculty of Electrical Engineering, Czech
  Technical University in Prague,  Zikova 4, 166 27, Prague}
\email{schlump@math.tamu.edu}

 \thanks{The second named  author  was supported by the National Science Foundation under Grant Numbers DMS--1160633
 and DMS--1464713.}

\begin{abstract} We define 
two metrics $d_{1,\alpha}$ and $d_{\infty,\alpha}$ on each  Schreier family $\cS_\alpha$, $\alpha<\omega_1$, with which we prove the following  metric characterization
of reflexivity of a Banach space $X$:
$X$ is reflexive if and only if there is an $\alpha<\omega_1$, so that there is no mapping 
$\Phi:\cS_\alpha\to X$ for which
$$  cd_{\infty,\alpha}(A,B)\le \|\Phi(A)-\Phi(B)\|\le C  d_{1,\alpha}(A,B) \text{ for all $A,B\in\cS_\alpha$.}$$
Secondly, we prove for  separable and reflexive Banach spaces $X$, and certain countable ordinals $\alpha$ that 
$\max(\Sz(X),\Sz(X^*))\le \alpha$ if and only if $(\cS_\alpha, d_{1,\alpha})$ does not bi-Lipschitzly embed into $X$. Here $\Sz(Y)$ denotes the Szlenk index of a Banach space $Y$.
\end{abstract}
\maketitle

\tableofcontents

\section{Introduction and statement of the main results}\label{S:1}

In this paper we seek a  metric characterization of  reflexivity of Banach spaces. By a metric characterization of   a property of a Banach  space
we mean a characterization which refers only to the metric structure of  that space but not its linear structure.
In 1976 Ribe \cite{Ri} showed  that  two Banach spaces, which are uniformly homeomorphic,   have uniformly linearly isomorphic finite-dimensional subspaces. In particular this means that the finite dimensional or {\em local}  properties  of a Banach space are determined by its metric structure.
 Based on this result Bourgain \cite{Bo} suggested the ``Ribe Program'', which asks  to find  metric descriptions of  finite-dimensional invariants of Banach spaces. 
 In \cite{Bo}  he  proved the following characterization of super reflexivity: 
  a Banach space $X$ is super reflexive if and only if the binary trees $B_n$ of length at most  $n$, $n\in\N$, endowed with their graph metric, are not uniformly bi-Lipschitzly embedded into $X$.
  A binary tree of length at most  $n$, is the set $B_n=\bigcup_{k=0}^n\{-1,1\}^k$, with the {\em graph } or {\em shortest path metric} 
  $d(\sigma,\sigma')=i+j  -2\max\{ t\ge 0: \sigma_s=\sigma'_s: s=1,2,\ldots,t\} $, for $\sigma=(\sigma_s)_{s=1}^i\not=\sigma'=(\sigma_s')_{s=1}^j$ in $\bigcup_{k=0}^n\{-1,1\}^k$. A new and shorter proof of this result was 
  recently obtained by Kloeckner in  \cite{Kl}. In \cite{Ba1} 
   Baudier extended this result and proved that a Banach space $X$ is super reflexive, if and only if the infinite binary tree $B_\infty=\bigcup_{n=0}^\infty \{-1,1\}^n$ (with the graph distance)
   does not bi-Lipschitzly embed into $X$. Nowadays this result can be deduced from Bourgain's result and a result of  Ostrovskii's  \cite[Theorem 1.2]{Os1} which states that a locally finite metric space $A$ embeds
   bi-Lipschitzly into a Banach space $X$ if all of its finite subsets uniformly 
   bi-Lipschitzly embed into $X$.
   In \cite{JS} Johnson and Schechtman characterized superflexivity, using the  {\em Diamond Graphs}, $D_n$, $n\in\N$, and proved that Banach space $X$ is super reflexive if and only if  the $D_n$, $n\in\N$, do  not 
   uniformly bi-Lipschitzly embed into $X$. There are several other  local properties, \ie properties of the finite dimensional subspaces of Banach spaces,  for which metric characterizations were found. The following are some examples:
   Bourgain, Milman and Wolfson \cite{BMW} characterized having non trivial type  using {\em Hammond cubes} (the sets $B_n$, together with the $\ell_1$-norm),  and Mendel and Naor \cite{MN1,MN2} present
   metric characterizations of Banach spaces with type $p$, $1<p\le 2$, and cotype $q$, $2\le q<\infty$.
   % and Lee, Naor, Peres establish a metric characterization of $q$-convexity \cite{LNP}. 
   For a more extensive account on the Ribe program we would like refer the reader  to the  survey articles \cite{Ball,Na}
   and the book \cite{Os2}.
  
  Instead of only asking for metric characterizations of local properties, one can also ask for metric characterizations of other properties of Banach space, properties which might not be determined by 
   the finite dimensional subspaces. A result in this direction was obtained by Baudier, Kalton and Lancien  in \cite{BKL}. They showed that a  reflexive Banach space $X$ has a renorming, which is  {\em asymptoticaly uniformly convex} (AUC)
   and {\em asymptoticaly uniformly  smooth}  (AUS), if and only if the  countably  branching trees of length $n\in\N$, $T_n$ do not uniformly bi-Lipschitzly into $X$.
   Here $T_n=\bigcup_{k=0}^n\N^k$, together with the graph metric, \ie $d( a, b)= i+j-\max\{ t\ge 0: a_s=b_s, s=1,2\ldots t\} $, for $a=(a_1,a_2,\ldots a_i)\not=  b=(b_1,b_2,\ldots b_j)$ in $T^n$.
    Among the many equivalent  conditions for a reflexive Banach space $X$ to be   (AUC)- and (AUS)-renormable (see \cite{OS2}) one of them states that  $\Sz(X)=\Sz(X^*)=\omega$,
    where $\Sz(Z)$ denotes the {\em Szlenk index} of a Banach space $Z$ (see Section \ref{S:6} for the definition and properties of the Szlenk index).  In \cite{DKLR} 
    Dilworth, Kutzarova, Lancien and Randrianarivony, showed that a separable  Banach space $X$ is reflexive and (AUC)- and (AUS)-renormable if and only $X$ admits an equivalent norm for which $X$ has {\em  Rolewicz's $\beta$-property}.
    According to \cite{Ku1}  a Banach space $X$ has Rolewicz's $\beta$-property if and only  if 
  $$\bar \beta_X(t)= 1-\sup\Big\{ \inf\Big\{\frac{\|x- x_n\|}{2}: n\ge1\Big\}: (x_n)_{n=1}^\infty\subset B_X,\,\, \text{sep} \big[(x_n)\big]\ge t,\,\ x\in B_X\Big\}>0,$$
  for all $t>0$, where $\text{sep}\big[(z_n)\big]=\inf_{m\not=n} \|z_m-z_n\|$,  for a sequence $(z_n)\subset X$.   The function $\bar\beta_X$ is called the {\em $\beta$-modulus of $X$}.   
    Using the equivalence between the  positivity of the $\beta$-modulus and the property that a separable Banach space is reflexive and (AUC)- and (AUS)-renormable, Baudier \cite{Ba3} was able to establish a new and shorter proof of the above cited result from   \cite{BKL}.  Metric descriptions of other non-local Banach space properties, for example the Radon-Nikod\'ym  
   Property, can be found in \cite{Os3}.
   
   In our paper we would like to concentrate on metric descriptions of the property that a Banach space is reflexive, and subclasses of reflexive Banach spaces. In \cite{Os3},
   Ostrovskii established a {\em submetric} characterization of reflexivity. Let $T$ be the set of all pairs $(x,y)$ in $\ell_1\times \ell_1$, for which
   $\|x-y\|_{1}\le 2 \|x-y\|_s$, where $\|\cdot\|_1$ denotes the usual norm on $\ell_1$ and  $\|\cdot\|_s$ denotes the {\em summing norm}, \ie 
   $\|z\|_s=\sup_{k\in\N}\big|\sum_{j=1}^k z_j\big|$, for $z=(z_j)\in \ell_1$.  \cite[Theorem 3.1]{Os3} states that  a Banach space $X$ is not reflexive if and only if
    there is a map $f:\ell_1\to X$ and a number $0<c\le 1$,  so that $c\|x-y\|_1\le \| f(x)-f(y)\|\le \|x-y\|_1$ for all $(x,y)\in T$. In Section \ref{S:13} we will formulate a
    similar result, using a discrete  subset of $\ell_1\times\ell_1$,    witnessing the same phenomena.
     Recently  Proch\'azka \cite[Theorem 3]{Pr} formulated an interesting metric description of reflexivity. He constructed a unifomly discrete metric space $M_R$
      with the following properties: If $M_R$ bi-Lipschitzly embeds into a Banach space $X$ with distortion less than $2$, then $X$ is non reflexive. The {\em distortion}
      of a bi-Lipschitz embedding $f$ of one metric into another is the product of the Lipschitz constant of $f$ and the   
      Lipschitz constant of $f^{-1}$. Conversely, if $X$ is non reflexive, then there exists a renorming $|\cdot|$  of $X$, so that $M_R$   embeds into $(X,|\cdot|)$ isometrically.

  Our paper has the  goal to find a metric characterization of  reflexivity. An optimal result would be a statement, similar to Bourgain's result, of the form
   ``all members of a certain family $(M_i)_{i\in I}$ of metric spaces embed uniformly bi-Lipschitzly into a space $X$ if and only if $X$ is not reflexive''.
   In the language, introduced by Ostrovskii \cite{OS2}, this would mean that $(M_i)_{i\in I}$ {\em is a family of test spaces for reflexivity}.
    Instead, our result will be of the form (see Theorem A below),``there is a family of sets $(M_i)_{i\in I}$, and for $i\in I$, there are metrics 
    $d_{\infty,i}$ and $d_{1,i}$ on $M_i$, with the property, that  a given space $X$ is non reflexive if and only if there are injections $\Phi_i: M_i\to X$  and $0<c\le 1$ so that 
    $cd_{\infty,i}(x,y)\le \|x-y\|\le d_{1,i}(x,y)$, for all $x,y\in M_i$''. In Section \ref{S:13}  we will  discuss the difficulties, in arriving to a characterization of reflexivity of the first form.  
 Nevertheless, if we restrict ourselves to the class of reflexive spaces we 
   arrive to  a metric characterization     for the 
   {\em complexity} of a given space, which we measure by the Szlenk index, using  test spaces.  Roughly speaking,  the higher 
  the Szlenk index is of a given Banach space, the more averages of a given weakly null sequence one has to take to arrive to a norm null sequence. For a precise formulation of this  statement 
   we refer to Theorem \ref{T:6.3}.  For  the class of separable and reflexive spaces we will introduce an uncountable  family of metric spaces  $(M_\alpha)_{\alpha<\omega_1}$ for which we will show that the higher 
   the complexity of a given reflexive and separable space $X$ or its dual $X^*$ is, the more members of  $(M_\alpha)_{\alpha<\omega_1}$ can be  uniformly bi-Lipschitzly embedded into $X$.

The statements of  our main results are as follows. The  definition of the {\em Schreier families $\cS_\alpha\subset[\N]^{<\omega}$,} for $\alpha<\omega_1$,  will be recalled in Section \ref{S:2}, the {\em Szlenk index }$\Sz(X)$ for a Banach space $X$ in Section \ref{S:6},   and the two 
 metrics $d_{1,\alpha}$ and $d_{\infty,\alpha}$ on $\cS_\alpha$ will be defined in Section \ref{S:8}. The statements of  our main results are as follows.
  \begin{thmA}  A separable Banach space $X$ is reflexive if and only if there is an $\alpha<\omega_1$ for which  there does not exist a map
  $\Phi: \cS_\alpha\to X$, with the property that for some numbers $C\ge c>0$
  \begin{equation}\label{E:1.1}
    cd_{\infty,\alpha}(A,B)\le \|\Phi(A)-\Phi(B)\|\le C  d_{1,\alpha}(A,B) \text{ for all $A,B\in\cS_\alpha$.}
    \end{equation}
 \end{thmA}
 \begin{defin} Assume that $X$ is a Banach space, $\alpha<\omega_1$, and  $C\ge c>0$. We call  a map
 $\Phi: \cS_\alpha\to X$, with the property that  for all $A,B\in\cS_\alpha$,
   \begin{equation}\label{E:1.2}
    cd_{\infty,\alpha}(A,B)\le \|\Phi(A)-\Phi(B)\|\le C  d_{1,\alpha}(A,B)   \end{equation}
    a
    {\em   $c$-lower-$d_{\alpha,\infty}$ and $C$-upper-$d_{\alpha,1}$ embedding  of $\cS_\alpha$ into $X$}. If $\cA$ is a subset of $\cS_\alpha$, and 
    $\Phi:\cA\to X$, is a map which satisfies \eqref{E:1.2} for all $A,B\in \cA$, we call it a
    {\em   $c$-lower-$d_{\alpha,\infty}$ and $C$-upper-$d_{\alpha,1}$ embedding  of $\cA$ into $X$}.
 \end{defin}
 
  Our next result extends one direction (the ``easy direction'') of the main result of \cite{BKL} to  spaces with higher order Szlenk indices. As in \cite{BKL} reflexivity is not needed here.
 \begin{thmB} Assume that $X$ is a separable Banach space and that $\max\big(\Sz(X),\Sz(X^*)\big)> \omega^\alpha$, for some countable ordinal $\alpha$.
 Then $(\cS_\alpha,d_{1,\alpha})$ embeds bi-Lipschitzly into $X$ and $X^*$.
 \end{thmB}
 
We will deduce one direction of Theorem A from James's characterization of reflexive Banach spaces \cite{Ja0}, and show  that for any non reflexive Banach space
$X$   and any $\alpha<\omega_1$ there is a map $\Phi_\alpha:\cS_\alpha \to X$ which satisfies \eqref{E:1.1}.
The converse will follow from the following  result.
\begin{thmC} Assume that $X$ is a reflexive and separable Banach space. Let $  \xi<\omega_1$ and put $\beta=\omega^{\omega^\xi}$.
If  for some numbers $C>c>0$ there exists
a  $c$-lower-$d_{\alpha,\infty}$ and $C$-upper-$d_{\alpha,1}$ embedding of $\cS_{\beta^2}$ into $X$,
    then $\Sz(X)> \omega^{\beta}$ or $\Sz(X^*)>\beta$.
\end{thmC}
Theorem C, and thus the missing part of Theorem A, will be shown in Section  \ref{S:12}, Theorem \ref{T:12.6}.
Combining Theorems B and C, we obtain the following characterization of certain bounds of the   Szlenk index of $X$ and its dual $X^*$. This result
extends the main result of \cite{BKL} to  separable and reflexive  Banach spaces  with higher order Szlenk indices.
\begin{cor}\label{C:1.1} Assume that $\omega<\alpha<\omega_1$ is an ordinal for which $\omega^\alpha=\alpha$. Then the following statements are equivalent
for a separable and reflexive space $X$
\begin{enumerate}
\item[a)]  $\max(\Sz(X),\Sz(X^*))\le \alpha$,
\item[b)] $(\cS_\alpha,d_1)$ is not bi-Lipschitzly embeddable into $X$.
\end{enumerate}
\end{cor}  
 Corollary \ref{C:1.1} and a result in \cite{OSZ1} yield the  following corollary.   
We thank  Christian Rosendal who pointed it out to us.
\begin{cor}\label{C:1.2} If $\alpha\kle\omega_1$, with $\alpha=\omega^\alpha$, the class of all Banach spaces $X$ for which
$\max(\Sz(X),\Sz(X^*))\kleq\alpha$ is Borel in the Effros-Borel structure of closed subspaces of $C[0,1]$.
\end{cor}

A proof of Corollaries \ref{C:1.1} and \ref{C:1.2} will be presented  at the end of Section \ref{S:12}.
 For the proof of our main results  we will need to introduce some  notation and to make several preliminary observations. The reader who is at first only interested to understand our main results  will only need the definition of 
 the {\em Schreier families} $\cS_\alpha$, $\alpha<\omega_1$, given in Subsection \ref{SS:2.2},  the definition of {\em  repeated averages}  stated in the beginning of of Section \ref{S:4}, and the definition of the two metrics $d_{1,\alpha}$ and $d_{\infty,\alpha}$ on $\cS_\alpha$ introduced in Section \ref{S:8}.
 
\section{Regular Families,  Schreier Families   and  Fine Schreier Families}\label{S:2}

In this section we will first recall the definition of  general {\em  Regular Subfamilies of }$[\N]^{<\omega}$. Then we recall the definition of  the {\em Schreier Families}  $\cS_\alpha$ 
and the {\em Fine Schreier Families } $\cF_{\beta,\alpha}$,  $\alpha\le \beta<\omega_1$ \cite{AA}. The recursive definition of both families depends on choosing for 
every limit ordinal a sequence $(\alpha_n)$, which increases to $\alpha$. In order to ensure that later our proof  will work out, we will need 
 $(\alpha_n)$ to satisfy certain conditions.

\subsection{Regular Families in $[\N]^{<\omega}$}\label{SS:2.1}
For a set $S$ we  denote the subsets, the finite subsets, and the  countably infinite subsets by $[S]$,   $[S]^{<\omega}$, and $[S]^\omega$. 
 We always write subsets of $\N$ in increasing order. Thus,  if we write $A=\{a_1,a_2,\ldots a_n\}\in[\N]^{<\omega}$, or  $A=\{a_1,a_2,\ldots \}\kin[\N]^{\omega}$ we always assume that $a_1\kle a_2\kle \ldots $.
Identifying the elements  of $[\N]$ in the usual way with elements of $\{0,1\}^\omega$, we consider on $[\N]$ the product topology of the discrete topology on $\{0,1\}$. Note that it   follows for a sequence 
$(A_n)\ksubset [\N]^{<\omega}$ and  $A\kin[\N]^{<\omega}$, that $(A_n)$ converges to $A$ if and only if for all
$k\kge \max A$ there is an $m$ so that $A_n\cap[1,k]\keq A$, for all $n\kge m$.

 For $A\kin [\N]^{<\omega}$ and $B\kin[\N]$ we write $A\kle B$, if $\max(A)\kle \min(B)$.
  As a matter of convention we put $\max(\emptyset)=0$ and $\min(\emptyset)=\infty$,
  and thus $A\kle\emptyset$ and $\emptyset\kgr  A$ is true for all $A \in [\N]^{<\omega}$.
  For $m\in\N$ we write $m\le A$ or $m<A$, if 
 $m\le \min(A)$, or $m<\min(A)$, respectively.

We denote by $\preceq$ the partial order of {\em extension} on $[\N]^{<\omega}$, \ie $A=\{a_1,a_2,\ldots a_l\}\preceq B=\{b_1,b_2,\ldots b_m\} $ if 
$l\le m$ and $a_i=b_i$, for $i=1,2\ldots l$, and we write  $A\prec B$ if $A\preceq B$  and $A\not= B$.

 We say that $\cF\ksubset[\N]^{<\omega}$ is 
 {\em closed under taking  restrictions} or {\em  is a subtree of $[\N]^{<\omega}$} if $A\kin\cF$ whenever $A\prec B$ and  $B\kin\cF$,
\emph{hereditary} if $A\kin\cF$ whenever $A\ksubset B$ and $B\kin\cF$,
and $\cF$ is called \emph{compact }if it is compact in the product
topology. Note that a  family  which is closed under restrictions is compact if and only if it  is {\em well founded}, \ie if it does  not
contain strictly ascending chains with respect to extensions.
Given $n,\ a_1\kle \dots\kle
a_n,\ b_1\kle \dots\kle b_n$ in $\N$ we say that $\{b_1,\dots,b_n\}$
is \emph{a spread }of $\{a_1,\dots,a_n\}$ if $a_i\kleq b_i$ for $i\keq1,\dots,n$. A family $\cF\ksubset[\N]^{<\omega}$ is called
\emph{spreading }if every spread of every element of $\cF$ is also in
$\cF$.  We sometimes have to pass from a family $\cF\subset [\N]^{<\omega}$ to the subfamily 
$\cF\cap[N]^{\omega}= \{ A\in \cF: A\subset N\}$, where $N\subset \N$ is infinite. Of course we might lose the property that $\cF$ is spreading.
 Nevertheless $\cF\cap[N]^{\omega}$ will be called  {\em spreading relative to $N$}, if  with $A\in  \cF\cap[N]^{<\omega}$ every spread of $A$ which is a subset of $N$ is in $\cF$.
 Note that if $\cF\subset [N]^{<\omega}$ is spreading relatively to $N=\{n_1,n_2,\ldots \}\in [\N]^{\omega}$, then 
 $$\tilde\cF =\big\{ A\in[\N]^{<\omega}:\{n_j:j\in A\}\in \cF\big\} ,$$
 is spreading.
 
A second way to pass to a sub families is the following. Assume that $\cF\subset[\N]^{<\omega}$ and $N=\{n_1,n_2,\ldots\}\in [\N]^{\omega}$, then we call 
$$\cF^N=\big\{ \{n_j:j\in A\} : A\in\cF\big\},$$
{\em the spread of $\cF$ onto $N$.}

A family $\cF\ksubset[\N]^{<\omega}$ is called 
{\em regular } if it is hereditary, compact, and spreading.
 Note that if  $\cF\ksubset[\N]^{<\omega}$ is compact, spreading and closed under restriction it is also hereditary and thus regular.
  Indeed if $B=\{b_1,b_2,\ldots, b_l\}\in \cF$ and $1\le i_1< i_2<\ldots< i_k\le l $, then
  $A=\{b_{i_1},b_{i_2}, \ldots, b_{i_k}\}$ is a spread of $B'=\{b_1,b_2,\ldots, b_k\}$ and  since
  $B'\in \cF$, it also follows that $A\in\cF$.

If   $\cF\ksubset[N]^{<\omega}$, for some infinite set $N\subset \N$, is called 
  {\em regular relatively to $N$}, if it is compact, spreading relatively  to $N$ and hereditary.  
  
If $\cF\subset [\N]^{<\omega}$ we denote the maximal elements  of $\cF$, \ie the elements $A\in\cF$ for which there is no $B\in\cF$ with $A\prec B$, by $\MAX(\cF)$. Note that if $\cF$ is compact every element in $\cF$ can be extended to a maximal element in $\cF$.

For $\cF\subset[\N]^{<\omega}$ and $A\in [\N]^{<\omega}$ we define 
$$\cF(A)=\{ B\in[\N]^{<\omega}: A< B, A\cup B\in \cF\}.$$
Note that if $\cF$ is compact, spreading, closed under restrictions, or hereditary, so is $\cF(A)$.

If $\cF\subset [\N]^{<\omega}$ is compact, we denote  by $\CB(\cF)$ its {\em Cantor Bendixson} index, which is defined as follows. 
We first define the {\em derivative} of $\cF$ by 
\begin{align*}
\cF'&=\Big\{A\in\cF:  \exists (A_n) \subset \cF\setminus\{A\}  \quad A_n\to_{n\to\infty } A\Big\} \\
      &=\Big\{A\in\cF: \exists (B_n)\subset \cF(A)\setminus\{\emptyset\} \quad B_n\to_{n\to\infty} \emptyset\Big\}
      =\cF\setminus \{ A\in\cF: A\text{ is isolated in }\cF\}.
\end{align*}  
Note that if $\cF$ is hereditary then $\cF'$ is hereditary, but also note that if $\cF$ is only closed under restrictions, this is not necessarily true for 
$\cF'$. Indeed, consider for example
$$\cF=\big\{\emptyset, \{1\},\{1,2\},\{1,2,n\},\, n>2\big\} .$$
Every maximal element $A$ of  $\cF$ is not in $\cF'$ and if $\cF$ is spreading,  then 
$$\cF'=\cF\setminus \MAX(\cF).$$
For $A\in[\N]^{<\omega}$ it follows that 
\begin{equation}\label{E:2.1} \cF'(A)=\big(\cF(A)\big)'.\end{equation}
Indeed 
\begin{align*}
B\in \cF'(A)&\iff  B>A \text{ and } A\cup B\in\cF'\\
&\iff  B>A \text{ and } \exists (C_n)\in \cF\ksetminus\{ A\cup B\} \quad C_n\to_{n\to\infty} A\cup B\\
&\iff B>A \text{ and } \exists (B_n)\in \cF(A) \ksetminus\{ B\} \,\,B_n\to_{n\to\infty}  B\iff B\in\big(\cF(A)\big)'\\
&[\text{Choose $B_n=C_n\setminus A$, for $n$ large enough}].
\end{align*}  
 
 By transfinite induction we define for each ordinal $\alpha$ the {\em $\alpha$-th derivative of $\cF$}, by
 $$\cF^{(0)}=\cF,\, \cF^{(\alpha)}= \big(\cF^{(\gamma)}\big)',\text{ if }\alpha=\gamma+1,
\text{ and }\cF^{(\alpha)}=\bigcap_{\gamma<\alpha}\cF^{(\gamma)}, \text{ if } \alpha \text{ is a limit ordinal.}$$

It follows that $\cF^{(\beta)}\subset \cF^{(\alpha)}$ if $\alpha\le \beta$. By transfinite induction \eqref{E:2.1} generalizes to
\begin{equation}\label{E:2.2} \cF^{(\alpha)}(A)=\big(\cF(A)\big)^{(\alpha)}, 
\text{ for all $A\in[\N]^{<\omega}$ and ordinal  $\alpha$.}
\end{equation}

Assume that $\cF\subset [\N]^{<\omega}$ is compact. Since $\cF$ is countable and since every countable and   compact metric space has isolated points, it follows that for some $\alpha<\omega_1$ 
the $\alpha$-th derivative of $\cF$ is empty and we define
$$\CB(\cF)=\min\{\alpha : \cF^{(\alpha)}=\emptyset\}.$$
$\CB(\cF)$ is always a successor ordinal. Indeed, if $\alpha<\omega_1$ is a limit ordinal and $\cF^{(\gamma)}\not=\emptyset $
for all $\gamma<\alpha$, it follows that $\cF^{(\alpha)}=\bigcap\cF^{(\gamma)}\not=\emptyset$.

\begin{defin}\label{D:2.1} For  $\cF$,$\cG\subset [\N]$ we define  
\begin{align}\label{E:2.1.1}
&\cF\sqcup_< \cG:=\big\{A\cup B: A\in \cF,\,\,B\in \cG\text{ and } A<B \big\}\\
\label{E:2.1.2}&\cF[\cG]:=\left\{ \bigcup_{i=1}^n B_i: \begin{matrix} n\in\N,\, B_1<B_2<\ldots< B_n,\, B_i\in\cG,\, i=1,2,\ldots n,\\
         \text{ and } \{\min(B_i):i=1,2\ldots n\}\in \cF\end{matrix}\right\}.
\end{align}
\end{defin}
It is not hard to see that if $\cF$ and $\cG$ are regular families so are $\cF  \sqcup_<\cG$ and $\cF[\cG]$.

\smallskip

\subsection{The Schreier families}\label{SS:2.2}
We define the {\em Schreier families} $\cS_\alpha\subset[\N]^{<\omega}$ by transfinite induction for all  $\alpha<\omega_1$,
as follows :
\begin{align}\label{E:2.3}
&\cS_0=\big\{ \{n\}: n\in\N\big\}\cup\{\emptyset\}
\intertext{if $\alpha=\gamma+1$, we let }
\label{E:2.4}&\cS_{\alpha}=\cS_1[\cS_\gamma]
    =\Big\{ \bigcup_{j=1}^n E_j: n\kleq\min(E_1),\, E_1< E_2<\ldots< E_n,\, E_j\kin\cS_{\gamma}, j=1,2,\ldots,n \Big\},\\
\intertext{and if $\alpha$ is a limit ordinal we choose a fixed sequence $\big(\lambda(\alpha,n):n\in\N\big)\subset[1,\alpha)$ which increases to $\alpha$ and put}
 \label{E:2.5}  &\cS_\alpha=\{  E: \exists \, k\le \min(E), \text{ with $E\kin\cS_{\lambda(\alpha,k)}$} \}.
 \end{align}

 An easy induction shows that $\cS_\alpha$ is a hereditary, compact and spreading family for all
$\alpha\kle\omega_1$. It is not very hard to see by transfinite induction that $\cS_\alpha$ is in the following very limited sense 
{\em backwards spreading}:
\begin{align}
\label{E:2.6}  &\text{If $A=\{a_1,a_2, \ldots,a_n\}\in\cS_\alpha$ , then }
  \text{$\{a_1,a_2,\ldots,a_{n-1},a_{n}-1\} \in S_\alpha$.}
\end{align} 
So, in particular, if $A\in\cS_\alpha\setminus \{\emptyset\}$ is not maximal, then $(A\cup\{k\})_{k>\max(A)}\subset \cS_\alpha$.

Secondly, by transfinite induction  we can easily  prove that $\cS_\alpha$ is ``almost'' increasing in $\alpha$, in the following sense:
\begin{prop}\label{P:2.2}
For all ordinals $\alpha<\beta<\omega$, there is an $n\in\N$ so that 
$$\cS_\alpha\cap[n,\infty)^{<\omega}\subset \cS_\beta.$$
\end{prop}

The following formula for $\CB(\cS_\alpha)$ is well know and can easily be shown by transfinite induction for all $\alpha<\omega_1$.
\begin{prop}\label{P:2.3} For $\alpha<\omega_1$ we have  $\CB(\cS_\alpha)=\omega^\alpha+1$.
\end{prop}

We now make  further assumptions on the approximating sequence $(\lambda(\alpha,n))\subset [1,\alpha)$, we had chosen to define
the Schreier family $\ S_\alpha$, for limit ordinals $\alpha<\omega_1$. And we will choose  $\big(\lambda(\alpha,n)\big)$  recursively.
Assume that $\alpha$ is a countable limit ordinal and that we have defined $\big(\lambda(\gamma,n)\big)$, for all limit ordinals $\gamma<\alpha$, and thus, $\cS_\gamma$ for all $\gamma<\alpha$.

Recall that  $\alpha$, can be  represented uniquely  in its {\em Cantor Normal Form }
\begin{equation}\label{E:2.7}
\alpha=\omega^{\xi_k} m_k+\omega^{\xi_{k-1}} m_{k-1}+\ldots+\omega^{\xi_1} m_1,
\end{equation}
where $\xi_k>\xi_{k-1}>\ldots \xi_1$, $m_k,m_{k-1},\ldots, m_1\in \N$, and, since $\alpha$ is a limit ordinal, $\xi_1\ge 1$.

We distinguish between three cases:

\noindent Case 1. $k\ge 2$ or $m_1\ge 2$. In that case  we put for $n\in\N$
 \begin{equation}\label{E:2.8}
 \lambda(\alpha,n)=\omega^{\xi_k} m_k+\omega^{\xi_{k-1}} m_{k-1}+\ldots+\omega^{\xi_1} (m_1-1) +\lambda(\omega^{\xi_1},n)
 \end{equation}
 Before considering the next cases we need to make the following observation:

\begin{prop}\label{P:2.4}
Assume that for all limit ordinals  $\gamma\kleq \alpha$  satisfying Case1 the approximating sequences $\big(\lambda(\gamma,n):n\kin\N\big)$  
satisfies the above condition \eqref{E:2.8}.  It follows for all $\gamma\le \alpha$,
with 
$$\gamma=\omega^{\xi_l}m_l+\omega^{\xi_{l-1}}m_{l-1}+ \ldots +\omega^{\xi_1}m_1,$$
being the Cantor Normal Form,  that
\begin{align}
\label{E:2.4.1} &\cS_\gamma= \cS_{\gamma_2}[\cS_{\gamma_1}], \text{ where for some $j=1,\ldots l$}\\
& \gamma_1 =\omega^{\xi_l}m_l+\omega^{\xi_{l-1}}m_{l-1}+ \ldots \omega^{\xi_j}m^{(1)}_j \text{ and }
\gamma_2= \omega^{\xi_j}m_j^{(2)}+\omega^{\xi_{j-1}}m_{j-1}+ \ldots \omega^{\xi_1}m_1,\notag\\
&\qquad \text{with $m^{(1)}_j,m^{(2)}_j\in\N\cup\{0\}$,   $m_j=m^{(1)}_j+m_j^{(2)}$.}\notag
\end{align}
\end{prop}
\begin{proof} We will show  \eqref{E:2.4.1}  by transfinite induction  for  all $\gamma\le \alpha$.
 Assume that  \eqref{E:2.4.1}  holds for all $\tilde \gamma<\gamma$. 
    If $\gamma=\omega^\xi$, then  
\eqref{E:2.4.1} is trivially satisfied. Indeed, in that case $\gamma=\gamma+ 0=0+\gamma$, are the only two choices to write  $\gamma$ as the sum of two ordinals,
 and we observe that $\cS_0[\cS_\gamma]=\cS_\gamma[\cS_0]=\cS_\gamma$.

It is left to verify \eqref{E:2.4.1} in the case that $l\ge 2$ or $m_l\ge 2$.
Let $\gamma=\gamma_1+\gamma_2$ be  a decomposition of $\gamma$ as in the statement of \eqref{E:2.4.1}. We can without loss of generality assume that 
$\gamma_2>0$.

If $\gamma_2=\beta+1$ for some $\beta$ (which implies that $\gamma$ itself is a successor ordinal)
it follows from the  induction hypothesis and  \eqref{E:2.4}  that 
$\cS_{\gamma_1+\beta+1}= \cS_1\big[\cS_\beta[\cS_{\gamma_1}]\big]$, so we need to show that 
\begin{equation*}  \cS_1\big[\cS_\beta[\cS_{\gamma_1}]\big] =\cS_{\beta+1}[\cS_{\gamma_1}].
\end{equation*}
If $A\in  \cS_1\big[\cS_\beta[\cS_{\gamma_1}]\big]$, we can write $A$ as
$A=\bigcup_{i=1}^m A_i$ with 
$m\kleq A_1\kle A_2\kle \ldots\kle A_n$ and $A_i\kin \cS_\beta[\cS_{\gamma_1}]$, for $i\keq1,\ldots, n$ which in turn means
that 
$A_i\keq\bigcup_{j=1}^ {m_i} A_{(i,j)}$, where $A_{(i,1)}\kle A_{(i,2)}\kle \ldots \kle A_{(i,l_i)}$,  $A_{(i,j)}\kin \cS_{\gamma_1}$, for $j\keq1,2,\ldots, l_i$,
and $\{\min(A_{(i,j)}:j=1,2,\ldots,l_i\}\in \cS_\beta$, for $i=1,2,\ldots , m$.
This means that  $\{ \min A_{(i,j)}:j=1,2,\ldots,l_i, i=1,2,\ldots, m\}$ is in $\cS_{\beta+1}$ and thus we conclude that 
$A\in \cS_{\beta+1}[\cS_{\gamma_1}]$. 
Conversely, we can   show in a similar way that $ \cS_{\beta+1}[\cS_{\gamma_1}]\subset  \cS_1\big[\cS_\beta[\cS_{\gamma_1}]\big] $.

If $\gamma_2$ is a limit ordinal we first observe that 
$$\lambda(\gamma, n)= \lambda({\gamma_1+\gamma_2},n )=\gamma_1+ \lambda(\gamma_2,n ).$$
If  $A\in \cS_{\gamma_1+\gamma_2}$ it follows that there is an $n\le \min A$ so that, using the induction hypothesis, we have   
$$A\in  \cS_{\gamma_1+\lambda(\gamma_2,n)}=\cS_{\lambda(\gamma_2,n)}[\cS_{\gamma_1}].$$
This means that $A=\bigcup_{j=1}^m A_j$ with, $A_1<A_2<\ldots< A_m$,
$\{\min(A_j):j=1,2,\ldots, m\}\in \cS_{\lambda(\gamma_2,n)}$ and  $A_j\in\cS_{\gamma_1}$, for $j=1,2,\ldots, m$.
Since $n\le\min(A)=\min(A_1)$,  it follows that $\{\min(A_j):j=1,2,\ldots,m\}\in \cS_{\gamma_2} $, and, thus, that $A\in \cS_{\gamma_2}[\cS_{\gamma_1}]$. 
Conversely, we can similarly show that if $A\in   \cS_{\gamma_2}[\cS_{\gamma_1}]$, then it follows that $A\in\cS_{\gamma_1+\gamma_2}$.
\end{proof}
 If Case 1 does not hold $\alpha$ must be of the form $\alpha=\omega^\gamma$.  
 
\noindent Case 2.  $\alpha=\omega^{\omega^\kappa}$, for some $\kappa<\omega_1$.
In that  cases  we make the following requirement on the sequence $\big(\lambda(\alpha,n):n\in\N\big)$:
\begin{equation}\label{E:2.10}
\cS_{\lambda(\alpha,n)}\subset \cS_{\lambda(\alpha, n+1)}, \text{ for all $n\in\N$}.
\end{equation} 
We can assure \eqref{E:2.10} as follows: first choose any sequence $\lambda'(\alpha,n)$, which increases to  $\alpha$.  Then we notice that Proposition \ref{P:2.2} yields that for a fast enough increasing sequence $(l_n)\subset \N$, it follows that 
$\cS_{\lambda'(\alpha,n)+l_n}\subset \cS_{\lambda'(\alpha,n+1)+l_{n+1}}$. Indeed, we first note that the only set $A\in\cS_\gamma$, $\gamma<\alpha$ which contains $1$, must be the singleton $A=\{1\}$. This follows easily by induction. Secondly we note that by \eqref{E:2.4} it follows 
that $[\{2,3,\ldots,n\}]\subset \cS_{\gamma+n}$, for each $\gamma<\alpha$, and $n\in\N$, which yields our claim.

The remaining case is the following

 \noindent Case 3.   $\alpha=\omega^{\omega^\kappa+\xi}$, where $1\le \xi\le \omega^\kappa$.

We first observe that in this case $\kappa$ and $\xi$ are uniquely defined.

\begin{lem}\label{L:2.5}
Let $\alpha$ be an ordinal number so that there are ordinal numbers $\kappa$, $\xi$ with $\xi \le \omega^\kappa$ and $\alpha = \omega^{\omega^\kappa + \xi}$.
Then for every $\kappa'$, $\xi'$ with $\xi' \le \omega^{\kappa'}$ so that $\alpha = \omega^{\omega^{\kappa'}+\xi'}$, we have $\kappa = \kappa'$ and $\xi = \xi'$.
\end{lem}

\begin{proof}
Let $\alpha = \omega^{\omega^\kappa + \xi} = \omega^{\omega^{\kappa'}+\xi'}$ be as above. By \cite[Theorem 41, \S7.2]{Su} $\omega^\kappa+\xi = \omega^{\kappa'}+\xi'$.
If $\kappa'<\kappa$, then $\omega^{\kappa'}+\xi' \le \omega^{\kappa'}2 < \omega^{\kappa'}\omega = \omega^{\kappa'+1} \le \omega^{\kappa} \le \omega^{\kappa} + \xi$, which is  a contradiction. We conclude that $\kappa\le \kappa'$, and therefore by interchanging the roles of $\kappa$ and $\kappa'$
  we obtain that $\kappa = \kappa'$. In conclusion, $\omega^\kappa + \xi = \omega^\kappa + \xi'$ and therefore $\xi = \xi'$ as well.
\end{proof}
We choose now  a sequence $(\theta(\xi,n))_n$ of ordinal numbers increasing to $\omega^\xi$, so that
\begin{equation}
\label{E:2.11}
\mathcal{S}_{\omega^{\omega^\kappa}\theta(\xi,n)}\subset \mathcal{S}_{\omega^{\omega^\kappa}\theta(\xi,n+1)}
\end{equation}
and define 
\begin{equation}\label{E:2.12} \lambda(\alpha,n) = \omega^{\omega^\kappa}\theta(\xi,n).
\end{equation} 
We  describe how \eqref{E:2.11} can be obtained. Start with an arbitrary sequence $(\theta'(\xi,n))_n$ increasing to $\omega^\xi$. 
We shall recursively choose natural numbers $(k_n)_{n\in\N}$, so that setting $\theta(\xi,n) = \theta'(\xi,n) + k_n$,
 \eqref{E:2.11} is satisfied. Assuming that $k_1,\ldots,k_n$ have been chosen choose $k_{n+1}$,, as in the argument yielding \eqref{E:2.10},  so that
$$\mathcal{S}_{\omega^{\omega^\kappa}\theta(\xi,n)}\subset \mathcal{S}_{\omega^{\omega^\kappa}\theta'(\xi,n+1) + k_{n+1}}.$$
We will show that this $k_{n+1}$ is the desired natural number, i.e. that
$$\mathcal{S}_{\omega^{\omega^\kappa}\theta(\xi,n)}\subset \mathcal{S}_{\omega^{\omega^\kappa}(\theta'(\xi,n+1) + k_{n+1})}.$$
First note that using finite induction and  Proposition \eqref{P:2.4}, it is easy to verify that for  $\gamma<\alpha$, with $\gamma=\omega^\xi$, for some $\xi<\omega_1$, and $n\kin\N$
\begin{equation}
\label{E:2.13}
\cS_{\gamma\cdot n} = \underbrace{\cS_{\gamma}[\cS_{\gamma}\cdots\cS_{\gamma}[\cS_{\gamma}]]]}_{n\text{-times}}.
\end{equation}
and thus
\begin{align*}
\mathcal{S}_{\omega^{\omega^\kappa}(\theta'(\xi,n+1) + k_{n+1})} &=
 =\cS_{\omega^{\omega^\kappa} \Theta'(\xi,n+1)+\omega^{\omega^\kappa} k_{n+1}}=
 \underbrace{\cS_{\omega^{\omega^\kappa}}[\cdots [\cS_{\omega^{\omega^\kappa}}}_{k_{n+1}\text{-times}}[\cS_{\omega^{\omega^\kappa}\theta'(\xi,n+1)}]]] \\
&\supset\underbrace{\cS_{1}[\cdots [\cS_{1}}_{k_{n+1}\text{-times}}[\cS_{\omega^{\omega^\kappa}\theta'(\xi,n+1)}]]]
= \mathcal{S}_{\omega^{\omega^\kappa}\theta'(\xi,n+1) + k_{n+1}} \supset \mathcal{S}_{\omega^{\omega^\kappa}\theta(\xi,n)}.
\end{align*}
We point out that the sequence $(\theta(\xi,n))_n$ also depends on $\alpha$.

 \begin{prop}\label{P:2.6}
Assuming the approximating sequences $\big(\lambda(\alpha,n):n\in\N\big)$ satisfy for all limit ordinals $\alpha$ the above conditions.  It follows for all $\gamma<\omega_1$,
with 
$$\gamma=\omega^{\xi_l}m_l+\omega^{\xi_{l-1}}m_{l-1}+ \ldots +\omega^{\xi_1}m_1,$$
being the Cantor Normal Form,  that
\begin{align}
\label{E:2.6.1}& \cS_{\lambda(\gamma,n)}\subset  \cS_{\lambda(\gamma,n+1)} \text{ for all $n\in\N$, if $\gamma$ is a limit ordinal.}\\
\label{E:2.6.2} &\cS_\gamma= \cS_{\gamma_2}[\cS_{\gamma_1}], \text{ where for some $j=1,2,\ldots ,l$}\\
& \gamma_1 =\omega^{\xi_l}m_l\kplus\omega^{\xi_{l-1}}m_{l-1}\kplus \ldots\kplus\omega^{\xi_j}m^{(1)}_j \text{ and }
\gamma_2= \omega^{\xi_j}m_j^{(2)}\kplus\omega^{\xi_{j-1}}m_{j-1}\kplus \ldots\kplus \omega^{\xi_1}m_1,\notag\\
&\qquad \text{with $m^{(1)}_j,m^{(2)}_j\in\N\cup\{0\}$,   $m_j=m^{(1)}_j+m_j^{(2)}$.}\notag
\intertext{and 
if $\beta = \omega^{\omega^\kappa}$ and $\gamma$ is a limit ordinal with $\gamma \le \beta$, then}
&\label{E:2.6.3}\text{there is a sequence $(\eta(\gamma,n))_n$ increasing to $\gamma$
so that
 $\lambda(\beta\gamma,n) = \beta\eta(\gamma,n)$}
\end{align}
(this sequence $(\eta(\gamma,n))_n$ can depend on $\beta$).
\end{prop}
\begin{proof}
We first will prove \eqref{E:2.6.1}, \eqref{E:2.6.2}  simultaneously for all $\gamma<\omega_1$.
Assume that our claim is true all $\tilde\gamma<\gamma$. 
\eqref{E:2.6.2} follows from Proposition \ref{P:2.4}.

If $l=m_1=1$ we deduce \eqref{E:2.6.1} from the choice of $\lambda(\gamma,n)$, $n\in\N$, in that case (see \eqref{E:2.10}, \eqref{E:2.11} and \eqref{E:2.12}). If  $l\ge 2$ or $m_2\ge 2$ we deduce from \eqref{E:2.4.1}  and the induction hypothesis that 
\begin{align*}
\cS_{\lambda(\gamma,n)}&=\cS_{\omega^{\xi_k}m_k+\ldots+ \omega^{\xi_2}{m_2}+ \omega^{\xi_1}{m_1}+\lambda(\omega^{\xi_1},n)}\\
&=
\cS_{\lambda(\omega^{\xi_1},n)}[\cS_{\omega^{\xi_k}m_k+\ldots+ \omega^{\xi_2}{m_2}+ \omega^{\xi_1}{m_1}}]\\
&\subset
\cS_{\lambda(\omega^{\xi_1},n+1)}[\cS_{\omega^{\xi_k}m_k+\ldots+ \omega^{\xi_2}{m_2}+ \omega^{\xi_1}{m_1}}]=
\cS_{\lambda(\gamma,n+1)},\end{align*}  
which verifies \eqref{E:2.6.1} also in that case

To verify \eqref{E:2.6.3} let $\kappa<\omega_1$ with  $\beta=\omega^{\omega^\kappa}\ge \gamma$.
Recall that by \eqref{E:2.12} $\lambda(\omega^{\omega^\kappa+\xi_1},n) = \omega^{\omega^\kappa}\theta(\xi_1,n)$. For each $n$, define $\eta(\gamma,n) = \omega^{\xi_l}m_l+\omega^{\xi_{l-1}}m_{l-1}+ \cdots +\omega^{\xi_1}(m_1-1) + \theta(\xi_1,n)$.
We will show that $(\eta(\gamma,n))_{n\in\N}$ has the desired property.  Note that the Cantor Normal Form of $\beta\gamma$ is $\beta\gamma = \omega^{\omega^\kappa+\xi_l}m_l+\omega^{\omega^\kappa+\xi_{l-1}}m_{l-1}+ \cdots +\omega^{\omega^\kappa+\xi_1}m_1$.
Hence, by \eqref{E:2.8}:
\begin{eqnarray*}
\lambda(\beta\gamma,n) &=&  \omega^{\omega^\kappa+\xi_l}m_l+\omega^{\omega^\kappa+\xi_{l-1}}m_{l-1}+ \ldots +\omega^{\omega^\kappa+\xi_1}(m_1-1) + \lambda(\omega^{\omega^\kappa+\xi_1},n)\\
&=& \omega^{\omega^\kappa+\xi_l}m_l+\omega^{\omega^\kappa+\xi_{l-1}}m_{l-1}+ \ldots +\omega^{\omega^\kappa+\xi_1}(m_1-1) + \omega^{\omega^\kappa}\theta(\xi_1,n)\\
&=& \omega^{\omega^\kappa}(\omega^{\xi_l}m_l+\omega^{\xi_{l-1}}m_{l-1}+ \ldots +\omega^{\xi_1}(m_1-1) + \theta(\xi_1,n))
 = \beta\eta(\gamma,n).
\end{eqnarray*}
\end{proof}

\begin{rem}
The proof of Proposition \ref{P:2.6}, in particular the definition of $(\eta(\gamma,n))_n$, implies the following.
Let $\xi$ be a countable ordinal number and $\gamma \le \beta = \omega^{\omega^\xi}$ be  a limit ordinal number.
If $\gamma = \omega^{\xi_1}$, then
\begin{equation}
\label{E:2.15}
\eta(\gamma,n) = \theta(\xi_1,n),\text{ for all }n\in\N.
\end{equation}
Otherwise, if the Cantor normal form of $\gamma$ is
$$\gamma=\omega^{\xi_l}m_l+\omega^{\xi_{l-1}}m_{l-1}+ \ldots +\omega^{\xi_1}m_1$$
and $\gamma_1 =\omega^{\xi_l}m_l+\omega^{\xi_{l-1}}m_{l-1}+ \ldots \omega^{\xi_j}m^{(1)}_j$, $\gamma_2= \omega^{\xi_j}m_j^{(2)}+\omega^{\xi_{j-1}}m_{j-1}+ \ldots \omega^{\xi_1}m_1$, with $m^{(1)}_j,m^{(2)}_j\in\N\cup\{0\}$,   $m_j=m^{(1)}_j+m_j^{(2)}$,
then we have
\begin{equation}
\label{E:2.16}
\eta(\gamma,n) = \gamma_1 + \eta(\gamma_2,n),\text{ for all }n\in\N.
\end{equation}

\end{rem}

\begin{cor}\label{C:2.5} If $\alpha<\omega_1$ is a limit ordinal  it follows that 
\begin{equation}\label{E:2.5.1}   \cS_\alpha=\big\{ A\in[\N]^{<\omega}\setminus\{\emptyset\}: A\in \cS_{\lambda(\alpha,\min(A))}\big\}\cup\{\emptyset\}.
\end{equation}
\end{cor} 
\begin{remark}  If we had defined  $\cS_\alpha$ by \eqref{E:2.5.1}, for limit ordinals $\alpha\kle\omega_1$ 
where $\big(\lambda(\alpha,n):n\kin\N)$ is any sequence increasing to $\alpha$, then we would not have ensured that the 
family $\cS_\alpha$ is a regular family.
\end{remark}

\subsection{The  fine Schreier families}\label{SS:2.3}
We will now define the {\em Fine Schreier Sets}.  For that we will also need to choose appropriate approximating sequences for  limit ordinals.
 We will  define them as a doubly indexed family $\cF_{\beta, \alpha}\subset [\N]^{<\omega}$, $\alpha\le \beta<\omega_1$. Later in the proof of Theorems A and C, we will fix $\beta$, depending on the Banach space $X$ we are considering.
\begin{defin}
\label{D:2.9}
For a countable ordinal number $\xi$ and $\beta = \omega^{\omega^{\xi}}$, we recursively define an hierarchy of families of finite subsets of the natural numbers $(\cF_{\beta,\gamma})_{\gamma\le\beta}$ as follows:
\begin{itemize}
 \item[(i)] $\cF_{\beta,0} = \{\emptyset\}$,
 \item[(ii)] if $\gamma<\beta$ then $\cF_{\beta,\gamma+1} = \{\{n\}\cup F: F\in\cF_{\beta,\gamma},n\kin\N\}$ (i.e., $\cF_{\beta,\gamma+1} =\cF_{\beta,1}\sqcup_<\cF_{\beta,\gamma} $), and
 \item[(iii)] if $\gamma\le\beta$ is a limit ordinal number, then $\cF_{\beta,\gamma} = \cup_{n\in\N}(\cF_{\beta,\eta(\gamma,n)}\cap[n,\infty)^{<\omega})$, where  $(\eta(\gamma,n))_n$ is the sequence provided by 
 Proposition \ref{P:2.6}  (and  depends on $\beta$).
\end{itemize}

\end{defin}

\begin{rem}
It can be easily shown by transfinite induction that each family $\cF_{\beta,\gamma}$ is regular.  In the literature fine Schreier families are usually defined  recursively as a singly indexed family $(\cF_\alpha)_{\alpha<\omega_1}$ of subsets of $[\N]^{<\omega}$.
 In that case $\cF_{0}$ and $\cF_\alpha$ are defined  for successor ordinals as in (i) and (ii). And if  $\alpha$ is a limit ordinal is defined as in (iii) where we do not assume that the approximating sequence $(\eta(\alpha,n))_{n\in\N}$  does not depend of any $\beta\ge \alpha$.
\end{rem} 

Let $\xi$ be a countable ordinal number and $\xi_1\le\omega^\xi$. If  and $\beta = \omega^{\omega^\xi}$ and $\gamma = \omega^{\xi_1}$,  it follows by \eqref{E:2.15} that $\eta(\gamma,n)=\theta(\gamma,n)$ for $n\in\N$.
The choice of $(\theta(\xi_1,n))_{n\in\N}$ may be done so that alongside \eqref{E:2.11},
we also have
\begin{equation}
\label{E:2.17}
\cF_{\beta,\eta(\gamma,n)} = \cF_{\beta,\theta(\xi_1,n)} \subset \cF_{\beta,\theta(\xi_1,n+1)} = \cF_{\beta,\eta(\gamma,n+1)}.
\end{equation}
This can be achieved by possibly adding to  $\theta'(\xi_1,n)$  a large enough  natural number.

\begin{prop}\label{P:2.10}
Let $\xi$ be a countable ordinal number and $\beta = \omega^{\omega^\xi}$. Assume that for all limit ordinals limit ordinals $\gamma\le\beta$ the approximating sequence $(\eta(\gamma,n))_n$ satisfies   conditions
\eqref{E:2.15} and \eqref{E:2.16}, and for
the case $\gamma = \omega^{\xi_1}$ the approximating sequence $(\theta(\xi_1,n))_n$ satisfies condition \eqref{E:2.17}.
Then, for all $\gamma\le\beta$, whose Cantor normal form is 
$$\gamma=\omega^{\xi_l}m_l+\omega^{\xi_{l-1}}m_{l-1}+ \cdots +\omega^{\xi_1}m_1$$
we have that that
\begin{equation}
\label{E:2.11.1} \cF_{\beta,\eta(\gamma,n)}\subset  \cF_{\beta,\eta(\gamma,n+1)} \text{ for all $n\in\N$, if $\gamma$ is a limit ordinal,}
\end{equation}
and if for some for some $1\le j\le l$, $\gamma_1 =\omega^{\xi_l}m_l+\omega^{\xi_{l-1}}m_{l-1}+ \cdots +\omega^{\xi_j}m^{(1)}_j$ and $\gamma_2= \omega^{\xi_j}m_j^{(2)}+\omega^{\xi_{j-1}}m_{j-1}+ \cdots +\omega^{\xi_1}m_1$ with
$m^{(1)}_j,m^{(2)}_j\in\N\cup\{0\}$,   $m_j=m^{(1)}_j+m_j^{(2)}$, then
\begin{equation}
\label{E:2.11.2} \cF_{\beta,\gamma} = \cF_{\beta,\gamma_2}\sqcup_<\cF_{\beta,\gamma_1}.
\end{equation}
\end{prop}
%%%% Proof may be omitted in submitted version
\begin{proof} We will show \eqref{E:2.11.1} and \eqref{E:2.11.2} simultaneously by transfinite induction  for  all $\gamma\le\beta$.
 Assume that \eqref{E:2.11.1} and \eqref{E:2.11.2}  hold for all $\tilde \gamma<\gamma$. 
    If $\gamma=\omega^\kappa$, then  \eqref{E:2.11.1} follows from \eqref{E:2.15} and \eqref{E:2.17}, while 
\eqref{E:2.11.2} is trivially satisfied. Indeed, in that case $\gamma=\gamma+ 0=0+\gamma$, are the only two choices to write  $\gamma$ as the sum of two ordinals,
 and we observe that $\cF_{\beta,0}\sqcup_<\cF_{\beta,\gamma} = \cF_{\beta,\gamma}\sqcup_<\cF_{\beta,0} = \cF_{\beta,\gamma}$.

If $\gamma$ is  a limit ordinal (thus $\xi_1>0$) and either $l\ge 2$ or $m_l\ge 2$, then \eqref{E:2.11.1} follows from the inductive assumption. Indeed,
 for $n\in\N$  it follows that 
$\eta(\gamma,n)=\gamma' + \eta(\omega^{\xi_1}, n)$ with 
$\gamma'=\omega^{\xi_l}m_l+\omega^{\xi_{l-1}}m_{l-1}+ \ldots +\omega^{\xi_1}(m_1-1)$ and thus 
(note that $\gamma'+\eta(\omega^{\xi_1},n)<\gamma$)
$$\cF_{\beta,\eta(\gamma,n)}=\cF_{ \beta,\eta(\omega^{\xi_1}, n)}\sqcup_<\cF_{\beta,\gamma'}\subset \cF_{\beta, \eta(\omega^{\xi_1}, n+1)}\sqcup_<\cF_{\beta,\gamma'}=
 \cF_{\beta,\eta(\gamma,n+1)}.$$

It is left to verify \eqref{E:2.11.2} in the case that $l\ge 2$ or $m_l\ge 2$.
Let $\gamma=\gamma_1+\gamma_2$ be  a decomposition of $\gamma$ as in the statement of \eqref{E:2.11.2}. We can without loss of generality assume that 
$\gamma_2>0$.

If $\gamma_2=\gamma_2'+1$ for some $\gamma_2'$ (which implies that $\gamma$ itself is a successor ordinal)
it follows from the  inductive assumption and       Definition \ref{D:2.9} (ii)  that 
\begin{align*}
\cF_{\beta,\gamma}&=\cF_{\beta,\gamma_1+\gamma_2' + 1} = \cF_{\beta,1}\sqcup_<\cF_{\beta,\gamma_1+\gamma_2'} =  \cF_{\beta,1}\sqcup_<\left(\cF_{\beta,\gamma_2'}\sqcup_<\cF_{\beta,\gamma_1}\right)\\
&=\left(\cF_{\beta,1}\sqcup_<\cF_{\beta,\gamma_2'}\right)\sqcup_<\cF_{\beta,\gamma_1} = \cF_{\beta\gamma_2'+1}\sqcup_<\cF_{\beta,\gamma_1} = \cF_{\beta,\gamma_2}\sqcup_<\cF_{\beta,\gamma_1}.
\end{align*}

If $\gamma_2$ is a limit ordinal then recall that by \eqref{E:2.16} we have $\eta(\gamma, n)= \gamma_1 + \eta({\gamma_2},n )$.
If  $A\in \cF_{\beta,\gamma}$ it follows that there is an $n\le\min A$ so that, using the inductive assumption, we have   
$$A\in\cF_{\beta,\eta(\gamma,n)} =  \cF_{\beta,\gamma_1+\eta(\gamma_2,n)}=\cF_{\beta,\eta(\gamma_2,n)}\sqcup_<\cF_{\beta,\gamma_1}.$$
This means that $A=A_1\cup A_2$ with, $A_1<A_2$, $A_1\in\cF_{\beta,\eta(\gamma_2,n)}$, and $A_2\in\cF_{\beta,\gamma_1}$.
If $A_1\neq\emptyset$, then $\min(A_1) = \min (A) \ge n$, i.e. $A_1\in\cF_{\beta,\gamma_2}$ and hence $A\in\cF_{\beta,\gamma_2}\sqcup_<\cF_{\beta,\gamma_1}$.
If on the other hand $A_1 = \emptyset$, then $A\in\cF_{\beta,\gamma_1}\subset \cF_{\beta,\gamma_2}\sqcup_<\cF_{\beta,\gamma_1}$.

Conversely, we can similarly show that if $A\in   \cF_{\beta\gamma_1}\sqcup_<\cF_{\beta,\gamma_1}$, then  $A\in\cF_{\beta,\gamma}$.
\end{proof}

\begin{cor}
\label{C:2.12}
Let $\xi$ be a countable ordinal number and $\gamma\le\beta = \omega^{\omega^\xi}$ be a limit ordinal number. Then
\begin{equation}
\label{E:2.12.1}
\cF_{\beta,\gamma} = \{F\in[\N]^{<\omega}: F\in\cF_{\beta,\eta(\gamma,\min(F))}\}\cup\{\emptyset\}. 
\end{equation}
\end{cor}
The following formula of the Cantor Bendixson index of $\cS_\alpha$ and $\cF_{\beta,\alpha}$ can be easily shown by transfinite induction.
 \begin{prop}\label{P:2.13} For any $\alpha,\kappa<\omega_1$, with  $\alpha\le \beta=\omega^{\omega^\kappa}$,
$$\CB(\cS_\alpha)=\omega^{\alpha}+1 \text{ and } \CB(\cF_{\beta,\alpha})=\alpha+1.$$
Moreover, assuming $\omega^\alpha\le \beta$,   for every $M\in[\N]^{\omega}$, there is an $M\in[N]^{\omega}$  so that
$$\cS_\alpha^N\subset \cF_{\beta, \omega^\alpha}\text{ and }\cF_{\beta, \omega^\alpha}^N\subset \cS_\alpha.$$
\end{prop}

The main result in \cite{Ga} states that if $\cF$ and $\cG$ are two hereditary subsets of $[\N]$, then for any $M\in[\N]^{\omega}$ there is an $N\in[M]^{\omega}$ so that 
 $\cF\cap[N]^{<\omega}\subset \cG$ or $\cG\cap[N]^{<\omega}\subset \cF$. Together with  Proposition \ref{P:2.13} this yields

  \begin{prop}\label{P:2.14} For any $\alpha,\gamma,\kappa<\omega_1$, with  $\omega^\alpha\le \beta=\omega^{\omega^\kappa}$,  and any $M\in[N]^{<\omega}$, there is an $N\in[M]^{<\omega}$ so that 
  \begin{align*}&\cS_\gamma^N\subset \cS_\alpha\cap[N]^{<\omega}\subset\cF_{\beta,\omega^\alpha}, \text{ if $\gamma<\omega^\alpha$, and }
                            \cF_{\beta,\omega^\alpha}^N\subset \cF_{\beta,\omega^\alpha}\cap[N]^{<\omega}\subset \cS_\gamma, \text{ if $\gamma>\omega^\alpha$.}
  \end{align*}
  \end{prop}

\subsection{Families indexed by subsets of $[\N]^{<\omega}$}\label{SS:2.4}
We consider families of the form $(x_A:A\kin\cF)$ in some set $X$ indexed over  $\cF\subset[\N]^{<\omega}$. If $\cF$ is a tree, \ie closed under restrictions, such a family 
 is called  an {\em indexed tree}.   Let us also assume that $\cF$ is spreading. 
  The passing to  a {\em pruning } of such an indexed tree is what corresponds to passing to  subsequences for sequences.
  Formally speaking we define a pruning of $(x_A:A\in\cF)$ as follows.
   Let
  $\pi:\cF\to \cF$ be an order isomorphism with the property that if $F\in\cF$ is not maximal, then for any 
  $n\in \N$, so that $n>\max(A)$ and $A\cup\{n\}\in \cF$,  $\pi(A\cup\{n\})$ is of the form $\pi(A)\cup \{ s_n\}$, where $(s_n)$ is a sequence which  increases with $n$. 
  We call then the family $(x_A:A\in\pi(\cF))$ a {\em pruning } of $(x_A:A\in\cF)$.
  Let
 $\xt_A= x_{\pi(A)}$ for $A\in \cF$. $(\xt_A:A\in \cF)$ is then simply a {\em relabeling } of the family 
 $(x_A: A\in \pi(\cF))$, and we  call it also  a pruning of $(x_A:A\in\cF)$.
 It is important to note that the branches of a pruning of an indexed tree  $(x_A:A\in\cF)$,  are a subset of the  branches 
 of the original tree   $(x_A:A\in\cF)$. Here a branch of $(x_A:A\in\cF)$, is a set of the form 
 $$\xb_F=(x_{\{a_1\}},x_{\{a_1,a_2\}}, \ldots x_{\{a_1,a_2, \ldots,a_l\}}) \text{ for  $F=\{a_1,a_2,\ldots,a_l\}\in \cF$}.$$
 Also the nodes of a pruned tree, namely the sequences of the form $(\xt_{A\cup\{n\}}: A\cup\{n\}\in\cF)$, with $A\kin\cF$ not maximal, are subsequences
  of the nodes of the original tree. 
  
  Let us finally mention, how we usually choose prunings inductively. Let $\{A_n:n\in\N\}$ be a {\em consistent enumeration} of $\cF$. By this we mean 
   that if $\max(A_m)<\max(A_n)$ then $m<n$.
  Thus, we also have 
  if $A_m\prec A_n$, then $m<n$, and if $A_m=A\cup\{s\}\in \cF$ and $A_n=A\cup\{t\}\in\cF$, for some (non maximal) $A\in\cF$ and $s<t$ in $\N$, then $m<n$.
  Of course, then $A_1=\emptyset$ and  $\pi(\emptyset)=\emptyset$ assuming now that $\pi(A_1)$, $\pi(A_2), \ldots ,\pi(A_m)$ have been chosen,
  then $A_{m+1} $ must be of the form $A_m= A_l\subset \{k\} $, with $l<m=1$. Moreover if, $k>\max (A_l)+1$ and if $A_l\cup\{k-1\}\in \cF$ then $A_l\cup\{k-1\}= A_j$ with 
  $l<j<m+1$, and $\pi(A_j)= \pi(A_l)\cup\{s\}$ for some $s$ has already been chosen. Thus, we need to choose $\pi(A_{m+1})$ to be of the form $\pi(A_l)\cup\{t\}$,
  where, in case  that $A_l\cup\{k-1\}\in \cF$, we need to choose $t>s$.

\begin{prop}\label{P:3.1} Assume that $\cF\subset[\N]^{<\omega}$ is compact. Let    $r\in\N$   and 
$f:\MAX(\cF)\to \{1,2,\ldots,  r\}$.  Then for every $M\in [\N]^{\omega}$  there exists an $N\in[M]^{\omega}$ and an $i\in\{1,2,\ldots, r\}$  so that 
$$ \MAX( \cF)\cap [N]^{\omega} \subset\{ A\in\MAX(\cF): f(A)=i\}.$$ 
\end{prop}
Proposition \ref{P:3.1} could be deduced from Corollary 2.5 and Proposition 2.6 in \cite{AKT}. To make the paper as self-contained as possible we want to give a short proof.

\begin{proof}[Proof of \ref{P:3.1}]  Without loss of generality we can assume that $r=2$. We prove our assumption by induction on the Cantor Bendixson index of 
$\cF$. If $\CB(\cF)=1$, then $\cF$ can only consist  of finitely many sets. Since $\cF$ is spreading it follows that $\cF=\{\emptyset\}$ 
and our claim is trivial. Assume that our claim is true  for all regular  families $\cE$, with $\CB(\cE)<\beta$.  Let 
$\cF=\cF_1\cup\cF_2$, where $\cF$ is a regular family with $\CB(\cF)=\alpha+1$, where $\alpha$ is the predecessor of $\beta$, if $\beta$ is a successor ordinal (and thus $\beta=\alpha+1$),
and $\alpha=\beta$, if $\beta$ is a limit ordinal, and $f:\cF\to \{1,2\}$ with
$\cF_1=\{ A\in\cF : f(A)=1\}$ and $\cF_2=\{ A\in\cF : f(A)=2\}$, and let $M\in[\N]^{\omega}$.

First we observe that there is a cofinite subset $M'$ of $M$, so that $\CB(\cF(\{m\}))\le\alpha$, for all $m\in M'$. Indeed, if that were not true,
and thus $\CB(\cF(\{m\}))=\alpha+1$ for all $m$ in an infinite subsets $M'$ of $M$, we could choose for  every $m\in M'$ an
element $A_m\in \big(\cF(\{m \})\big)^{(\alpha)}= \big(\cF^{(\alpha)} \big)(\{m\})$. Thus $\{m\}\cup A_m\in \cF^{(\alpha)}$.
But this would imply that $\emptyset=\lim_{m\in M',m\to\infty}\{m\}\cup A_m\in \cF^{(\alpha+1)}$, which is a contradiction.

Since $\CB(\cF(\{m\}))$ has to be a a successor ordinal we deduce that $\CB(\cF(\{m\}))<\beta$, for $m\in M$, and can therefore apply our induction hypothesis,  and  choose inductively natural numbers $m_1<m_2<m_3<\ldots $ and infinite sets $N_0=M'\supset N_1\supset N_2\ldots $,  and elements $c_1,c_2,\ldots \in\{1,2\}$, 
  so that  $m_j\keq\min N_{j-1}\kle\min(N_j)$ and $\cF(\{m_j\}) \cap [N_j]^{<\omega}\subset \{ A\in\cF: f(A)=c_j\}$, for all $j\in\N$.
  Indeed, if $N_{j-1}$ has been defined we let $m_j=\min(N_{j-1})$, and apply our induction hypothesis 
  to the family $\cF(\{m_j\})$, and the coloring $f_j:\MAX\big(\cF(\{m_j\})\big)\to \{1,2\}$, $B\mapsto f(\{m_j\}\cup B)$, and the set $N_{j-1}\setminus\{m_j\}$,
  to obtain an infinite set $N_j\subset N_{j-1}\setminus\{m_j\}$.

  Then take a $c$ for which $\{j:c_j=c\}$ is infinite and   $N= \{m_j: c_j=c\}$. If $A=\{a_1,a_2,\ldots , a_l\} \in \MAX(\cF)\cap [N]^{<\omega}$, 
  Then $a_1=m_j$ for some $j\in\N$, with $c_j=c$, and $\{a_2,a_3,\ldots, a_n\}\in \MAX(\cF(\{a_1\})\cap [N_j]^{<\omega}$,
  and thus $f(A)=c$, which verifies our claim.
\end{proof}

\section{Repeated averages on  Schreier sets}\label{S:4}

We recall   {\em repeated averages } defined on maximal sets of $S_\alpha$, $\alpha<\omega_1$ (c.f. \cite{AG}). As in our previous sections we will 
assume that $\cS_\alpha$ is recursively defined using the conditions made in Subsection \ref{SS:2.2}.
We first need the following
  characterization of maximal elements of $\cS_\alpha$, $\alpha<\omega_1$, which can be easily proven by transfinite induction using 
  for the limit ordinal case Corollary \ref{C:2.5}.

\begin{prop}\label{P:4.1}
Let $\alpha<\omega_1$ then
\begin{enumerate}
\item  $A\in\MAX(\cS_{\alpha+1})$ if and only if  $A=\bigcup_{j=1}^n A_j$, with $n=\min(A_1)$ and $A_1<A_2<\ldots< A_n$ are in $\MAX(\cS_\alpha)$.
 In this case the $A_j$ are unique.
 \item If $\alpha$ is a limit ordinal then  $A\kin\MAX(\cS_{\alpha})$ if and only if 
  $A\kin \MAX(\cS_{\lambda(\alpha,\min(A))}\!)$.
 \end{enumerate}
\end{prop}

For each $\alpha<\omega_1$ and each 
$A\kin \MAX(\cS_\alpha)$ we will now introduce  an element   $z_{(\alpha,A)}\in S_{\ell_1^+}$ with $\supp(z_{(\alpha,A)})=A$,
which we will call {\em repeated average of complexity $\alpha$ on $A\in \MAX(\cS_\alpha)$}.
If $\alpha=0$ then $\MAX(\cS_0)$ consists of singletons and for $A=\{n\}\in \MAX(\cS_\alpha)$ we put $z_{(0,\{n\})}=e_n$, the  $n$-th element of the unit basis of $\ell_1$. 
Assume for all $\gamma<\alpha$ and all $A\kin\MAX(\cS_\gamma)$ we already introduced
$z_{(\gamma, A)}$ which we write as 
$z_{(\gamma, A)}=\sum_{a\in A} z_{(\gamma, A)}(a) e_a$, with $z_{(\gamma, A)}>0$ for all $a\kin A$. If 
$\alpha=\gamma+1$ for some $\gamma<\omega_1$ and if $A\kin \MAX(\cS_\alpha)$ we  write 
by Proposition \ref{P:4.1} (i)  $A$ in a unique way  as  $A=\bigcup_{j=1}^n A_j$, with $n=\min A$ and $A_1<A_2<\ldots< A_n$ are maximal 
in $\cS_{\gamma}$. We then define
\begin{equation}\label{E:4.1}  z_{(\alpha, A)} =\frac1n \sum_{j=1}^n  z_{(\gamma, A_j)}=\frac1n \sum_{j=1}^n \sum_{a\in A_j} z_{(\gamma, A_j)}(a) e_a,\end{equation}
and thus 
\begin{equation}\label{E:4.2}  
z_{(\alpha, A)}(a)=\frac1n z_{(\gamma, A_j)} (a)\text{ for $j=1,2,\ldots, n$ and $a\in A_j$}.\end{equation}
If $\alpha$ is a limit ordinal and $A\in\MAX( \cS_{\alpha})$ then, by Corollary \ref{C:2.5}, $A\in \cS_{\lambda(\alpha,\min(A))}$, and  we  put
\begin{equation}\label{E:4.3}
z_{(\alpha,A)}=z_{(\lambda(\alpha,\min(A)),A)}=\sum_{a\in A} z_{(\lambda(\alpha,\min(A)),A)}(a) e_a.\end{equation} 
The following result was, with slightly different notation, proved in \cite{AG}.
\begin{lem}\label{L:4.2} {\rm  \cite[Proposition 2.15]{AG}}
For all $\vp>0$, all $\gamma<\alpha$, and all $M\in[\N]^{\omega}$, there is an $N=N(\gamma,\alpha, M,\vp)\in[M]^{\omega}$, so that  
$ \sum_{a\in A}z_{(\alpha,B)}(a)<\vp$
for all 
$B\in \MAX(\cS_\alpha\cap[N]^{\omega})$ and  $A\in\cS_\gamma$.
\end{lem}

The following Proposition will be proved by transfinite induction.
\begin{prop}\label{P:4.3}
Assume $\alpha<\omega_1$ and $A\in \cS_\alpha$ (not necessarily maximal). 
If $B_1,B_2$ are two extensions of $A$ which both are maximal in $\cS_\alpha$ then it follows
$$z_{(\alpha,B_1)}(a)=z_{(\alpha,B_2)}(a)\text{ for all $a\in A$.}$$
\end{prop}
\begin{rem} Proposition \ref{P:4.3} says the following: If  $\alpha<\omega_1$  and $A=\{a_1,a_2,\ldots, a_n\}$ is in $\MAX( \cS_\alpha)$,
then $z_{(\alpha,A)}(a_1)$ only depends on $a_1$, $z_{(\alpha,A)}(a_2)$, only depends on $a_1$ and $a_2$ etc.  
\end{rem}

\begin{proof}[Proof of Proposition \ref{P:4.3}] Our claim is trivial for $\alpha=0$, assume that $\alpha=\gamma+1$ and our claim is true for $\gamma$,
 and let $A\in\cS_{\gamma+1}$. Without loss of generality $A\not=\emptyset$, otherwise we would be done. Using Proposition \ref{P:4.1}, we can find an integer  $1\le l\le \min A$,  sets $A_1,A_2,\ldots, A_{l-1}\in \MAX(S_\gamma)$, 
 and $A_l\in S_\gamma$ (not necessarily maximal in $S_\gamma$) so that  $A_1<A_2<\ldots< A_l$ and
 $A=\bigcup_{j=1}^l A_j.$
 By Proposition \ref{P:4.1}, any extension of $A$ to a maximal element in $S_\gamma$ will then be of the form
 $B=\bigcup_{j=1}^l A_j \cup \bigcup_{j=l}^{\min(A)} B_j$,  
 where $A_l<B_l<B_{l+1}<\ldots<B_{\min(A)}$ ($B_l$ may be empty, in which case $A_l<B_{l+1}<\ldots< B_{\min(A)}$), so that $A_l\cup B_l$, $B_{l+1},\ldots, B_{\min(A)} $ are in $\MAX(\cS_\gamma)$. 
 No matter how we extend $A$ to a maximal element  $B$ in $\cS_{\gamma+1}$, the  inductive formula \eqref{E:4.1} yields 
 $$z_{(\gamma+1, B)}(a)=\frac1{\min(A)} z_{(\gamma,A_j)}(a)\text{ whenever for some $j=1,2,\ldots,l-1$  we have } a\in A_j.$$
 In the case that $a\in A_l$, then, by our induction hypothesis, $z_{\gamma,A_l\cup B_l}(a) $ does not depend on the choice of $B_l$, and 
 $$z_{(\gamma+1,B)}(a) =\frac1n z_{(\gamma,A_l\cup B_l)}(a)\text{ whenever } a\in A_l$$
 Thus, in both cases,   the value of $z_{(\gamma+1,B)}(a) $ does  not depend on how we extend $A$ to a maximal element $B$ in $\cS_{\gamma+1}$.
 
If $\alpha$ is  a limit ordinal and $A\in\cS_\alpha$ is not maximal, we also  can assume that $A\not=\emptyset$, and thus
it follows from the formula  \eqref{E:2.5.1} in Corollary \eqref{C:2.5}   that  $A\in S_{\lambda(\alpha,\min(A))}$. 
For any two extension $B$ of  $A$ into a  maximal set of $\MAX(\cS_\alpha)$, 
 it follows from Proposition \ref{P:4.1} that
 $B$ is  maximal in  $S_{\lambda(\alpha,\min(A))}$, and  that $z_{(\alpha,B)}= z_{(\lambda(\alpha,\min(A),B)}$. Thus, also in this case our claim follows from the induction hypothesis.
\end{proof}

Using Proposition \ref{P:4.3} we can define consistently $z_{(\alpha, A)}\in B_{\ell_1^+}$  for any  $\alpha<\omega_1$ and  any $A\in S_\alpha$
by 
$$z_{(\alpha, A)}=\sum_{a\in A} z_{(\alpha, B)}(a) e_a, \text{ where $B$ is any maximal extension of  $A$ in $\MAX(\cS_\alpha)$}.$$
In particular this implies the following recursive definition of $z_{(\alpha,A)}$:
If $A\in\cS_{\alpha+1}\setminus\{\emptyset\}$ we can  write  in a unique way $A$ as $A=\bigcup_{j=1}^n A_n$, where $A_1<A_2<\ldots< A_n$, $A_j\in \MAX(\cS_\alpha)$, for 
$j=1,2,\ldots, n-1$ and $A_n\in S_\alpha\setminus\{\emptyset\}$,
and note that 
\begin{equation}\label{E:4.4}
 z_{(\alpha+1,A)}=\frac1{\min(A)}\sum_{j=1}^d z_{(\alpha, A_j)}
\end{equation}
and if $\alpha$ is a limit ordinal
 \begin{equation}\label{E:4.5}
 z_{(\alpha,A)}= z_{(\alpha, \lambda(\alpha,\min(A))}.
\end{equation}

For $D\in\cS_\alpha$ define $\zeta(\alpha,A)= z_{(\alpha, D)}(\max(D)$. For $A\in\cS_\alpha$ it follows therefore 
$$z_{(\alpha, D)}=\sum_{D\preceq A} \zeta(\alpha, D)e_{\max(D)}.$$
We also 
 put $\zeta(\alpha,\emptyset)=0$ and $e_{\max(\emptyset)}=0$.

By transfinite induction we can easily show the following estimate fro $1\le \alpha<\omega_1$
\begin{equation}
\label{E:11.5}
\zeta(\alpha,A)\le \frac{1}{\min A}.
\end{equation}
From Proposition \ref{P:2.6} we deduce the following formula for $z_{(\alpha,A)}$
\begin{prop}\label{P:4.4} Assume $\alpha<\omega_1$ and that its Cantor Normal Form is 
$$\alpha=\omega^{\xi_l} m_l + \omega^{\xi_{l-1}} m_{l-1}+\ldots +\omega^{\xi_1} m_1.$$
Let  $j=1,2\ldots l$ and $m^{(1)}_j,m^{(2)}_j\in\N\cup\{0\}$, with 
 $m^{(1)}_j+m^{(2)}_j=m_j$. Put
 \begin{align*}
 \gamma_1&= \omega^{\xi_l} m_l + \omega^{\xi_{l-1}} m_{l-1}+\ldots +\omega^{\xi_{j-1}} m_{j-1}+\omega^{\xi_{j}} m_{j}^{(1)}\\
 \gamma_2&=\omega^{\xi_j} m_{j}^{(2)}+\omega^{\xi_{j-1}} m_{j-1}+\ldots+\omega^{\xi_1} m_1.
 \end{align*}
 For $A\in\MAX(\cS_\alpha)$ we use Proposition \ref{P:2.6} and write  
 $A=\bigcup_{j=1}^n A_j$, where 
 $A_j\in\cS_{\gamma_1}$ for $j=1,2,\ldots, n$, $A_1<A_2<\ldots< A_n$, and
 $B=\{ \min(A_j): j=1,2,\ldots, n\}\in \cS_{\gamma_2}$.
 
 Then it follows that $A_j\in \MAX(\cS_{\gamma_1})$, for $j=1,2,\ldots, n$, 
 $B\in\MAX(\cS_{\gamma_2})$ and
 \begin{equation}\label{E:4.4.1}
  z_{(\alpha,A)}=\sum_{j=1}^n z_{(\gamma_2,B)}(\min(A_j)) z_{(\gamma_1, A_j)}.
 \end{equation}
 In other words, if $\emptyset \prec D\preceq A$, and thus $D=\bigcup_{j=1}^{i-1} A_j\cup \At_i$,
 for some $0\le i<n$, and $\emptyset\prec \At_i\preceq A_i$, then
 \begin{equation}\label{E:4.4.1a}
 \zeta(\alpha,D)= \zeta(\gamma_2,\{\min(A_j): j=1,2,\ldots i\})\cdot \zeta(\gamma_1,\At_i).
  \end{equation}
\end{prop}

\begin{proof}
We prove by transfinite induction for all $\beta<\omega_1$,
  with  Cantor Normal Form
  $$\beta=\omega^{\xi_j} m_{j}+\omega^{\xi_{j-1}} m_{j-1}+\ldots+\omega^{\xi_1} m_1,$$
  the following
  
  \noindent{\bf Claim:} If  $\gamma<\omega_1$ has  Cantor Normal Form 
  $$ \gamma=\omega^{\xi_l} \mt_{l}+\omega^{\xi_{l-1}}\mt_{l-1}+\ldots+\omega^{\xi_j} \mt_j,$$
 where  $\mt_j$ could possibly be vanishing, 
    and if  $A=\bigcup_{i=1}^n A_i\in\MAX(\cS_{\gamma+\beta})=\MAX\big(\cS_{\beta}[\cS_{\gamma}]\big)$, where 
 $A_i\in\cS_{\gamma}$ for $i=1,2,\ldots, n$, $A_1<A_2<\ldots <A_n$, and
 $B=\{ \min(A_i): i=1,2,\ldots, n\}\in \cS_{\beta}$, then 
   it follows that $A_i\in \MAX(\cS_{\gamma})$ for $i=1,2,\ldots, n$, 
 $B\in\MAX(\cS_{\gamma_1})$ and 
 \begin{equation}\label{E:4.4.2}
  z_{(\alpha,A)}=\sum_{i=1}^n z_{(\beta,B)}(\min(A_i)) z_{(\gamma, A_i)}
 \end{equation}

For $\beta=0$ the claim is trivial and for  $\beta=1$, our claim follows from  Proposition  \ref{P:4.1} and the definition of $z_{(\gamma+1, A)}$ for 
$A\in  \MAX(\cS_{\gamma+1})$.

 Assume now that the  claim  is true for all $\tilde\beta<\beta$, and that $\gamma<\omega_1$ has above form and let 
  $A=\bigcup_{i=1}^n A_i\in\MAX(\cS_{\gamma+\beta})$, where 
 $A_i\in\cS_{\gamma}$ for $i=1,2,\ldots, n$, $A_1<A_2<\ldots <A_n$, and
 $B=\{ \min(A_i): i=1,2,\ldots, n\}\in \cS_{\beta}$.

First we note that \eqref{E:2.6} implies that the $A_i$ are maximal  in $\cS_\gamma$. Indeed,
if for some $i_0=1,2,\ldots, n$, $A_{i_0}$ is not  maximal in $\cS_\gamma$, then if $i_0=n$, it would directly follow
that $A$ cannot be maximal in $\cS_{\gamma+\beta}$ and if $i_0<n$
we could define
$\tilde A_i=A_i$, for $i=1,2,\ldots,, l-1$, 
$\tilde A_{i_0}= A_{i_0}\cup\{ \min(A_{i_0+1})\}$, 
$\tilde A_{i}= (A_{i}\cup \{\min (A_{i+1})\} )\setminus \{ \min(A_{i})\}$,
for $i=i_0,i_0+1,\ldots,l-1$,
and $\tilde A_{l}=A_l\setminus\{\min A_{l}\}$.
Then, by  \eqref{E:2.6} and the fact that the Schreier families are spreading,  $A=\bigcup_{i=1}^n\tilde A_i$ is also a decomposition of elements of $\cS_\gamma$ with $\tilde B=\{ \min(\tilde A_i):i=1,2,\ldots,n\}\in\cS_\beta$. But now $\tilde A_n$ is not maximal in 
$\cS_\gamma$ and we get again a contradiction.

It is also easy to see that $B$ is maximal in $\cS_\beta$.

To verify \eqref{E:4.4.2} we first assume that $\beta$ is a successor ordinal, say $\beta=\alpha+1$.
Then we can write $B$ as $B=\bigcup_{i=1}^m B_i$, where $m=\min(B) =\min(A)$, 
$B_1<B_2<\ldots <B_m$, and $B_i\in\MAX(\cS_\alpha)$, for $i=1,2,\ldots, m$.
We can write $B_i$ as 
$$B_i=\{ \min(A_s): s=k_{i-1}+1, k_{i-1}+2,\ldots, k_i\},$$
with $k_0=0<k_1<\ldots< k_m=n$. We put
 $\Ab_i=\bigcup_{s=k_{i-1}+1}^{k_i} A_s\in \cS_{\gamma+\alpha}=\cS_\alpha[\cS_\gamma]$,
for $i=1,2, ,\ldots, m$.
From the definition of $z_{(\beta+1,B)}$ and from the induction hypothesis  we deduce now that 
\begin{align*}
z_{(\gamma+\alpha+1,A )}&= \frac1m \sum_{i=1}^m  z_{(\gamma+\alpha, \Ab_i)}\\
&= \frac1m \sum_{i=1}^m  \sum_{s=k_{i-1}+1}^{k_i}\!\!\!z_{(\alpha,B_i)}(\min(A_s)) z_{(\gamma, A_s)}
=\sum_{s=1}^n z_{(\beta,B)}(\min(A_s))z_{(\gamma, A_s)},
\end{align*}
which proves the claim if $\beta$ is a successor ordinal.

If $\beta$ is a limit ordinal it follows from Corollary \ref{C:2.5}, our definition  of
 $z_{(\beta,B)}$ and $z_{(\gamma+\beta,A)}$,
 and our choice of the approximating sequence $(\lambda(\gamma+\beta),n)$ that 
\begin{align*}
z_{(\gamma+\beta,A)}&=
z_{(\lambda(\gamma+\beta,\min(A)),A)}
=z_{(\gamma+\lambda(\beta,\min(B)),A)}\\
&=\sum_{j=1}^n z_{(\lambda(\beta,\min(B)),B)}(\min(A_j))z_{(\gamma,A_j)}=
\sum_{j=1}^n \!\!z_{(\beta,B)}(\min(A_j))z_{(\gamma,A_j)}
\end{align*}
which proves our claim also in the limit ordinal case.
\end{proof}

If $\alpha<\omega_1$ and $A\in \MAX(\cS_\alpha)$ then  $z_{(\alpha,A)}$ is an element of $S_{\ell_1}\cap \ell_1^+$ and can  therefore be seen as probability on $\N$, as well as $A$.
We denote the expected value of a function $f$ defined on $A$ or on all of $\N$ as $\E_{(\alpha,A)}(f)$.

As done in \cite{Sch}, we deduce the following statement from Lemma \ref{L:4.2}.
\begin{cor}\label{C:4.5}{\rm\cite[Corollary 4.10]{Sch} } For each  $\alpha<\omega_1$ and  $A\in\MAX( \cS_\alpha)$  let $f_A:A\to [-1,1]$  have the property that 
$\E_{\alpha,A} (f_A)\ge \rho$,  for some fixed number  $\rho\in [-1,1]$.
For  $\delta>0$ and $M\in[\N]^{\omega}$ put 
$$\cA_{\delta,M}= \left\{ A\in \cS_\alpha\cap[M]^{<\infty}:\begin{matrix} \exists B\kin \MAX(S_\alpha\cap[M]^{<\infty}), \\ A\ksubset B,  \text{ and } f_B(a)\ge \rho-\delta \text{ for all } a\kin A\end{matrix}\right\}.$$
Then $\CB(\cA_{\delta,M})=\omega^\alpha+1$.  
\end{cor}
We finish this section with an observation, which will be needed later.

 \begin{defin}\label{D:4.5a}
If $A\subset \N\setminus\{\emptyset\}$ is finite, we can write it in a unique way as  a union $A=\bigcup_{j=1}^d A_j$, where  $A_1<A_2<\ldots< A_d$,
and $A_j\in\MAX(\cS_1)$ if $j=1,2,\ldots,d-1$, and $A_d\in\cS_1\setminus\{\emptyset\}$.  We call $(A_j)_{j=1}^d$ the {\em optimal $\cS_1$-decomposition of $A$},
and we define  
$$l_1(A)=\min(A_{d})-   \#A_{d}.$$
 For $A=\emptyset $ we put $l_1(A)=0$.
 \end{defin}
The significance of this number and its connection to  the repeated averages is explained in the following Lemma.
\begin{lem}\label{L:4.6} Let $\alpha\in[1,\omega_1)$, $A\in\cS_\alpha$ and let $(A_j)_{j=1}^d$ be its optimal $\cS_1$-decomposition.
\begin{enumerate}
\item $l_1(A)=0$ if and only if  $A=\emptyset$ or $A_d\in\MAX(\cS_1)$.
\item If $A\in \MAX(\cS_\alpha)$ then $A_d\in\MAX(\cS_1)$ and, thus, $l_1(A)=0$. 
\item If $l_1(A)>0$, then for all $\max(A)<k_1<k_2<\ldots <k_{l_1(A)}$ it follows that 
$A\cup\{k_1,k_2,\ldots,k_{l_1(A)}\}\in \cS_\alpha$ and 
$$\zeta(\alpha,A\cup\{k_1,k_2,\ldots, k_i\})=\zeta(\alpha,A)\text{ for all $i=1,2\ldots l_1(A)$}.$$
\item If $m> l_1(A)$ and  $\max(A)<k_1<k_2<\ldots <k_m$, have the property that 
$A\cup\{k_1,k_2,\ldots k_m\}\in\cS_\alpha$ then
$$\zeta(\alpha, A\cup\{k_1,k_2,\ldots k_i\})\le \frac1{k_{l_1(A)+1}}.$$
\item If $A\not= \emptyset$, then 
 $$ \sum_{D\preceq A, l_1(D')=0} \zeta(\alpha, D)\le \frac1{\min(A)} \text{ and }  \sum_{D\preceq A, l_1(D)=0} \zeta(\alpha, D)\le \frac1{\min(A)} $$
\end{enumerate} 
(recall that $D'=D\setminus\{\max D\}$ for $D\in[\N]^{<\omega}\setminus\{\emptyset\}$ and $\emptyset'=\emptyset$).
\end{lem}

\begin{proof} We proof (i) through (v) by transfinite induction for all $\alpha\in[1,\omega_1)$. For $\alpha=1$, (i), (ii), (iii)  and (v) follow from the definition of $\cS_1$ and the definition
of $\zeta(\alpha,A)$, for $A\in \cS_1$, while (iv) is vacuous in that case.
Assume our claim is true for some  $\alpha<\omega_1$, and let $A\in \cS_{\alpha+1}$.  Without loss of generality we can assume that $A\not=\emptyset$. Indeed,
if $A=\emptyset$, then  (i) is clear (ii), (iii), and (v) are vacuous, while (iv) follows easily by induction
from the fact that always $\zeta(\alpha,A)\le\frac1{\min(A)}$ if $A\in\cS_\alpha\setminus\{\emptyset\}$. By the definition of $\cS_{\alpha+1}$,
$A$ can be written in a unique way as $A=\bigcup_{j=1}^n B_j$ where $B_j\in\MAX(\cS_\alpha)$, for $j=1,2,\ldots, n-1$, and $B_n\in\cS_\alpha$.
For $j=1,2,\ldots,n$ let $(A_{j,i})_{i=1}^{c_j}$ be the optimal $\cS_1$-decomposition of $B_j$.
From  the induction hypothesis (ii) it follows that $(A_{j,i})$ are maximal in $\cS_1$, for $j<n$ or for $j=n$ and $i<c_n$.
Therefore it follows that $(A_{j,i}:j=1,2,\ldots,n, i=1,2\ldots c_j)$ (appropriately ordered) is the optimal $\cS_1$-decomposition  of $A$ and it follows   $l_1(A)=l_1(B_n)$,
and $A_d=A_{n,c_n}$.

We can deduce (i) from the induction hypothesis. If $A\in\MAX(S_{\alpha+1})$,  then, in particular, $B_n\in \MAX(S_\alpha)$, and thus $l_1(A)=l_1(B_n)=0$.
Conversely, if $l_1(A)=l_1(B_n)=0$, then $A_d=A_{n,c_n}\in \MAX(S_1)$. This proves (ii) for $\alpha+1$.

If $l_1(A)\kgr 0$ and  $\max(A)\keq\max(B_n)\kle k_1\kle k_2\kle\ldots\kle k_{l_1(A)}$, then it follows from the fact that $l_1(A)\keq l_1(B_n)$ and  our induction hypothesis that 
$B_n\cup\{k_1,k_2,\ldots, k_{l_1(A)}\}\kin \cS_\alpha$ and 
$$\zeta(\alpha, B_n\cup\{k_1,k_2,\ldots, k_i\})= \zeta(\alpha,B_n).$$
Therefore $A\cup  \{k_1,k_2,\ldots, k_{l_1(A)}\}\in\cS_{\alpha+1}$ and, using our recursive formula, we obtain
$$\zeta(\alpha, A\cup\{k_1,k_2,\ldots, k_i\})= \frac1{\min(A)} \zeta(\alpha, B_n\cup\{k_1,k_2,\ldots, k_i\})=\frac1{\min(A)} \zeta(\alpha, B_n)=\zeta(\alpha,A),$$
which verifies (iii).
In order to show (iv) let  $m>l_1(A)$ and $\max(A)<k_1<k_2<\ldots< k_m$, is such that $A\cup\{k_1,k_2,\ldots,k_m\} \in \cS_{\alpha+1}$ we distinguish between two cases. Either
$B_n\cup\{k_1,k_2,\ldots,k_m\} \in \cS_\alpha$. In that case we deduce from the induction hypothesis 
$$\zeta\big(\alpha, A\cup\{k_1,k_2,\ldots, k_m\}\big)= \frac1{\min(A)} \zeta(\alpha, B_n\cup\{k_1,k_2,\ldots, k_m\})\le \frac1{k_{l_1(A)+1}}.$$
Or we can write $A\cup\{k_1,k_2,\ldots,k_m\}$ as
$A\cup\{k_1,k_2,\ldots,k_m\}= \bigcup_{j=1}^n B_j \cup\bigcup_{j=n}^p B'_j,$
where $p>n$,  $B_n< B_n'<B'_{n+1}<\ldots<   B'_p$, 
$B_n\cup B_n'\in \MAX(\cS_\alpha )$, $B'_{n+1}, \ldots, B'_{p-1}\in \MAX(\cS_\alpha)$, and $B'_p\in\cS_\alpha\setminus\{\emptyset\}$.
Let $s\le m$ such that $k_s=\min(B_p)$. Then  $s>l_1(B_n)$, and  $l_1(B_n\cup B_n')=0$. It follows therefore from 
 \eqref{E:4.4} and the induction hypothesis  that
 $$\zeta(\alpha+1,A\cup\{k_1,k_2,\ldots,k_m\})=\frac1{\min(A)}\zeta(\alpha, \{k_s,k_{s+1},\ldots,k_{m}\} )\le \frac1{k_s}\le \frac1{k_{l_1(A)+1}},$$
This proves (iv) in both cases.

Finally, to verify (v) we   observe, that  by the induction hypothesis and \eqref{E:4.4}
\begin{align*}
 \sum_{D\preceq A, l_1(D')=0} \zeta(\alpha+1, D) &=\frac1{\min( A)} \sum_{j=1}^n\sum_{D\preceq B_j, l_1(D')=0} \zeta(\alpha, D)\\
&\le  \frac1{\min A} \sum_{j=1}^n \frac1{\min(B_j)}\le \frac{n}{\min(A)} \frac1{\min(A)}\le \frac1{\min(A)},
\end{align*}
which proves the first part of (v), while the second follows in the same way.

If   $\alpha<\omega_1$ is a limit ordinal and assuming that our claim is true for all $\gamma<\alpha$ we proceed as follows.
For $A\in\cS_\alpha$, we a can assume again that $A\not=\emptyset$ and  it follows from Corollary \eqref{C:2.5} that $A\in \cS_{\lambda(\alpha, \min(A))}$ and, by Proposition \ref{P:4.1} $A$ is maximal
in $\cS_\alpha$ if and only if it is maximal in $\cS_{\lambda(\alpha, \min(A))}$. Therefore (i) through (v) follow from our claim being true for $\lambda(\alpha,\min(A))$.  
\end{proof}

\begin{rem}
Recall, that if $\beta = \omega^{\omega^\xi}$ is a countable ordinal number and $\gamma<\beta$, then by \eqref{E:2.13} we have $\cS_{\beta(\gamma+1)} = \cS_\beta[\cS_{\beta\gamma}]$.
An argument very similar to what was used in the proof of Lemma \ref{L:4.6}, implies the following: if $B_1<\cdots<B_d$ are in $\MAX(\cS_{\beta\gamma})$ so that $\bar B = \{\min(B_j): 1\le j\le d\}$ is
a non-maximal $\cS_\beta$ set, $D = \cup_{j=1}^dB_j$ and $C\in\cS_{\beta\gamma}$ with $D<C$, then
\begin{equation}
\label{E:4.6}
l_1( C ) = l_{1}(D\cup C).
\end{equation}
\end{rem}

\begin{cor}
\label{C:4.7}
Let $A=\{a_1,a_2,\ldots, a_l\}$ and $\tilde A =\{\at_1,\ldots, \at_{\lt}\}$  be two sets in $[\N]^{<\omega}$ whose optimal $\cS_1$-decompositions  $(A_j)_{j=1}^d$ and $(\tilde A_j)_{j=1}^d$, respectively, have the same length 
and satisfy $\min(A_j)=\min(\tilde A_j)$, for $j=1,2,\ldots, d$.

Then it follows for $\alpha<\omega_1$ that $A\in\cS_\alpha$ if and only if $\At\in \cS_\alpha$ and in that case  for $D\preceq A$ and $\Dt\preceq \At$, with $\#D=\#\Dt$, it follows that 
$\zeta(\alpha,D)=\zeta(\alpha, \Dt)$.
\end{cor}

\begin{proof}
We prove this lemma by transfinite induction on $\alpha$. If $\alpha = 1$ then, $A_1=A$ and $\At_1=\At$, and $a_1=\at_1$, and thus
$\zeta(1,D)=\zeta(1, \Dt)$ for all $D\preceq A$ and $\Dt\prec \At$.

Assume that
the conclusion holds for some $\alpha$ and let $A\in\cS_{\alpha+1}$, $\tilde A\in[\N]^{<\omega}$ satisfy the assumption.
 Let $A = \cup_{i=1}^pC_i$,
where $C_1<\cdots<C_{p-1}$ are in $\MAX(\cS_\alpha)$, whereas $C_p\in\cS_\alpha$ and $p\le \min(A)$.
 Write also $\At$ as $\At = \cup_{i=1}^p \Ct_i$,
where $\Ct_1<\cdots<\Ct_p$, and the $\Ct_j $ are chosen such that $\#C_j=\#\Ct_j$ for $j=1,2,\ldots, p-1$.

From Lemma \ref{L:4.6} (ii) It follows that for some sequence $0=d_0<d_1<d_2<\ldots< d_p=d$
the sequence $(A_j)_{j=d_{i-1}+1}^{d_i}$ is the optimal $\cS_1$-decomposition of $C_i$ for $i=1,2,\ldots,p$.
Now we can first deduce 
 $(\At_j)_{j=1}^{d_1}$ is the  optimal $\cS_1$-decomposition of $\Ct_1$, then 
 deduce that   $(\At_j)_{j=d_1+1}^{d_2}$is the  optimal $\cS_1$-decomposition of $\Ct_2$, and so on.
 We are therefore in the position to apply the induction hypothesis and deduce that for all
 $i=1,2,\ldots, p$, and $D\preceq C_i$ and $\Dt\preceq \Ct_i$ it follows that 
 $\zeta(\alpha, D)= \zeta(\alpha,\Dt)$. Our claim follows therefore from our recursive formula \eqref{E:4.4}.
 
 As usual in the case that $\alpha $ is a limit ordinal the verification follows easily from the definition of $\cS_\alpha$ in that case.
\end{proof}

\begin{lem}\label{L:4.8}
Let $X$ be a Banach space, $\alpha$ be a countable ordinal number, $B\in\MAX(\cS_\alpha)$ and $(x_A)_{A\preceq B}$ be vectors in $B_X$. Then
\begin{equation}
\label{E:4.8.1}
\Big\|\sum_{A\preceq B}\zeta(\alpha,A)x_A - \sum_{A\preceq B}\zeta(\alpha,A')x_A\Big\| \le \frac{2}{\min(B)}.
\end{equation}\end{lem}

\begin{proof} Using, Lemma \ref{L:4.6}  part (iii), then  part (v) we obtain
\begin{subequations}
\begin{eqnarray}
\Big\|\sum_{A\preceq B}\!\!&&\zeta(\alpha,A)x_A - \sum_{A\preceq B}\zeta(\alpha,A')x_A\Big\|\nonumber\\
&\le& \Big\|\sum_{\substack{A\preceq B\\ l_1(A')\not=0}}\left(\zeta(\alpha,A) - \zeta(\alpha,A')\right)x_A\Big\| 
 \kplus \Big\|\!\!\sum_{\substack{A\preceq B\\ l_1(A')=0}}\zeta(\alpha,A)x_A\Big\| \kplus \Big\|\!\!\sum_{\substack{A\preceq B\\ l_1(A')=0}}\zeta(\alpha,A')x_A\Big\|\nonumber\\
&{\le} & \sum_{\substack{A\preceq B\\ l_1(A')=0}}\zeta(\alpha,A) +\sum_{\substack{A\preceq B\\ l_1(A')=0}}\zeta(\alpha,A') \le  \frac{2}{\min(B)}.\nonumber
\end{eqnarray}
\end{subequations}
\end{proof}

\section{Trees and their indices}\label{S:5}

Let $X$ be an arbitrary set. We set $X^{<\omega}\keq
\bigcup_{n=0}^\infty X^n$, the set of all finite sequences in $X$,
which includes the sequence of length zero denoted by $\emptyset$. For
$x\kin X$ we shall write $x$ instead of $(x)$, \ie we identify $X$
with sequences of length~$1$ in $X$. A \emph{tree on $X$ }is a 
non-empty
subset $\cF$ of $X^{<\omega}$ closed under taking initial segments: if
$(x_1,\dots,x_n)\kin \cF$ and $0\kleq m\kleq n$, then
$(x_1,\dots,x_m)\kin\cF$. A tree $\cF$ on $X$ is \emph{hereditary }if
every subsequence of every member of $\cF$ is also in $\cF$.

Given $\xb\keq (x_1,\dots,x_m)$ and $\yb\keq(y_1,\dots,y_n)$ in
$X^{<\omega}$, we write $(\xb,\yb)$ for the concatenation of $\xb$
and $\yb$:
$$
(\xb,\yb)=(x_1,\dots,x_m,y_1,\dots,y_n)\ .
$$
Given $\cF\ksubset X^{<\omega}$ and $\xb\kin X^{<\omega}$, we let
$$
\cF(\xb)= \{\yb\kin X^{<\omega}:\,(\xb,\yb)\kin\cF\}\ .
$$
Note that if $\cF$ is a tree on $X$, then so is $\cF(\xb)$
(unless it is empty). Moreover,
if $\cF$ is hereditary, then so is $\cF(\xb)$ and $\cF(\xb)\ksubset
\cF$.

Let $X^\omega$ denote the set of all (infinite) sequences in $X$. Fix
$S\ksubset X^\omega$. For 
 a subset 
 $\cF$ of $X^{<\omega}$ \emph{the $S$-derivative
  $\cF_S'$ of $\cF$ }consists of all  $\xb=(x_1,x_2, \ldots,x_l)\kin
X^{<\omega}$ for which there is a sequence $(y_i)_{i=1}^\infty\kin S$
with $(\xb,y_i)\kin\cF$ for all $i\kin\N$.
 
 Note that if $\cF$ is a hereditary tree then it follows that
$\cF'_S\ksubset\cF$ and that $\cF'_S$ is also a hereditary  tree
(unless it is empty).

We then define
higher order derivatives $\cF^{(\alpha)}_S$ for ordinals
$\alpha\kle\omega_1$ by recursion as follows.

\begin{align*}
  \cF^{(0)}_S = \cF, 
  \cF_S^{(\alpha+1)} = \big( \cF_S^{(\alpha)}\big)'_S\text{ for }\alpha\kle\omega_1\text{ and }
  \cF_S^{(\lambda)} = \bigcap _{\alpha <\lambda}
  \cF_S^{(\alpha)} \text{ for limit ordinals
  }\lambda\kle\omega_1.
\end{align*}

It is clear that $\cF^{(\alpha)}_S\ksupset \cF_S^{(\beta)}$  if 
$\alpha\kleq\beta$ and that $\cF_S^{(\alpha)}$ is a hereditary tree
(or the empty set) 
 for all
$\alpha$, whenever $\cF$ is a hereditary tree. An easy induction also shows that
$$
\big(\cF(\xb)\big)_S^{(\alpha)} = \big(\cF_S^{(\alpha)}\big)
(\xb)\qquad \text{for all }\xb\kin X^{<\omega},\ \alpha\kle\omega_1\ .
$$

We now define \emph{the $S$-index $I_S(\cF)$ of $\cF$ }by
$$
I_S(\cF) = \min \big\{\alpha\kle\omega_1:\,\cF_S^{(\alpha)}\keq \emptyset
\big\}
$$
if there exists $\alpha\kle\omega_1$ with $\cF_S^{(\alpha)}\keq
\emptyset$, and $I_S(\cF)\keq\omega$ otherwise.

\begin{rem}
  If $\lambda$ is a limit ordinal and
  $\cF_S^{(\alpha)}\kneq\emptyset$ for all $\alpha\kle\lambda$, then
  in particular $\emptyset\kin\cF_S^{(\alpha)}$ for all
  $\alpha\kle\lambda$, and hence
  $\cF_S^{(\lambda)}\kneq\emptyset$. This shows that $I_S(\cF)$ is
  always a successor ordinal.
\end{rem}

\begin{exs}\label{Ex:5.1}
  1.~A  family $\cF\ksubset [\N]^{<\omega}$ can be thought of
  as a tree on $\N$: a set $F\keq\{m_1,\dots,m_k\}\kin[\N]^{<\omega}$ is
  identified with $(m_1,\dots,m_k)\kin\N^{<\omega}$ (recall that
  $m_1\kle\dots\kle m_k$ by our convention of always listing the
  elements of a subset of $\N$ in increasing order).

  Let $S$ be the set of all strictly increasing sequences in $\N$. In
  this case the $S$-index of a compact,  family
  $\cF\ksubset[\N]^{<\omega}$ is nothing else but the Cantor-Bendixson
  index of $\cF$ as a compact topological space, which we will continue to denote
  by $\CB(\cF)$. We will also use the term Cantor-Bendixson
  derivative instead of $S$-derivative and use the notation
  $\cF'$ and $\cF^{(\alpha)}$.
   
    \noindent
  2.~If $X$ is an arbitrary set and $S\keq X^\omega$, then the
  $S$-index of a tree $\cF$ on $X$ is what is usually called \emph{the
  order of $\cF$ }(or \emph{the height of $\cF$}) denoted by
  o$(\cF)$. Note that in this case the $S$-derivative of $\cF$
  consists of all finite sequences $\xb\kin X^{<\omega}$ for which
  there exists $y\kin X$ such that $(\xb,y)\kin\cF$.

  The function o$(\cdot)$ is the largest index: for any
  $S\ksubset X^\omega$ we have o$(\cF)\kgeq I_S(\cF)$.
\end{exs}

We say that $S\ksubset X^\omega$ \emph{contains diagonals }if  for
any sequence $(\xb_n)$ in $S$ with $\xb_n\keq
(x_{n,i})_{i=1}^\infty$ there exist $i_1\kle i_2\kle\dots$ in $\N$
so that $(x_{n,i_n})_{n=1}^\infty$ belongs to~$S$. In particular every
subsequence of every member of $S$  belongs to $S$.
If $S$ contains diagonals, then the $S$-index of a tree on $X$ may be
measured via the Cantor-Bendixson index of the fine Schreier families
$\big(\cF_{\alpha}\big)_{\alpha<\omega_1}$.

\begin{prop}
  \label{P:5.2}\cite[Proposition 5]{OSZ2}
  Let $X$ be an arbitrary set and let $S\ksubset X^\omega$. If $S$
  contains diagonals, then for a hereditary tree $\cA$ on $X$ and for a countable
  ordinal $\alpha$ the following are equivalent.
  \begin{enumerate}
  \item[a)]
    $\alpha<I_S(\cA)$.
  \item[b)]
    There is a family $\big(x_F\big)_{F\in\cF_\alpha\setminus\{\emptyset\}} \ksubset \cA$ such that for $F=(m_1,m_2,\ldots,m_k)\in \cF_\alpha$ the branch 
    $\xb_F= \big( x_{\{m_1\}}, x_{\{m_1,m_2\}}, \dots,x_{\{m_1,m_2,\dots,m_k\}} \big)$ is in  $\cA$,  and 
    $\big(x_{F\cup\{n\}}\big) _{n>\max F}$ is in~$S$, if $F$ is not maximal in $\cF_\alpha$.
  \end{enumerate}
\end{prop}
\begin{defin}\label{D:5.3}  Let $\cF\subset [\N]^{<\omega}$ be regular, $S$  a set of sequences in the set $X$, and
 and $(x_A:A\in \cF)$  a tree in $X$ indexed by $\cF$. We call $(x_A:A\in \cF)$  an {\em $S$- tree} if for 
 every non maximal $A\in \cF$ the sequence 
 $(x_{A\cup\{n\}}: n\in \N, \text{ with } A\cup\{n\}\in \cF)$ is a sequence in $S$.
 
 If $X$ is a Banach space and  $S$ are the $w$-null sequences, we call   $(x_A:A\in \cF)$
 a  {\em $w$-null tree}. Similarly we define $w^*$-null trees in in $X^*$.
\end{defin}

\begin{remark}\label{R:5.4}
 In the case of $X=\N$ and $S=[\N]^\omega$   we deduce from Proposition \ref{P:5.2} that
if $\cA\subset[\N]^{<\omega}$ is hereditary and  compact, then $\CB(\cA)> \alpha$ if and only if
there is  an order isomorphism
$\pi:\cF_{\alpha} \to \cA$,  so that 
for all  $A\in\cF_\alpha\setminus \MAX(\cF_\alpha)$ and $n>\max(A)$ 
it follows that  $\pi(A\cup \{n\})=\pi(A)\cup\{s_n\}$, where $(s_n)$ is an increasing sequence
in $\{s\in \N: s>\max\pi(A)\}$.

\end{remark}

\begin{exs}\label{Ex:5.5}
  1.~\emph{The weak index. }Let $X$ be a separable Banach space. Let
  $S$ be the set of all weakly null sequences in $B_X$, the unit
  ball of $X$. We call the $S$-index of a tree $\cF$ on $X$
  \emph{the weak index of $\cF$ }and we shall denote it by
  $I_w(\cF)$. We shall use the term \emph{weak derivative }instead of
  $S$-derivative and use the notation $\cF'_{w}$ and
  $\cF^{(\alpha)}_w$.

  When the dual space $X^*$ is separable, the weak topology on
  the unit ball $B_X$, or on any bounded subset  of $X$ is metrizable. Hence in this case the set
  $S$ contains diagonals and
  Proposition~\ref{P:5.2} applies.

\noindent
2.\emph{The weak$^*$ index}.  We can define the weak$^*$ index similarly to the weak index. If $X$ is a separable Banach space
 the $w^*$ topology on $B_X^*$  is metrizable. This implies that the set $S$ 
 of all $w^*$ null sequences in $B_{X^*}$ is diagonalizable.
  We  call the $S$ index of a tree $\cF$ on $X^*$ 
 \emph{the weak$^*$ index of $\cF$ }and we shall denote it by
  $I_{w^*}(\cF)$. We shall use the term \emph{weak$^*$ derivative }instead of
  $S$-derivative and use the notation $\cF'_{w^*}$ and
  $\cF^{(\alpha)}_{w^*}$.

%%% Next paragraph could be deleted
  \noindent
  3.~\emph{The block index. }Let $Z$ be a Banach space with an FDD
  $E\keq(E_i)$. A \emph{block tree of $(E_i)$ in $Z$ }is a tree $\cF$ such that every element of $\cF$ is a (finite) block
  sequence of $(E_i)$. Let $S$ be the set of all  infinite
  block sequences of $(E_i)$ in $B_Z$. We call the $S$-index of a block
  tree $\cF$ of $(E_i)$ \emph{the block index of $\cF$ }and we shall
  denote it by $I_{bl}(\cF)$. We shall use the term \emph{block
  derivative }instead of $S$-derivative and use the notation
  $\cF'_{\mathrm{bl}}$ and $\cF^{(\alpha)}_{\mathrm{bl}}$. Note that the set $S$
  contains diagonals, and hence
  Proposition~\ref{P:5.2} applies.

  Note that $(E_i)$ is a shrinking FDD of $X$ if and only if every
  element of $S$ is weakly null. In this case we have
  \begin{equation}
    \label{E:5.1}
      I_{\mathrm{bl}}(\cF) \leq I_w(\cF)
  \end{equation}
  for any block tree $\cF$ of $(E_i)$ in $Z$. The converse is false in
  general, but it is true up to perturbations and without the
  assumption that $(E_i)$ is shrinking.
\end{exs}

\section{The Szlenk index}\label{S:6}

 Here we recall the definition and basic properties of the Szlenk
index, and  prove further properties that are relevant
for our purposes. 
 
Let $X$ be a separable  Banach space, and let $K$ be
a non-empty subset of $X^*$. For $\vp\kge 0$ set
$$K'_\vp
=\Big\{ x^*\in X^*: \exists (x^*_n)\subset K\,\, w^*-\lim_{n\to\infty} x^*_n= x^*\text{ and } \|x^*_n- x^*\|>\vp\Big\}.$$ 
and
define $K^{(\alpha)}_\vp$ for each countable ordinal $\alpha$ by
recursion as follows:
\begin{align*}
  K^{(0)}_\vp = K,
  K^{(\alpha+1)}_\vp = \big( K^{(\alpha)}_\vp
  \big)'_\vp\text{ for }\alpha\kle\omega_1, \text{ and }
  K^{(\lambda)}_\vp = \bigcap _{\alpha <\lambda}
  K^{(\alpha)}_\vp \text{ for limit ordinals }\lambda\kle\omega_1.
\end{align*}
Next, we associate to $K$ the following ordinal indices:
$$
\eta(K,\vp) = \sup \{
\alpha\kle\omega_1:\,K^{(\alpha)}_\vp\kneq\emptyset \}\ ,
\text{ 
and }
\eta(K) = \sup _{\vp>0} \eta(K,\vp)\ .
$$
Finally, we define the {\em Szlenk index $\Sz(X)$ of }$X$ to be
$\eta(B_{X^*})$, where $B_{X^*}$ is the unit ball of $X^*$.
\begin{rem}
  The original definition   of the derived sets $K'_\vp$ in ~\cite{Sz} is slightly different. It might lead to different values of $\Sz(K,\vp)$, but 
  to same values of $\Sz(K)$, and in particular to the same values of $\Sz(X)$. 
 \end{rem}

Szlenk used his index to show that there is no separable, reflexive
space universal for the class of all separable, reflexive spaces. This
result follows immediately from the following properties of the
function $\Sz(\cdot)$.
\begin{thm}\cite{Sz}
  \label{T:6.1}
  Let $X$ and $Y$ be separable Banach spaces.
  \begin{enumerate}  \item
    $X^*$ is separable if and only if $\Sz(X)\kle\omega_1$.
  \item
    If $X$ isomorphically embeds into $Y$, then $\Sz(X)\kleq \Sz(Y)$.
  \item
    For all $\alpha\kle\omega_1$ there exists a separable, reflexive
    space with Szlenk index at least $\alpha$.
  \end{enumerate}
\end{thm}
We also recall the following observation of \cite{AJO}  about the form of  the Szlenk index of a Banach space with separable dual.
\begin{prop}\label{P:6.2}\cite[Theorems 3.22 and 4.2]{AJO} If $X$ has a separable dual then there is an $\alpha<\omega_1$ with 
$\Sz(X)=\omega^\alpha$.
\end{prop}
The following Theorem combines several equivalent  description of the Szlenk index of a separable space not containing $\ell_1$

\begin{thm}\label{T:6.3} Assume that $X$ is a separable space not containing $\ell_1$ and $\alpha<\omega_1$. The following conditions are equivalent.
\begin{enumerate}
\item$\Sz(X)>\omega^\alpha$.
\item There is an $\vp>0$ and  a tree $(z^*_A:A\in \cS_\alpha)\subset B_{X^*}$,  so that 
for any non maximal $A\in \cS_\alpha$ we have
\begin{align}
 \label{E:6.3.1} &w^*-\lim_{n\to\infty} z^*_{A\cup\{n\}}=z^*_A\text{ and }
 \big\|z^*_A-z^*_{A\cup\{n\}}\big\|> \vp \text{ for $n>\max(A)$.}
\end{align}  
\item
 There is an $\vp>0$,  a tree $(z^*_A:A\in\cS_\alpha)\subset B_X^*$ and  a $w$-null tree $(z_A:A\in \cS_\alpha)\subset B_X$, so that 
\begin{align} 
  \label{E:6.3.2}  &z^*_B(z_A)>\vp ,\text{ for all $A,B\in\cS_\alpha\setminus\{\emptyset\}$, with $A\preceq B $,}\\
  \label{E:6.3.3}  &\big|z^*_{B}(z_A)\big| < \vp/2 \text{ for all $A,B\in \cS_\alpha\setminus\{\emptyset\}$,  with $A\not\preceq B$,}
 \intertext{and for all non maximal $A\in\cS_\alpha$ we have}
 \label{E:6.3.4} &w^*-\lim_{n\to\infty } z^*_{A\cup\{n\}}= z^*_A.
 \end{align}  
 \item There is an $N\in[\N]^{\omega}$, an $\vp>0$   and  a $w$-null tree tree $(x_A: A\in \cS_\alpha\cap [N]^{<\omega})$, so that 
 for every maximal $B$ in $\cS_\alpha\cap[N]^{<\omega}$ we have 
 $$\Big\|\sum_{A\preceq B} \zeta(\alpha,A) x_A\Big\|\ge \vp.$$

 \item There is an $\vp>0$, so that  $I_w(\cF_\vp)>\omega^\alpha$, where 
  $$\cF_\vp= \left \{(x_1,x_2,\ldots,x_l)\subset S_X:   \forall (a_j)_{j=1}^l\ksubset[0,1] \quad\Big\|\sum_{j=1}^l a_j x_j\Big\|\ge \vp \sum_{j=1}^l a_j           \right\}.$$
  \item There is an $\vp>0$ so that  $I_{w^*}(\cG_\vp)>\omega^\alpha$, where 
  $$\cG_\vp= \left \{(x^*_1,x^*_2,\ldots,x^*_l)\subset B_{X^*}:   \|x^*_j\|\ge \vp \text{ and } \Big\|\sum_{i=1}^j x^*_i\Big\|\le 1, \text{ for $j=1,2,\ldots, l$}          \right\}.$$
\end{enumerate}
\end{thm} 
\begin{proof} In order to show  ``(i)$\Rightarrow$(ii)'' we first prove the following 
\begin{lem} \label{L:6.4}Let $X$ be a separable Banach space and  $K\subset X^*$ be $w^*$-compact, $0<\vp<1$,  and $\beta<\omega_1$. Then for  every $x^*\in K^{(\beta)}_{\vp}$ there is  a family
 $(z^*_{(x^*,A)}:A\in \cF_\beta)\subset K$, so that 
 \begin{align}
 \label{E:6.4.1} &z^*_{(x^*,A)}\in K \text{ and }z^*_{(x^*,\emptyset)}=x^*,\\
 \label{E:6.4.2} &\text{if $A$ is not maximal in $\cF_\beta$, then $\|z^*_{(x^*,A)}\kminus z^*_{(x^*,A\cup\{n\})}\|>\vp$}
  \text{ for all $n>\max(A) $,}\\
   \label{E:6.4.3}  &\text{if $A$ is not maximal in $\cF_\beta$, then } z^*_{(x^*,A)}=w^*-\lim_{n\to\infty} z^*_{(x^*,A\cup\{n\}).}
 \end{align}
\end{lem}
\begin{proof} We will prove our claim by transfinite induction for all $\beta<\omega_1$. Let us first assume that $\beta=1$. For $x^*\in K'_{\vp}$ 
  choose  a sequence $(x^*_n)$ which $w^*$-converges to $x^*$, with  $\|x^*-x^*_n\|>\vp$, for $n\in\N$. Thus,  we can choose 
$z_{(x^*,\emptyset)}= x^*$, and  $z_{(x^*,\{n\})}= x^*_n$. This choice satisfies \eqref{E:6.4.1},  \eqref{E:6.4.2}, and   \eqref{E:6.4.3}, for $\beta=1$ 
 (recall that $\cF_1=\big\{\{n\}: n\kin\N\big\} \cup\{\emptyset\}$). 
 
 Now assume that our claim is true for all $\gamma<\beta$.  First assume that $\beta$ is a successor ordinal and let 
 $\gamma<\omega_1$ so that $\beta=\gamma+1$. Let $x^*\in K^{(\gamma+1)}_{\vp}$. Thus there is a sequence 
 $(x^*_n)\subset K^{(\gamma)}_{\vp} $ which $w^*$-converges to $x^*$, with $\|x_n^*-x^*\|>\vp$, for $n\kin\N$. 
 By our induction hypothesis we can choose for each $n\in\N$ a family 
 $\big(z^*_{(x^*_n,   A)} :A\in\cF_\gamma\big)$ satisfying  \eqref{E:6.4.1},  \eqref{E:6.4.2}, and   \eqref{E:6.4.3}, for 
 $\gamma$ and $x^*_n$ instead of $x^*$. 
 For every  $A\in\cF_{\gamma+1}$ it follows  that $A\setminus\{\min A\}\in \cF_\gamma$ and we define
 $z^*_{(x^*,\emptyset)}:= x^*$ and  for $A\in\cF_{\gamma+1}\setminus\{\emptyset\}$ 
 $$z^*_{(x^*,A)} :=z^*_{(x^*_{\min A}, A\setminus \{\min(A)\})}.$$
 It is then easy to see that  $\big(z^*_{(x^*,A)}:A\in\cF_{\gamma+1}\big)$ satisfies 
  \eqref{E:6.4.1},  \eqref{E:6.4.2}, and   \eqref{E:6.4.3}.

  Assume that $\beta<\omega_1$ is a limit ordinal and let $\big(\mu(\beta,n)  :n\kin\N\big)\subset(0,\beta)$, be the sequence
  of ordinals increasing to $\beta $, used to define $\cF_\beta$.  We abbreviate 
  $\beta_n=\mu(\beta,n)$, for $n\kin\N$.
   Let  $x^*\in K^{(\beta)}_{\vp}=\bigcap_{\gamma<\beta} K^{(\gamma)}_{\vp}$. Since $\beta_n+1<\beta$, we can use for each 
   $n\kin\N$ our induction hypothesis and choose 
   a family
   $$\big( z^*_{(n,x^*,A)}: A\in \cF_{\beta_n+1}\big) \subset X^*,$$
   satisfying   \eqref{E:6.4.1},  \eqref{E:6.4.2}, and   \eqref{E:6.4.3}, for 
 $\beta_n+1$. In particular it follows that 
 $$w^*-\lim_{j\to\infty}  z^*_{(n,x^*,\{j\})}=x^*, \text{ for all $n\kin\N$}.$$
 Since the $w^*$-topology  is metrizable on $K$ we can find an increasing sequence $(j_n:n\in\N)$ in $\N$, $j_n>n$, for $n\in\N$, so that 
 $$  w^*-\lim_{n\to\infty}  z^*_{(n,x^*,\{j_n\})}=x^*.$$
 Consider for $n\in\N$ the set
 $$\cF_{\beta_n+1}(j_n)= \big\{ A\in [\N]^{<\omega}: j_n\kle\min A\text{ and }\{j_n\}\cup A\in \cF_{\beta_n+1}\big\}= 
   \{A\in \cF_{\beta_n} : j_n\kle\min A  \}.$$
   Since $\cF_{\beta_n}$ is spreading, for $n\kin\N$, we can choose $L_n=\{l^{(n)}_1,l^{(n)}_2,\ldots\}\in [\N]^\omega$ so that
       $$\cF_{\beta_n}^{L_n}= \big\{ \{l^{(n)}_{a_1}, l^{(n)}_{a_2}, \ldots, l^{(n)}_{a_m}\}:   \{a_1,a_2,\ldots, a_m\}\in \cF_{\beta_n}\big\}\subset \cF_{\beta_n+1}(j_n).$$
    We define the map
    $$\phi_n:  \cF_{\beta_n}\to  \cF_{\beta_n+1}(j_n),   \{a_1,a_2,\ldots, a_m\}\mapsto \{l^{(n)}_{a_1}, l^{(n)}_{a_2}, \ldots, l^{(n)}_{a_m}\}.$$
    Then we put  for $A\in\cF_\beta$
    $$z^*_{(x^*, A)}=\begin{cases}  x^* &\text{if $A=\emptyset$,}\\
                                                       z^*_{(x^*, \{j_n\})}  &\text{if $A=\{n\}$ for some $n\in\N$, and}\\
                                                       z^*_{(x^*, \{j_n\}\cup \phi_n( B))} & \text{if $A=\{n\}\cup B$ for some $n\in\N$ and $B\in \cF_{\alpha_n}\setminus \{\emptyset\}$,}
                                                       \end{cases}$$
 which has the wanted property also in that case.                                                     
\end{proof}
We now continue with our proof of  ``(i)$\Rightarrow$(ii)'' of Theorem \ref{T:6.3}. 
Assuming now that $\Sz(X)>\beta=\omega^{\alpha}$, it follows that $\big[B_{X^*}\big]_\vp^{(\beta)}\not=\emptyset$ for some $\vp>0$. We  choose  $x^*\in\big[B_{X^*}\big]_\vp^{(\beta)}$ 
and apply Lemma \ref{L:6.4} to obtain a tree in $B_{X^*}$  indexed by $\cF_{\omega^\alpha}$ satisfying the conditions \eqref{E:6.4.1}, \eqref{E:6.4.2}, and 
 \eqref{E:6.4.3}. Now  Proposition \ref{P:2.13} and a  relabeling of the tree, yields (ii).

\noindent ``(ii)$\Rightarrow$(iii)''
 For $A\kin\cS_\alpha\ksetminus\{\emptyset\}$  we  define $A'=A\setminus\{\max(A)\}$.
 Now let $(A_m)_{m\in\N} $ be a consistent ordering of $\cS_\alpha$ (see Subsection \ref{SS:2.4}) . We write $A<_{lin}B$ or $A\le_{lin}B$ if $A=A_m$ and $B=A_n$, for $m<n$ or $m\le n$, respectively.

Let $\vp>0$ and $(z^*_A: A\kin\cS_\alpha)\subset B_{X^*}$, so that  \eqref{E:6.3.1} is  satisfied. 
Then choose for each $A\in\cS_\alpha\setminus \{\emptyset\}$  an element $x_A\in S_X$ so that 
$(z^*_{A}-z^*_{A'})(x_A)>\vp$. 

Let $0<\eta<\vp/8$, and let $\big(\eta(A): A\in\cS_\alpha\big)\ksubset(0,1)$ satisfy the following conditions
\begin{align} 
\label{E:6.4.5}&\text{$\big(\eta(A)\big)$  is decreasing with respect to the linear ordering $<_{lin}$,}\\
\label{E:6.4.6} &\sum_{A\in\cS_\alpha} \eta(A)<\eta,\\
\label{E:6.4.7}&\sum_{B\in\cS_\alpha, B>_{lin} A}  \eta(B)<{\eta(A)}, \text{ for all $A\in\cS_\alpha$, and }\\
\label{E:6.4.8} &\eta(A_m) <\frac12\frac{\eta}{m+2}, \text{for all $m\in\N$.}
\end{align}

Since $X$ does not contain a copy of $\ell_1$ we can apply Rosenthal's $\ell_1$ Theorem, and assume, possibly   after passing to a pruning, that for each non maximal $A\kin\cS_\alpha$ the sequence 
$(x_{A\cup\{n\}})_{n>\max(A)})$ is  weakly  Cauchy.
Since   $(z^*_{A\cup\{n\}}- z^*_A)_{n>\max(A)} $ is $w^*$-null we can assume, possibly after passing to a further pruning,  that  
 $(z^*_{A\cup\{n\}}- z^*_A)(x_{A\cup\{n-1\}})<\eta(A\cup\{n\})$, for all non maximal $A\in\cS_\alpha$ and  $n>1+\max(A)$.

Let  $z_\emptyset\keq0$. For a non maximal element $A\in\cS_\alpha$ and $n\kgr 1+\max(A)$ let
 $$z_{A\cup\{n\}}\keq\frac12(x_{A\cup\{n\}}-x_{A\cup\{n-1\}})$$
  and  
$$z_{A\cup\{\max(A)+1\}}=x_{A\cup\{\max(A)+1\}}.$$
Then the families $(z_A:A\in\cS_\alpha)$ and $(z^*_A:A\in\cS_\alpha)$ are in $B_X$ and $B_{X^*}$ respectively, 
 $(z_A:A\in\cS_\alpha)$ is weakly null and  $(z^*_A:A\in\cS_\alpha)$ satisfies \eqref{E:6.3.4}.

 Moreover it follows that 
\begin{equation}\label{E:6.4.9}
(z^*_{A}- z^*_{A'})(z_A)\ge \frac\vp2 - \eta(A),\text{ for all $A\in\cS_\alpha\setminus \{\emptyset\}$.}
\end{equation}

% For any non maximal  $B\in\cS_\alpha$,  there is a subsequence  $(z_{B\cup\{m_k\}},z^*_{B\cup\{m_k\}})_{k>\max(B)}$ of  
%   $(z_{B\cup\{n\}},z^*_{B\cup\{n\}})_{n>\max(B)}$ so that for all $k\in\N$
 %  \begin{align*}
  % \big|(z^*_{B\cup\{m_k\}}-z^*_B)(z_A)\big|< \eta(B\cup\{k\})\text{ and }&\big|z^*_A(z_{B\cup\{m_k\}})\big| <\eta(B\cup\{k\})\\ 
  % &\text{for all $A,B\in\cS_\alpha$, with $A\le_{lin} B\cup\{m_{k-1}\}$.}
  % \end{align*} 
  
  Since $w-\lim_{n\to\infty} z_{B\cup\{n\}}=0$ and $w^*-\lim_{n\to\infty} z^*_{B\cup\{n\}} =z^*_B$, for every non maximal $B\in\cS_\alpha$
   we can, after passing again to a pruning,  assume that  
      \begin{align}\label{E:6.4.10}
   \big|(z^*_{B}-z^*_{B'})(z_A)\big|< &\eta(B)\text{ and }  \big|z^*_{A}(z_B)\big| <\eta(B)
     \text{ for all $A,B\kin\cS_\alpha$, with   $A<_{lin} B$.} 
   \end{align}

We are left with verifying   \eqref{E:6.3.2} and  \eqref{E:6.3.3}  for $\vp/4$ instead of $\vp$.

To show \eqref{E:6.3.2} let  $A,B\kin\cS_\alpha\setminus\{ \emptyset\}$, with $A\preceq B$. We choose  $l\kin\N$ and $A\keq B_0\prec B_1\prec B_2\prec\ldots \prec B_l\keq B$,
so that $B_j'\keq B_{j-1}$, for $j\keq1,2,\ldots, l$, 
and deduce from \eqref{E:6.4.9}  and \eqref{E:6.4.10} 
\begin{align*}
z^*_B(z_A)&= \sum_{j=1}^l (z^*_{B_j}-z^*_{B_j'})(z_A) + (z^*_{A} -z^*_{A'})(z_A) +z^*_{A'}(z_A)\\
 &\ge \frac\vp2- \sum_{j=1}^l\eta(B_j)- 2\eta(A)
 >\frac\vp2-2\eta > \frac\vp4.\end{align*}

In order to show \eqref{E:6.3.3} let $A,B\in\cS_\alpha\setminus\{ \emptyset\}$, with $A\not\preceq B$.  We choose   $l\in\N$ and $\emptyset=B_0\prec B_1\prec B_2\prec\ldots \prec B_l=B$
so that $B'_j=B_{j-1}$, for $j=1,2\ldots, l$ and, since for every $j=1,2,\ldots, l$, either $A<_{lin} B_j$ or $B_j<_{lin} A$ we  deduce from \eqref{E:6.4.10}
and the conditions \eqref{E:6.4.6} and  \eqref{E:6.4.8} on  $\eta(\cdot)$, that 
\begin{align*}
|z^*_B(z_A)|&\le\Big| \sum_{j=1}^l (z^*_{B_j}-z^*_{B_j'})(z_A)\Big| +|z^*_\emptyset(z_A)|  \\
&\le \sum_{ A<_{lin}B_{j}}\big|(z^*_{B_j}-z^*_{B_j'})(z_A)\big| + \sum_{ A>_{lin}B_j} \big(|(z^*_{B_j}(z_A)|+|z^*_{B_{j-1}}(z_A)|\big)+\eta(A)\\
&\le \sum_{ A<_{lin} B_{j}} \!\eta(B_j) +2\sum_{B_{j}<_{lin} A}\! \eta(A) +\eta(A)\le (2\#\{j\kleq l: B_j\kle_{\text{lin}} A\}\kplus2)\eta(A)<\frac{\vp}4
\end{align*} 
which verifies  \eqref{E:6.3.3} and finishes the proof of our claims.

\noindent ``(iii)$\Rightarrow$(iv)'' Let $\vp>0$, $(z^*_A:A\in\cS_\alpha)$ , and   $(z_A:A\in\cS_\alpha)$ satisfy the condition in (iii). Then it follows for a
maximal $B\in\cS_\alpha$ that
$$\Big\|\sum_{A\preceq B} \zeta(\alpha, A) z_A\Big\|\ge  \sum_{A\preceq B} \zeta(\alpha, A) z^*_B(z_A)\ge \vp \sum_{A\preceq B} \zeta(\alpha, A)=\vp,$$
which proves our claim.

\noindent ``(iv)$\Rightarrow$(v)''  Assume that $N\in[\N]^{\omega}$, $\vp>0$, and $(x_A: A\in\cS_\alpha\cap[N]^{<\omega})\subset B_X$ satisfy  (iv).
For $B\in\MAX(\cS_\alpha\cap[N]^{<\omega})$ put 
$y_B=\sum_{A\prec B} \zeta(\alpha,A) x_A$,  and choose $y^*_B\in S_{X^*}$ so that 
 $y_B^*(y_B)=\|y_B\|>\vp$. 
  
  For $B\in \MAX(\cS_\alpha)\cap [N]^{<\omega}$, we define $f_B: B\to [-1,1], \quad b\mapsto y^*_B(x_{\{a\in B, a\le b\}})$.
 From Corollary \ref{C:4.5} it follows now  for $\delta=\vp/2$ that 
 $\CB(\cA_{\delta,N})= \omega^{\alpha}+1$
 where 
 $$\cA_{\delta,N}= \left\{ A\in \cS_\alpha\cap[N]^{<\infty}:\begin{matrix} \exists B\kin \MAX(S_\alpha\cap[N]^{<\infty}), \\ A\ksubset B,  \text{ and } f_B(a)\ge \delta \text{ for all } a\kin A\end{matrix}\right\}.$$
 But from Proposition \ref{P:5.2} and the Remark thereafter we deduce 
 that  there is a order isomorphism $\pi:\cF_{\omega^\alpha}\to \cA_{\delta,N}$ so that 
  for every non maximal $A\in \cF_{\omega^\alpha}$ and any $n>\max(A)$ it follows that   $\pi(A\cup\{n\})= \pi(A)\cup\{s_n\}$, for some increasing sequence  $(s_n)\subset \N$.  
 Putting $z_A=x_{\pi(A)}$ it follows that $(z_A)_{A\in\cF_{\omega^\alpha}}$ is a weakly null tree and for every 
 $A=\{a_1,a_2,\ldots,a_l\}$  it follows that $(z_{\{a_1,a_2,\ldots,a_i\}})_{i=1}^l\in \cF_\delta$. Applying again Proposition \ref{P:5.2} yields (v).
  
 \noindent ``(v)$\iff$(i)'' This  follows from \cite[Theorem 4.2]{AJO} where it was shown that   $\Sz(X)=\sup_{\vp>0} I_w(\cF_\vp)$, if $\ell_1$ does not embed into $X$.

\noindent ``(ii)$\iff $(vi)'' Follows from  Proposition \ref{P:2.14} and an application of  Proposition \ref{P:5.2}  to the tree $\cG_\vp$ on $B_{X^*}$, and 
 and $S=\{(x^*_n)\subset B_X^*:  w^*-\lim_{n\to\infty}=0\}$.
\end{proof}
\begin{rem} We note that in the implication (i)$\Rightarrow$(ii) the assumption that $\ell_1$ does not embed into $X$ was not needed.
\end{rem}
We will also need the following {\em dual version } of  Theorem \ref{T:6.3}.
\begin{prop}\label{P:6.5}  Assume that $X$ is a Banach space whose dual  $X^*$ is separable, with $\Sz(X^*)>\omega^\alpha$. 
Then there is an $\vp>0$, a tree  $(z_A:A\in \cS_\alpha)\subset B_X$, and a $w^*$-null tree 
$(z^*_A:A\in\cS_\alpha)\subset B_X^*$  so that 
\begin{align} 
  \label{E:6.5.1}  &z^*_A(z_B)>\vp ,\text{ for all $A,B\in\cS_\alpha\setminus\{\emptyset\}$, with $A\preceq B $,}\\
  \label{E:6.5.2}  &\big|z^*_{A}(z_B)\big| < \frac{\vp}{2}\text{ for all $A,B\in \cS_\alpha\setminus\{\emptyset\}$,  with $A\not\preceq B$.}
 \end{align}  
\end{prop}
\begin{proof}
Recall that, as stated above, the implication (i)$\Rightarrow$(ii) of Theorem \ref{T:6.3} holds even if the space, to which the theorem is applied, contains $\ell_1$.
Applying this implication to $X^*$, we find $\vp>0$ and $\{z^{**}_A: A\kin\cS_\alpha\}\subset B_{X^{**}}$, so that \eqref{E:6.3.1}  is satisfied.
Then choose for each $A\in\cS_\alpha\setminus \{\emptyset\}$  an element $x^*_A\in S_X$ so that 
$(z^{**}_{A}-z^{**}_{A'})(x^*_A)>\vp$.  Again let $(A_n)$ be a consistent enumeration of $\cS_\alpha$, and write $A_m<_{\text{lin}} A_n$ if $m<n$.
We also assume that $(\eta(A):A\in\cS_\alpha)\subset (0,1)$, has the property that 
\begin{equation}\label{E:6.5.3}
\sum_{A\in\cS_\alpha}\eta(A)< \frac{\vp}{32}.
\end{equation}
After passing to  a first pruning we can assume that for all non maximal $A\in\cS_\alpha$  the sequence $(x^*_{A\cup\{n\}})$ $w^*$-converges, and that 
for any $B\in\cS_\alpha$ the sequance $z^{**}_B( x^*_{A\cup\{n\}})$ converges so some number $r_{A,B} $ (for fixed $A,B\in\cS_\alpha$ we only need to pass to a subsequence
of $(A\cup\{n\}: n\in\N, A\cup\{n\}>_{\text{lin}}B)$). 
Since   $(z^{**}_{A\cup\{n\}}- z^{**}_A)_{n>\max(A)} $ is $w^*$-null we can assume, after passing to a second pruning,
that 
we have
\begin{equation}
\label{TEMP1New}
\big|(z^{**}_{B}- z^{**}_{B'})(x^*_{A})\big|< \eta(B) \text{  for all  $A,B\in\cS_\alpha$,  with $A<_{\text{lin}} B$},
\end{equation}

We put  $z^*_\emptyset=0$ and for any non maximal element $A\in\cS_\alpha$
$$z^*_{A\cup\{\max(A)+1\}}=x^*_{A\cup\{\max(A)+1\}}  \text{ and } z^*_{A\cup\{n\}}=\frac12(x^*_{A\cup\{n\}}-x^*_{A\cup\{{n-1}\}}) \text{ if  $n>\max(A)+1$.}$$
It follows that   $(z^*_A:A\kin\cS_\alpha)$ is a $w^*$-null tree in $B_{X^*}$
and that  for any  $A\kin\cS_\alpha$ 
\begin{equation*}
(z^{**}_A-z^{**}_{A'})(z^*_{A})>\frac\vp2-\frac{\eta(A)}2.
\end{equation*}
Since  $z^{**}_B(x^*_{A\cup\{n\}})$ converges to $r_{A,B}$ we can assume, after passing to a third pruning, that
\begin{equation}\label{TEMP2New}
|z_A^{**}(z_B^*)| <  \frac{\eta(B)}2 \text{ whenever $A<_{\text{lin}} B$.}
\end{equation}
and hence
\begin{equation}\label{E:6.5.8New}
 z^{**}_A(z^*_A)>\frac\vp2- \eta(A)\text{ for all $A\in\cS_\alpha$}.
\end{equation}

Since $w^*$-$\lim_nz^{**}_{A\cup\{n\}} = z^{**}_{A}$  we can assume, by passing to a further pruning, that 
\begin{equation}\label{E:6.5.10New}
\big|(z^{**}_{B}-z^{**}_{B'})(z^*_A)\big|< \eta(B) \text{ for all $A,B\in\cS_\alpha$, with $A<_{\text{lin}} B$.}
\end{equation}

Since $B_X$ is $w^*$-dense in $B_{X^{**}}$ we can choose, for every $A\in \cS_\alpha$, a vector $z_A\in B_X$ so that  
   \begin{align}
      \label{E:6.5.12New} &|z^*_{A}(z_B)-z^{**}_B(z^*_A)|<\eta(B), \text{ for all $A,B\in\cS_\alpha$, with  $A\le_{\text{lin}} B$.}
    \end{align}
    
Combining \eqref{E:6.5.10New} and \eqref{E:6.5.12New} we obtain that for all $A$ and $B$ in $\cS_{\alpha}$, with $A<_{\text{lin}}B $ we have
\begin{equation}
\label{TEMP5New}
|z_A^*(z_B) - z_{B'}^{**}(z_A^*)| \leqslant |z_A^*(z_B) - z_{B}^{**}(z_A^*)| + |(z_B^{**} - z_{B'}^{**})(z^*_A)| <2\eta(B) .
\end{equation}

Using that $(z^*_A:A\kin\cS_\alpha)$ is a $w^*$-null tree, we can pass to a further pruning, so that  
\begin{equation}
\label{E:6.5.13New}
|z_B^*(z_A)| <  \eta(B), \text{ for all $A,B\in \cS_\alpha$ with $A<_{\text{lin}} B$}.
\end{equation}
 We deduce from \eqref{E:6.5.10New} and \eqref{E:6.5.12New}, for $A,B\in\cS_\alpha$, with $A\leq_{\text{lin}} B'(<_{\text{lin}} B)$  that
 \begin{align}\label{E:6.5.14aNew}
 | z^*_A(z_B\kminus z_{B'})|&\le |z^*_A(z_B)\kminus z^{**}_B(z^*_A) |+|z^*_A(z_{B'})\kminus z^{**}_{B'}(z^*_A) |+|(z^{**}_B\kminus z^{**}_{B'})(z^*_A)|\\
 &\le2\eta(B)+\eta(B').\notag
  \end{align}   
By \eqref{TEMP5New}, \eqref{TEMP2New}, and \eqref{E:6.5.13New} for  $A,B\in\cS_\alpha$, with $B'<_{\text{lin}} A<B$ we obtain
\begin{align}
\label{TEMP6New}
| z^*_A(z_B\kminus z_{B'})| \leqslant |z_A^*(z_B) \kminus z_{B'}^{**}(z_A^*)| + |z_{B'}^{**}(z_A^*)| + |z_A^*(z_{B'})| \le 2\eta(B)\kplus2\eta(A).
\end{align}

We now claim that the  families $\{ z_A:A\in\cS_\alpha\}$ and $\{ z^*_A:A\in\cS_\alpha\}$ satisfy  
\eqref{E:6.5.1}  and  \eqref{E:6.5.2}.

  In order to verify \eqref{E:6.5.1}, let $A,B\in\cS_\alpha \setminus \{\emptyset\}$, with $A\preceq B$.
  Then  let $k\in\N$ and $B_j\in\cS_\alpha$, for $j=0,1,2,\ldots, k$, be such that 
  $ A=B_0\prec B_1\prec B_2\prec\cdots\prec B_k = B$, and $B'_j=B_j$, for $j=1,2,\ldots ,k$. 
   \begin{align*} 
    z^*_A( z_B) 
    &= \sum_{j=1}^k   z^*_{ A}(z_{B_j}-z_{B_j'}) + z^*_{ A}(z_{ A})\\
    &\ge  z^{**}_{ A}(z^*_{ A})  -  |z^{**}_{ A}(z^*_{ A})-z^*_{ A}(z_{ A})|-| z^*_A(z_{B_1}-z_A)|- \sum_{j=2}^k  | z^*_{ A}(z_{B_j}-z_{B_j'}) |\\
        &> \frac{\vp}2- 3\eta(A)- 2  \sum_{j=1}^k\eta(B_j)\ge \frac{\vp}4.
           \quad(\text{By \eqref{E:6.5.8New}, \eqref{E:6.5.12New}, \eqref{E:6.5.14aNew} and   \eqref{TEMP5New}})
     \end{align*}
which yields  \eqref{E:6.5.1} if we replace $\vp$ by $\vp/4$.

 In order to verify \eqref{E:6.5.2}, let $A,B\in\cS_\alpha \setminus \{\emptyset\}$, with $A\not\preceq B$.
 If $A>_{\text{lin}}B $ we deduce our claim from \eqref{E:6.5.13New}. If $ A\le_{\text{lin}} B$, and thus  $ A<_{\text{lin}} B$,
 we choose $k\in\N$ and  
  $B_0\prec B_1\prec B_2\prec \cdots \prec B_k= B$, with $B_j'=B_{j-1}$, for $j=1,2,\ldots, k$,
  and $B_0<_{\text{lin}} A<_{\text{lin}}B_1$. 
 Applying \eqref{TEMP6New},\eqref{E:6.5.13New},  \eqref{E:6.5.14aNew}, and finally \eqref{E:6.5.3} we obtain
\begin{align*}
\big| z^*_A( z_B)\big| 
&\le |z_{ A}^*(z_{B_1} - z_{B_0})| + |z_{ A}^*(z_{B_0}) | + \Big|\sum_{j=2}^k z^*_{ A}(z_{B_j}-z_{B_j'})\Big| \\
&\le2\eta(B_1)\kplus\eta(B_0)\kplus3\eta(A)\kplus \sum _{j=2}^k(2\eta(B_j)+\eta(B_{j-1})) \le \frac\vp8
\end{align*}
which proves our claim.
\end{proof}
\begin{ex}\label{Ex:6.6} Let us construct an  example of   families $\big(z_A:A\in[\N]^{<\omega}\big)\subset B_{c_0} $ and 
$\big(z^*_A:A\in[\N]^{<\omega}\big)\subset B_{\ell_1} $ satisfying Proposition \ref{P:6.5}. Let $<_{\text{lin}}$ again be a linear consistent  ordering of $[\N]^{<\omega}$.
We first choose a family $(\At: A\in[\N]^{<\omega}\big)\subset [\N]^{<\omega}$ with the following properties 
 \begin{align} 
\label{E:6.6.1}&\text{$\At$ is  a spread of $A$, for each $A\in  [\N]^{<\omega}$,  }\\
\label{E:6.6.2}&\text{$A\prec B$ if and only if  $\At\prec \Bt$,}\\
 \label{E:6.6.3} &\text{if $A,B\kin  [\N]^{<\omega}$, $\emptyset\not=A<_{lin} B$, and  $C\in[\N]^{<\omega}$  is the maximal element in $[\N]^{<\omega}$, so}\\
&\text{that  $C\preceq A$ and $C\preceq B$, then $(\tilde A\setminus \tilde C)\cap(\tilde B\setminus \tilde C)=\emptyset $}.\notag
\end{align} 
We define for $A\in[\N]^{<\omega}$ 
$$ z_A=\sum_{a\in \At} e_a\text{ and } z^*_A=e^*_{\max(\At)},$$
where $(e_j)$  and $(e^*_j)$ denote the unit vector bases in $c_0$ and $\ell_1$ respectively.
It is now easy to verify that  the tree $(z^*_A)$ is $w^*$-null and that \eqref{E:6.5.1} is  satisfied for any $\vp\in (0,1)$. In order to verify \eqref{E:6.5.2}
let  $A,B\in [\N]^{<\omega}$ with $A\not \preceq B$. If $A>_{\text{lin}} B$, then $\max(\At)\not\in \Bt$, and our claim follows. If  $A<_{\text{lin}} B$
 let $C\in[\N]^{<\omega}$ be the maximal element for which $C\preceq A$ and $C\preceq B$. It follows that $C\prec A$, but also that $C\prec B$, which implies by 
 \eqref{E:6.6.3} that $\max(\At)\not\in \Bt$, and, thus our claim. 
 \end{ex}  

\section{Estimating certain convex  combinations of blocks using the Szlenk index}\label{S:7}

In this section we will assume that $X$ has an FDD $(F_j)$. This means that $F_j\subset X$ is a finite dimensional subspace of $X$, for $j\in\N$, and that every $x$ has a unique representation 
as sum $x=\sum_{j=1}^\infty  x_j$, with $x_j\in F_j$, for $j\in\N$.
For $x=\sum_{j=1}^\infty  x_j\in  X$ we call $\supp(x)=\{ j\in\N: x_j\not=0\} $ the {\em support of $x$ }(with respect to $(F_j)$),
and the smallest interval in $\N$ containing $\supp(x)$ is  called the {\em range  of $x$ }(with respect to $(F_j)$), and is denoted by 
$\ran(x)$. A (finite or infinite) sequence $(x_n)\subset X$ is called a {\em block  (with respect to $(F_j)$)}, if $x_n\not=0$, for all $n\in\N$, and
 $\supp(x_n)<\supp(x_{n+1})$, for all $n\in \N$ for which  $x_{n+1}$ is defined.
 
 We call an FDD {\em shrinking } if every bounded block $(x_n)_{n=1}^\infty$ is weakly null. As in the  case of bases,  $X^*$ is  separable, and thus  $\Sz(X)<\omega_1$,
  if $X$ has a shrinking FDD.

\begin{thm}\label{T:7.1}
Let $X$ be a Banach space with a shrinking FDD and $\alpha$ be a countable ordinal number with $\Sz(X)\le \omega^\alpha$.
Then for every $\vp>0$ and $M\in[\N]^\omega$, there exists $N\in[M]^\omega$ satisfying the following: for every $B = \{b_1,\ldots,b_d\}$ in $\cS_\alpha\cap[N]^{<\omega}$ and sequence $(x_i)_{i=1}^d$ in $B_X$,
with $\ran(x_j)\subset (b_{j-1},b_{j+1})$ for $j=1,\ldots,d$ (where $b_0 = 0$ and $b_{d+1} = \infty$), we have
\begin{equation}
\label{E:7.1.1}
\Big\|\sum_{j=1}^d\zeta(\alpha,{B}_j)x_j\Big\| < \vp,
\end{equation}
where ${B}_j = \{b_1,\ldots,b_j\}$ for $j=1,\ldots,d$.
\end{thm}

\begin{proof}
It is enough to find $N$ in $[M]^\omega$ so that \eqref{E:7.1.1} holds whenever $B\in\MAX(\cS_\alpha)\cap[N]^{<\omega}$. Indeed, if \eqref{E:7.1.1}
holds for all $B$ in $\MAX(S_\alpha)\cap[N]^{<\omega}$, then  for any $A\in\cS_\alpha\cap[N]^{<\omega}$ and $(x_i)_{i=1}^{\#A}$ satisfying the assumption of  Theorem \ref{T:7.1},
one may extend $A$ to any maximal set $B$ and extend the sequence $(x_k)_{k=1}^{\#A}$ by concatenating the zero vector $\#B-\#A$ times.
Towards a contradiction, we assume that such a set $N$ does not exist. Applying  Proposition \ref{P:3.1}  
to the Partition $(\cF,\cS_\alpha\setminus \cF)$ of $\cS_\alpha$, where
$$\cF=\left\{ B=\{b_1,b_2,\ldots ,b_n\}\kin \MAX(\cS_\alpha): 
\begin{matrix}\exists (x_j)\ksubset B_X,\, \ran(x_j)\ksubset (b_{j-1},b_{j+1}), \text{ for }j\keq1,2,\ldots, n\\
                             \big\|\sum_{j=1}^n \zeta(\alpha, \{b_1,b_2,\ldots, b_j\})x_j\big\|>\vp 
\end{matrix}\right\},$$
yields that there is $L$ in $[M]^\omega$, so that for all $B = \{b_1^B,\ldots,b^B_{d_B}\}$ in $\MAX(\cS_\alpha)\cap[L]^{<\omega}$,
there exists a sequence $(x_i^B)_{i=1}^{d_B}$ in $B_X$ with $\ran(x_j)\subset(b_{j-1}^B,b_{j+1}^B)$ for $j=1,\ldots,d_B$, so that
\begin{equation}
\label{E:7.1.2}
\Big\|\sum_{j=1}^{d_B}\zeta(\alpha,B_j^B)x_j^B\Big\| \ge \vp,
\end{equation}
where $B_j^B = \{b_1^B,\ldots,b_j^B\}$ for $j=1,\ldots,d_B$. For $A\preceq B$, if $A = B_j^B$, we use the notation $x_j^B = x_A^B$. Note that, under this notation, \eqref{E:7.1.2} takes the more convenient form
\begin{equation}
\label{E:7.1.3}
\Big\|\sum_{A\preceq B}\zeta(\alpha,A)x_A^B\Big\| \ge \vp
\end{equation}
and that 
\begin{equation}
\label{E:7.1.4}
\ran(x_{A'}^B)\subset\left(\max(A''),\max(A)\right), \;\text{for all}\;A\preceq B\;\text{with}\;A'\neq\varnothing,
\end{equation}
where $A''=(A')'$ and $\max(\emptyset) = 0$.

We will now apply several stabilization and perturbation arguments to show that we may assume that for $B\in\MAX(\cS_\alpha)$ and $A\preceq B$ the vector $x^B_{A'}$,   only depends on $A$., and will
the be renamed $x_A$.

Using a compactness argument, Proposition \ref{P:3.1} yields the following: if $A\in\cS_\alpha\cap[L]^{<\omega}$ is non-maximal with $A'\neq\emptyset$,
then for $\delta>0$, there is $L'\in [L]^\omega$ so that for all $D_1$, $D_2$ in $\MAX(\cS_\alpha(A))\cap[L']^{<\omega}$, we have $\|x_{A'}^{A\cup D_1} - x_{A'}^{A\cup D_2}\| < \delta$.
Combining the above with a standard diagonalization argument we may pass to a further infinite subset of $L$, and a perturbation of the block vectors $x_{A'}^B$ with $B\in\MAX(\cS_\alpha)\cap[L]^{<\omega}$ and $\{\min(B)\}\prec A\preceq B$
(and perhaps pass to a smaller $\vp$ in \eqref{E:7.1.3}), so that for every $B_1$, $B_2$ in $\MAX(\cS_\alpha)\cap[L]^{<\omega}$ and $A$, with $A'\neq\emptyset$, so that $A\preceq B_1$ and $A\preceq B_2$,
we have $x_{A'}^{B_1} = x_{A'}^{B_2}$. For every $A\in\cS_\alpha\cap[L]^{<\omega}$, we call this common vector $x_A$. Note that $x_A$ indeed depends on $A$ and not only on $A'$. For $A$ so that $A' = \emptyset$,
i.e. for those sets $A$ that are of the form $A = \{n\}$ for some $n\in L$, chose any normalized vector $x_A$ with $\supp(x_A) = \{n\}$. Note that, using
\eqref{E:7.1.4},
we have
\begin{equation}
\label{E:7.1.5}
\ran(x_A) \subset (\max(A''),\max(A)) \;\text{for all}\;A\in\cS_\alpha\cap[L]^{<\omega}\;\text{with}\;A'\neq\varnothing,
\end{equation}
where $\max(\emptyset) = 0$ and if $A' = 0$, i.e. $A = \{n\}$ for some $n\in\N$, then $\ran(x_A) = \{n\}$.

Furthermore, fixing $0 < \delta <\vp/12$ and passing to an infinite subset of of $L$, again denoted by $L$, satisfying $\min(L) \ge 1/\delta$,
\eqref{E:7.1.3} and \eqref{E:4.8.1} yield that for all $B\in\MAX(\cS_\alpha)\cap[L]^{<\omega}$
\begin{eqnarray}\label{E:7.1.6}
\Big\|\sum_{A\preceq B}\zeta(\alpha,A)x_A\Big\| &\ge& \Big\|\sum_{A\preceq B}\zeta(\alpha,A')x_A\Big\| - 2\delta\\
&=& \Big\|0x_{\{\min(B)\}} + \sum_{\{\min(B)\}\prec A\preceq B}\zeta(\alpha,A')x_{A'}^B\Big\| - 2\delta\nonumber\\
&=& \Big\|\sum_{A\prec B}\zeta(\alpha,A)x_{A}^B\Big\| - 2\delta
 \ge   \vp - 3\delta.\nonumber
\end{eqnarray}

For $B\in\MAX(\cS_\alpha\cap[L]^{<\omega})$ and $i=0,1,2$ define
\begin{align*}
B^{(i)} &= \{A\preceq B: \#A\!\!\!\!\mod 3 = i\}.
\end{align*}
By the triangle inequality, for some $0\le i(B)\le 2$, we have
$\|\sum_{A\in B^{(i(B))}}\zeta(\alpha,A)x_A\|\ge \vp/3 - \delta.$
By Proposition \ref{P:3.1}, we may pass to some infinite subset of $L$, again denoted $L$, so that for all $B\in\MAX(\cS_{\alpha}\cap[L]^{<\omega})$,
we have $i(B) = i_0$, for some common $i_0\in\{0,1,2\}$. We shall assume that $i_0 = 0$, as the other cases are treated similarly. Therefore,
for all $B\in\MAX(\cS_\alpha\cap[L]^{<\omega})$ we have
\begin{equation}
\label{E:7.1.7}
\Big\|\sum_{A\in B^{(0)}}\zeta(\alpha,A)x_A\Big\|\ge \frac{\vp}{3} - \delta.
\end{equation}
Lemma  \ref{L:4.6} parts  (iii) and (v) also imply the following. If $B\in\MAX(\cS_\alpha\cap[L]^{<\omega})$, then
\begin{equation}
\label{E:7.1.8}
\sum_{\substack{A\in B^{(0)}:\\ l_1(A')=0}}\zeta(\alpha,A) + \sum_{\substack{A\in B^{(0)}:\\ l_1(A'')=0 }}\zeta(\alpha, A)
+ \sum_{\substack{A\in B^{(0)}\\l_1(A''')=0}}\zeta(\alpha, A) \le \frac{3}{\min(B)} \le 3\delta.
\end{equation}
Hence, if for $B\in\MAX(\cS_\alpha\cap[L]^{<\omega})$ we set
$$\hat B^{(0)} = B^{(0)}\setminus\{A\preceq B: l_1(A')=0, \text{ or }l_1(A'')=0, \text{ or } l_1(A''')=0\},$$ then
\begin{equation}
\label{E:7.1.9}
\Big\|\sum_{A\in\hat B^{(0)}}\zeta(\alpha,A)x_A\Big\|\ge \frac{\vp}{3} - 4\delta.
\end{equation}

If $L = \{\ell_1,\ell_2,\ldots,\ell_k,\ldots\}$, define
$N = \{\ell_3,\ell_6,\ldots,\ell_{3k},\ell_{3(k+1)},\ldots\}. $

For each $A \in\cS_\alpha\cap[N]^{<\omega}$, with $(\#A)\!\!\!\!\mod 3 = 0$, we define $\tilde A\in[L]^{<\omega}$ as described bellow.
If $A =\{a_1,\ldots,a_d\}$, where $a_j = \ell_{3b_j}$ and $A_j = \{a_1,\ldots,a_j\}$ for $1\le j\le d$, we define the elements of a set $\At=
\{\at_1,\at_2,\ldots,\at_d\}$ in groups of three as follows. If $j\!\!\!\!\mod3=0$ we put
$$(\at_{j-2}, \at_{j-1},\tilde a_j )=
\begin{cases}
	(a_{j-2}, a_{j-1},	a_j)  & \mbox{if } l_1(A_j')=0 \mbox{ or }  l_1(A_{j-1}')=0 \mbox{ or } l_1(A_{j-2}')=0\\
		(\ell_{3b_{j}-2}\ell_{3b_{j}-1},a_j)&\mbox{if }l_1(A_j'), l_1(A_{j-1}'), l_1(A_{j-2}')\not=0.
	\end{cases}$$ 
It is not hard to see that $A$ and $\tilde A$ satisfy the assumptions of Lemma
\ref{C:4.7}, hence $\tilde A\in\cS_\alpha\cap[L]^{<\omega}$ and $\zeta(\alpha, A) = \zeta(\alpha, \tilde A)$.
Observe the following:
\begin{itemize}
\item[(a)] If $B\in\MAX(\cS_\alpha\cap[N]^{<\omega})$ and $A^{(1)}\prec A^{(2)}$ are in $\hat B^{(0)}$, then $\tilde A^{(1)}\prec\tilde A^{(2)}$.
\item[(b)] If $B\in\MAX(\cS_\alpha\cap[N]^{<\omega})$ and $A\in\hat B^{(0)}$,
\begin{equation}
\label{E:7.1.10}
\text{if }\max(A) = \ell_{3n}\text{ then we have }\max(\tilde A'') =\ell_{3n-2}.
\end{equation}
\end{itemize}
(a) is clear, while (b) follows from the fact that  $A\in \hat B^{(0)}$, implies that $d$ is divisible by $3$ and $l_1(A'')\not=0$.

We shall define a weakly null tree $(z_A)_{A\in\cS_\alpha\cap[N]^{<\omega}}$ so that for all $B\in\MAX(\cS_\alpha\cap[N]^\infty)$ we have
\begin{equation}
\label{E:7.1.11}
\Big\|\sum_{A\preceq B}\zeta(\alpha,A)z_A\Big\|\ge \frac{\vp}{3} - 4\delta.
\end{equation}
The choice of $\delta$ and statement (iv) of Theorem \ref{T:6.3} will yield a contradiction. 

For $A\in\cS_{\alpha}\cap[N]^{<\omega}$ define
$$z_A =
\left\{
	\begin{array}{ll}
		x_{\tilde A}  & \mbox{if }  \#A\!\!\!\!\mod3=0 \mbox{ and }l_1(A'),l_1(A''),l_1(A''')\not=0, \\
		0 &\mbox{else.}
	\end{array}
\right.$$

Let $C\in\MAX(\cS_{\alpha}\cap[N]^{<\omega})$ and, by (a) we can find $B$ be in $\MAX(\cS_{\alpha}\cap[L]^{<\omega})$ so that $\tilde A\preceq B$,
for all $A\preceq C$, with $\# A=0\mod 3$.  Then, one can verify that
$$\Big\|\sum_{A\preceq C}\zeta(\alpha,A)z_A\Big\| = \Big\|\sum_{A\preceq C}\zeta(\alpha,\tilde A) x_{\tilde A}\Big\| =
\Big\|\sum_{A\in\hat C^{(0)}}\zeta(\alpha, \tilde A) x_{\tilde A}\Big\| =
\Big\|\sum_{A\in\hat B^{(0)}}\zeta(\alpha, A) x_{A}\Big\|\ge\frac{\vp}{3} - 4\delta. $$

Let now $A\in\cS_\alpha\cap[N]^{<\omega}$ be non-maximal. We will show that w-$\lim_{n\in N}z_{A\cup \{n\}} = 0$. By the definition of
the vectors $z_A$, we need only treat the case in which $(\#A + 1) \!\!\!\!\mod 3 = 0$, i.e. when
$z_{A\cup\{n\}} = x_{\widetilde{A\cup\{n\}}}$ for all $n\in N$ with $n>\max(A)$. In this case,
by \eqref{E:7.1.5}, we deduce that if $\ell_{3n}\in N$, then
$$\min\supp(z_{A\cup\{\ell_{3n}\}}) = \min\supp(x_{\widetilde{A\cup\{\ell_{3n}\}}}) > \max(( \widetilde{A\cup\{\ell_{3n}\}})'') = \ell_{3n-2}$$
where the last equality follows by (b). Hence, $\lim_{n\in N}\min\supp(z_{A\cup\{n\}}) = \infty$.
The fact that the FDD of $X$ is shrinking completes the proof.
\end{proof}

\section{Two metrics on $\cS_\alpha$, $\alpha<\omega_1$}\label{S:8}
Since  $[\N]^{<\omega}$ with $\prec$ is a tree, with a unique root $\emptyset$, we could consider on $[\N]^{<\omega}$ the usual {\em tree distance}
which we denote by $d$:
For $A=\{a_1,a_2,\ldots,a_l\}\in[\N]^{<\omega}$, or  $B=\{b_1,b_2,\ldots, b_m\} $, we let $n=\max\{j\ge 0: a_i=b_i\text{ for } i=1,2,\ldots, j\}$, and then let 
$d(A,B)=l+m-2j$. But this distance will not lead to to the results we are seeking. Indeed, it was shown in \cite[Theorem 1.2]{BKL} that for any reflexive space $X$ $[\N]^{<\omega}$, with the graph metric embeds  
  bi-Lipschitzly  into $X$ iff and only if $\max(\Sz(X),\Sz(X^*))>\omega$.
 We will need {\em weighted graph metrics} on  $\cS_\alpha$.

\begin{enumerate}

\item The {\em weighted tree distance on $\cS_\alpha$}. For $A,B$  in $S_\alpha$ let $C$ be the  largest element in $\cS_\alpha$ (with respect to $\prec$)  so that
$C\preceq A$ and $C\preceq B$ (\ie $C$ is the common initial segment of $A$ and $B$) and then let 
$$d_{1,\alpha}(A,B)= \sum_{a\in A\setminus C} z_{(\alpha,A)}(a)+\sum_{b\in B\setminus C} z_{(\alpha,B)}(b)
=\sum_{C\prec D\preceq A}\zeta(\alpha, D) + \sum_{C\prec D\preceq B}\zeta(\alpha, D) . $$
\item The  {\em weighted  interlacing  distance on $\cS_\alpha$ }can be defined as follows. For $A,B\in \cS_\alpha$, say 
$A=\{a_1,a_2,\ldots,a_l\}$ and $B=\{b_1,b_2,\ldots, b_m\}$, with $a_1<a_2<\ldots <a_l$ and $b_1<b_2<\ldots<b_m$, we  put 
$a_0=b_0=0$ and  $a_{l+1}=b_{m+1}=\infty$, and  define 
$$d_{\infty,\alpha}(A,B)=
\max_{i=1,\ldots, m+1} \sum_{a\in A, b_{i-1}<a<b_i }z_{(\alpha,A)}(a)+
\max_{i=1,\ldots, l+1} \sum_{b\in B, a_{i-1}<b<a_i} z_{(\alpha,B)}(b).
$$ 
\end{enumerate}
\begin{rem} In order to explain $d_{ \infty,\alpha}$ let us take some  sets $A=\{a_1,a_2,\ldots, a_m\}$ and $B=\{b_1,b_2,\ldots, b_m\}$ in  $\cS_\alpha$ and 
and fix some $i\in\{0, 1,2,\ldots, m\}$. Now we measure how large the part of $B$ is which lies between $a_i$ and $a_{i+1}$ (as before $a_0=0$ and $a_{m+1}=\infty$)
by putting
$$ m_i(B):=\sum_{j, b_j\in(a_i,a_{i+1})} \zeta(\alpha,\{b_1,b_2,\ldots, b_j\}).$$
Then we define $m_j(A)$, for $j=1,2,\ldots, n$ similarly, and put
$$d_{\infty,\alpha}(A,B)=\max_{1\le i\le m} m_i(B)+\max_{1\le j\le n} m_j(A) .$$
We note that if $C$ is maximal so that $C\preceq A$ and $C\preceq B$, and if $A\setminus C< B\setminus C$, then
$d_{1,\alpha}(A,B)= d_{\infty,\alpha}(A,B)$.
\end{rem}
\begin{prop}\label{P:8.1}
The metric space $(\cS_\alpha, d_{1,\alpha,})$ is {\em  stable}, \ie 
for any sequences $(A_n)$ and $(B_n)$ in $\cS_\alpha$ and any ultrafilter $\cU$ on $\N$ it follows that 
$$\lim_{m\in \cU}\lim_{n\in\cU}  d_{1,\alpha}(A_m,B_n)=
\lim_{n\in\cU} \lim_{m\in\cU}\ d_{1,\alpha}(A_m,B_n),$$
while $(\cS_\alpha, d_{\infty,\alpha})$  is not stable.
\end{prop}
\begin{proof} Since $\cS_\alpha$ is compact in the product topology we can pass to  subsequences and assume 
that we can write  $A_n$ as 
$A_n=A\cup \At_n$,  with $\max(A)<\min(\At_n)\le\max(\At_n)<\min(\At_{n+1})$,
and $B_n$ as 
$B_n=A\cup \Bt_n$,  with $\max(B)<\min(\Bt_n)\le\max \Bt_n<\min \Bt_{n+1}$, for $n\in\N$.
We can also assume the sequences
$$(r^A_n)= \Big( \sum_{A\prec D\preceq \At_n} \zeta(\alpha,D) \Big)\text{ and }
(r^B_n) = \Big( \sum_{B\prec D\preceq \Bt_n}   \zeta(\alpha,D) \Big) $$
converge to some numbers $r^A$ and $r^B$ respectively.

Let $C\in\cS_\alpha$ be the maximal element for which $C\preceq A$ and $C\preceq B$.
For $m\in\N$ large enough and $n>m $, large enough (depending on $m$) it follows that

$$d_{1,\alpha}(A_m,B_n)= \sum_{C\prec D\preceq A}   \zeta(\alpha,D) +\sum_{C\prec D\preceq B}   \zeta(\alpha,D) +r^A_m +r^B_n.$$
 Also for $n\in\N$ large enough and $m>n $, large enough (depending on $n$) it follows  also that
$$d_{1,\alpha}(A_m,B_n)= \sum_{C\prec D\preceq A}   \zeta(\alpha,D) +\sum_{C\prec D\preceq B}   \zeta(\alpha,D) +r^A_m +r^B_n.$$
 Thus we have 
 \begin{align*}  \lim_{m\in \cU}\lim_{n\in\cU}  d_{1,\alpha}(A_m,B_n)&= \sum_{C\prec D\preceq A}   \zeta(\alpha,D) +\sum_{C\prec D\preceq B}   \zeta(\alpha,D) +r^A +r^B \\
&=\lim_{n\in\cU} \lim_{m\in\cU}\ d_{1,\alpha}(A_m,B_n).
\end{align*}

To prove our second claim we 
let   $(A_n)$ be a sequence of  nonempty elements in $\cS_\alpha$ with  
 $d(\emptyset,A_n)=\sum_{D\preceq A_n} \zeta(\alpha, D)<1/3 $, and 
 $ A_{n-1}<\min A_{n}$, for $n\in\N$.
Then we let $B\in \cS_\alpha $ with $d(\emptyset, B)=\sum_{D\preceq B} \zeta(\alpha, D) \in (1/3,1/2]$ , and then choose
$\Bt_n\in\cS_\alpha$ so that  $\max B<\min \Bt_n\le \max \Bt_n<\min \Bt_{n-1}$,  and so that $B_n=B\cup B_n'$ is maximal in $\cS_\alpha $, for   $n\in\N$.

Then it follows 
$$\lim_{m\in\cU}\lim_{n\in\cU} d_{\infty,\alpha}(A_m,B_n)= 1-d(\emptyset , B)\le 2/3 $$
and 
$$\lim_{n\in\cU}\lim_{m\in\cU} d_{\infty,\alpha}(A_m,B_n)\ge1.$$
 \end{proof}

We can now conclude one direction of   Theorem A from James' characterization of reflexive spaces. 

   \begin{prop}\label{P:8.2} If $X$ is a non reflexive Banach  space then  for any $0<c<1/4$ and  every $\alpha>0$ there is a map 
   $\Phi_\alpha:\cS_\alpha\to X$, so that 
      \begin{equation}\label{E:8.2.1}
    cd_{\infty,\alpha}(A,B)\le \|\Phi(A)-\Phi(B)\|\le d_{1,\alpha}(A,B) \text{ for all $A,B\in\cS_\alpha$.}
    \end{equation}
   \end{prop} 
   \begin{rem} Our argument will show that if $X$ is non reflexive there is a sequence $(x_n)\subset B_X$ so that for all
   $\alpha<\omega_1$ the map 
    $$\Phi_\alpha:\cS_\alpha\to X, \quad A\mapsto \sum_{D\preceq A} \zeta{(\alpha, D)}   x_{\max(D)},$$
    satisfies \eqref{E:8.2.1}.
   \end{rem}
  
   \begin{proof} Let $\theta$ be any number in $(0,1)$. Then by  \cite{Ja0} 
   there is a normalized  basic sequence in $X$  whose basic constant is at most $\frac 2\Theta$ satisfying
 \begin{equation}\label{E:8.2.2}
\Big\|\sum_{j=1}^\infty a_j x_j \Big\| \ge\theta  \sum_{j=1}^\infty a_j \text{ for all $(a_j)\in c_{00}$, $a_j\ge 0$, for all $j\in\N$.}\end{equation}
    Thus its bimonotonicity constant is at most $\frac4\Theta$, which means that for $m\le n$ the projection
    $$P_{[m,n]}: \overline{\spa(x_j)}\to  \overline{\spa(x_j)}, \quad \sum_{j=1}^\infty a_j x_j\mapsto \sum_{j=m}^n a_j x_j,$$
    has norm at most $\frac 4\Theta$.
    
  We define
$$
  \Phi:\cS_\alpha\to X,\quad A\mapsto \sum_{D\preceq A} \zeta{(\alpha, D)}   x_{\max(D)}.$$
  
   For $A,B\kin\cS_\alpha$ we let $C$ be the maximal  element in $\cS_\alpha$ for which $C\preceq A$ and $C\preceq B$. Then
      \begin{align*}
   \|\Phi(A)-\Phi(B)\|&=\Big\| \sum_{C\prec D\preceq A } \zeta{(\alpha, D)}  x_{\max D} -\sum_{C\prec D\prec B} \zeta{(\alpha,D)}  x_{\max(D)} \Big\|\\
         &\le \sum_{C\prec D\preceq A } \zeta{(\alpha, D)}  +\sum_{C\prec D\preceq B} \zeta{(\alpha,D)}  
   = d_{1,\alpha}(A,B).
   \end{align*}
   
   On the other hand, if we write $A=\{a_1,a_2,\ldots, a_l\}$ and  put $a_0=0$, and $a_{l+1}=\infty$, it follows for all $i=1,2,\ldots, l+1$,
 that 
 \begin{align*}
 \|\Phi(A)-\Phi(B)\| \ge \frac{\Theta}4 \big\|P_{(a_{i-1}, a_i)}\big(\Phi(A)-\Phi(B)\big)\big\|\ge 
    \frac{\Theta^2}4 \sum_{a_{i-1}<b<a_i} z_{(\alpha,B)}(b).
   \end{align*} 
   Similarly, if we write $B=\{b_1,b_2,\ldots, b_m\}$ and   put $b_0=0$, and $b_{m+1}=\infty$, it follows for all $j=1,2,\ldots, m+1$ that 
    \begin{align*}
 \|\Phi(A)-\Phi(B)\| \ge 
    \frac{\Theta^2}4 \sum_{b_{j-1}<a<b_j} z_{(\alpha,B)}(a).
   \end{align*} 
 Thus for any $i=1,2,\ldots, l$ and any  $j=1,2,\ldots, m$
 
 \begin{align*}
 \|\Phi(A)-\Phi(B)\| \ge 
  \frac{\Theta^2}8\Bigg[\sum_{a_{i-1}<b<a_i} z_{(\alpha,B)}(b)  +  \sum_{b_{j-1}<a<b_j} z_{(\alpha,B)}(a)\Bigg],
\end{align*}
   which implies our claim.
      \end{proof}
We finish this section with an observation which we will need later.
\begin{lem}\label{L:8.3}
Let $\xi$ and $\gamma$ be countable ordinal numbers with $\gamma < \beta = \omega^{\xi}$. Let $B_1 < \cdots <B_j$ be in $\MAX(\cS_{\beta\gamma})$, so that $\bar{B} = \{\min(B_j): 1\le j\le d\}$ is a non-maximal $\cS_\beta$ set
with $l_1(\bar{B}) > 0$ ($l_1(A)$ for $A\in[\N]^{<\omega}$ has been defined before Lemma \ref{L:4.6}) and set $D = \cup_{j=1}^d B_j\in\cS_{\beta(\gamma + 1)}$. Then for every $A$, $B$ in $\cS_{\beta\gamma}$ with $D<A$ and $D<B$ we have
\begin{align}
d_{1,\beta\gamma}(A,B) &= \frac{1}{\zeta(\beta,\bar{B})}d_{1,\beta(\gamma+1)}(D\cup A, D\cup B)\;\text{and}\label{E:8.3.1}\\
d_{\infty,\beta\gamma}(A,B) &= \frac{1}{\zeta(\beta,\bar{B})}d_{\infty,\beta(\gamma+1)}(D\cup A, D\cup B).\label{E:8.3.2}
\end{align}
\end{lem}

\begin{proof}
We will only prove \eqref{E:8.3.1}, as the proof of \eqref{E:8.3.2} uses the same argument. Let $C$ be the maximal element in $\cS_{\beta\gamma}$ so that $C\preceq A$ and $C\preceq B$.
Note that $\tilde C = D\cup C$ is the largest element of $\cS_{\beta(\gamma+1)}$ so that $\tilde C\preceq D\cup A$ and $\tilde C\preceq D\cup B$. Define $\tilde B_1 = \bar B \cup \{\min(A)\}$ and $\tilde B_2 = \bar B\cup \{\min(B)\}$ and observe that,
since $l_1(\bar B) > 0$, $\zeta(\beta,\tilde B_1) = \zeta(\beta,\tilde B_2) = \zeta(\beta,\bar B) $. Using \eqref{E:4.4.1} in Proposition \ref{P:4.4} we conclude the following:
\begin{equation}
\label{E:8.3.3}
\sum_{a\in (D\cup A)\setminus \tilde C}z_{(\beta(\gamma+1), D\cup A)}(a) = \zeta(\beta,\tilde B_1)\sum_{a\in A\setminus C} z_{(\beta\gamma,A)}(a) = \zeta(\beta,\bar B)\sum_{a\in A\setminus C} z_{(\beta\gamma,A)}(a)
\end{equation}
and similarly we obtain
\begin{equation}
\label{E:8.3.4}
\sum_{a\in (D\cup B)\setminus \tilde C}z_{(\beta(\gamma+1), D\cup B)}(a)  = \zeta(\beta,\bar B)\sum_{a\in B\setminus C} z_{(\beta\gamma,B)}(a).
\end{equation}
Applying \eqref{E:8.3.3} and \eqref{E:8.3.4} to the definition of the $d_{1,\alpha}$ metrics, the result easily follows.
\end{proof}

 \section{The Szlenk index and embeddings of $(\cS_\alpha, d_{1,\alpha})$ into $X$}\label{S:9}
 
In this section we show the following two results, Theorem \ref{T:9.1} and \ref{T:9.3}, which  establish a proof of Theorem B.

   \begin{thm}\label{T:9.1}   Let $X$ be a separable Banach space  and $\alpha$ a countable ordinal. 
   Assume that  $\Sz(X)\kgr\omega^\alpha$
   then $(\cS_\alpha, d_{1,\alpha})$ can be bi-Lipschitzly embedded into $X$ and $X^*$.
   \end{thm}
   Before proving Theorem \ref{T:9.1}  we want to first cover the case that $\ell_1$ embeds into  $X$.
   \begin{ex}\label{Ex:9.2} For each $\alpha<\omega_1$  we want to define  a bi-Lipschitz embedding of $(\cS_\alpha,d_{1,\alpha})$ into  a Banach space  $X$  and and its dual  $X^*$ under the assumption that  $\ell_1$ embeds into $X$.
   We first choose for every $A\in[\N]^{<\omega}$ a spread $\At$ of $A$ as in Example \ref{Ex:6.6}.    
   Then we define for $\alpha<\omega_1$
   $$\Phi: \cS_\alpha\to \ell_1, \qquad A\mapsto \sum_{D\preceq A} \zeta(\alpha,D) e_{\max (\Dt)}.$$
   Since for $A,B\in\cS_\alpha$ it follows that 
\begin{align*}
\|\Phi(A)-\Phi(B)\| &=\Big\| \sum_{C\prec D\preceq A} \zeta(\alpha,D) e_{\max (\Dt)}- \sum_{C\prec D\preceq B} \zeta(\alpha,D) e_{\max (\Dt)}\Big\|\\
&=
   \sum_{C\prec D\preceq A} \zeta(\alpha,D)+
   \sum_{C\prec D\preceq B} \zeta(\alpha,D)= d_{1,\alpha} (A,B),
   \end{align*}
   where $C\in\cS_\alpha$ is the maximal element for which $C\preceq A$ and $C\preceq B$, it follows that $\Phi$ is an isometric embedding 
   of $(\cS_\alpha, d_{1,\alpha})$ into $\ell_1$.
    
    Thus, if $\ell_1$ embeds into $X$ then  $(\cS_\alpha,d_{1,\alpha})$ bi-Lipschitzly embeds into $X$. Secondly $\ell_\infty$ is a quotient of $X^*$ in that case, and since $\ell_1$ embeds into $\ell_\infty$,
    it follows easily that $\ell_1$ embeds into $X^*$, and, thus, that   $(\cS_\alpha,d_{1,\alpha})$ also  bi-Lipschitzly embeds into $X^*$. 
    \end{ex}
   
  \begin{proof}[Proof of Theorem \ref{T:9.1}] 
  Because of Example  \ref{Ex:9.2} we can assume that $\ell_1$ does not embed into $X$.
  Thus,  we can  apply Theorem  \ref{T:6.3},  (i)$\iff$(iii)
    and obtain $\vp>0$,  a tree 
   $(z^*_A:A\in\cS_\alpha)\subset B_{X^*}$ and a $w$-null tree $(z_A:A\in\cS_\alpha)\subset B_X$, so that 
  \begin{align} 
  \label{E:9.1.1}  &z^*_B(z_A)>\vp ,\text{ for all $A,B\in\cS_\alpha\setminus\{\emptyset\}$, with $A\preceq B $,}\\
  \label{E:9.1.2}  &\big|z^*_{B}(z_A)\big| \le \frac\vp2 \text{ for all $A,B\in \cS_\alpha\setminus\{\emptyset\}$,  with $A\not\preceq B$,}
   \end{align}  
  Then we define
  $$\Phi:\cS_A\to X,\quad A\mapsto \sum_{D\preceq A}  \zeta(\alpha,D) z_D.$$
  If $A,B\in\cS_\alpha $ and $C\in \cS_\alpha$ is the maximal element of $\cS_\alpha$ for which $ C\preceq A$ and $ C\preceq B$, we note that 
  \begin{align*}
  \|\Phi(A)-\Phi(B)\| &= \Big\|\sum_{C\prec D \preceq A }\zeta(\alpha,D) z_D- \sum_{C\prec D \preceq B }\zeta(\alpha,D) z_D\Big\|\\
  &\le 
  \sum_{C\prec D \preceq A }\zeta(\alpha,D)+\sum_{C\prec D \preceq B }\zeta(\alpha,D)=d_{1,\alpha}(A,B).\end{align*}
Moreover we obtain
  \begin{align*}
  \|\Phi(A)-\Phi(B)\| &= \Big\|\sum_{C\prec D \preceq A }\zeta(\alpha,D) z_D- \sum_{C\prec D \preceq B }\zeta(\alpha,D)z_D\Big\|\\
  &\ge z^*_{A} \Big(
  \sum_{C\prec D \preceq A }\zeta(\alpha,D) z_D -  \sum_{C\prec D \preceq B }\zeta(\alpha,D)z_D\Big)\\
& \ge\vp  \sum_{C\prec D \preceq A }\zeta(\alpha,D)  - \frac\vp2 \sum_{C\prec D \preceq B }\zeta(\alpha,D).
   \end{align*} 
    Similarly we can show that 
  $$ \|\Phi(A)-\Phi(B)\| \ge \vp \sum_{C\prec D \preceq B }\zeta(\alpha,D) -  \frac{\vp}2 \sum_{C\prec D \preceq A }\zeta(\alpha,D),$$ 
  and thus 
  $$\|\Phi(A)-\Phi(B)\| \ge \frac\vp4  \Bigg(\sum_{C\prec D \preceq A }\zeta(\alpha,D)+ \sum_{C\prec D \preceq B }\zeta(\alpha,D)\Bigg)
  =\frac\vp4 d_{1,\alpha}(A,B).$$
  
  In order to define a Lipschitz embedding from $(\cS_\alpha,d_{1,\alpha})$ into $X^*$, we also can assume that $\ell_1$ does not embed into $X$. Because otherwise 
   $\ell_\infty$ is a quotient of $X^*$, and thus since $ell_1$ embeds into $\ell_\infty$, it is easy to see that $\ell_1$ embeds into $X^*$, and again   we let
  $$\Psi: \cS_\alpha\to X^*, \quad A\mapsto \sum_{D\preceq A} \zeta(\alpha, D) z^*_D.$$
  As in the case of $\Phi$ it is easy to see that $\Psi$ is a  Lipschitz function with constant not exceeding the value $1$.
   Again if $A,B\in\cS_\alpha $ let  $C\in \cS_\alpha$ be the maximal element of $\cS_\alpha$ for which $ C\preceq A$ and $ C\preceq B$.
    In the case that $C\prec A$, we let $C^+\in \cS_\alpha$ be the minimal element for which   $C\prec C^+\preceq A$.
  We note that $C^+\not\preceq D$ for any $D\in\cS_\alpha $ with $C\prec D\preceq B$, and it follows therefore that 
    \begin{align*}
  \|\Psi(A)-\Psi(B)\| &= \Big\|\sum_{C\prec D \preceq A }\zeta(\alpha,D) z^*_D- \sum_{C\prec D \preceq B }\zeta(\alpha,D)z^*_D\Big\|\\
                              &\ge \Bigg(\sum_{C\prec D \preceq A }\zeta(\alpha,D) z^*_D- \sum_{C\prec D \preceq B }\zeta(\alpha,D)z^*_D\Bigg)(z_{C^+})\\                                                      
                              &\ge \vp\sum_{C\prec D \preceq A }\zeta(\alpha,D) -\frac\vp2\sum_{C\prec D \preceq B }\zeta(\alpha,D).
                               \end{align*} 
  If $C=A$ we arrive trivially to the same inequality,  Similarly we obtain that 
    $$ \|\Psi(A)-\Psi(B)\|  \ge \vp\sum_{C\prec D \preceq B }\zeta(\alpha,D) -\frac\vp2\sum_{C\prec D \preceq A }\zeta(\alpha,D).$$
This yields
  \begin{align*}
  \|\Psi(A)-\Psi(B)\| &\ge\frac\vp4\Bigg(\sum_{C\prec D \preceq A }\zeta(\alpha,D)+
  \sum_{C\prec D \preceq B }\zeta(\alpha,D)\Bigg)
  =\frac\vp4 d_{1,\alpha}(A,B)\end{align*}
  which finishes the proof of our claim.
  \end{proof}

    \begin{thm}\label{T:9.3}
    Assume that $X$ is a Banach space having a separable dual $X^*$, with $\Sz(X^*)>\omega^\alpha$ then $(\cS_\alpha, d_{1,\alpha})$ can be bi-Lipschitzly  embedded 
    into $X$.
    \end{thm} 
    \begin{proof}  Applying Proposition \ref{P:6.5}  we obtain  $\vp>0$  and 
  trees $(z^*_A:A\in\cS_\alpha)\subset B_X^*$ and $(z_A:A\in \cS_\alpha)\subset B_X$, so that $(z^*_A:A\in\cS_\alpha)$ is  $w^*$-null   and
\begin{align}
  \label{E:9.3.1}  &z^*_A(z_B)>\vp ,\text{ for all $A,B\in\cS_\alpha\setminus\{\emptyset\}$, with $A\preceq B $,}\\
  \label{E:9.3.2}  &\big|z^*_{A}(z_B)\big| \le \frac{\vp}2, \text{ for all $A,B\in \cS_\alpha\setminus\{\emptyset\}$,  with $A\not\preceq B$,}
 \end{align}  

We define 

$$\Phi: \cS_\alpha\to X,\quad A\mapsto \sum_{D\preceq A} \zeta(\alpha, D) z_A.$$
Again it is easy to see that $\Phi$ is a Lipschitz function with constant not exceeding the value $1$.
Moreover if if $A,B\in\cS_\alpha $ let  $C\in \cS_\alpha$ be the maximal element of $\cS_\alpha$ for which $ C\preceq A$ and $ C\preceq B$, 
    In the case that $C\prec A$, we let $C^+\in \cS_\alpha$ be the minimal element for which   $C\prec C^+\preceq A$.
  We note that $C^+\not\preceq D$ for any $D\in\cS_\alpha $ with $C\prec D\preceq B$, and it follows therefore that 
    \begin{align*}
  \|\Phi(A)-\Phi(B)\| &= \Big\|\sum_{C\prec D \preceq A }\zeta(\alpha,D) z_D- \sum_{C\prec D \preceq B }\zeta(\alpha,D)z_D\Big\|\\
                              &\ge z^*_{C^+}\Bigg(\sum_{C\prec D \preceq A }\zeta(\alpha,D) z_D- \sum_{C\prec D \preceq B }\zeta(\alpha,D)z_D\Bigg)\\                                                      
                              &\ge \vp\sum_{C\prec D \preceq A }\zeta(\alpha,D) -\frac\vp2\sum_{C\prec D \preceq B }\zeta(\alpha,D)\eta(D).
                               \end{align*} 
  If $C=A$ we trivially deduce   that 
 $$   \|\Phi(A)-\Phi(B)\|    \ge \vp\sum_{C\prec D \preceq A }\zeta(\alpha,D) \kminus\sum_{C\prec D \preceq B }\zeta(\alpha,D)\eta(D).$$
Similarly we prove that 
 $$   \|\Phi(A)-\Phi(B)\|    \ge \vp\sum_{C\prec D \preceq B}\zeta(\alpha,D) -\sum_{C\prec D \preceq A }\zeta(\alpha,D)\eta(D).$$
 We obtain therefore that 
 \begin{align*}
   \|\Phi(A)-\Phi(B)\|    &\ge \frac12 \Big(\sum_{C\prec D \preceq A }(\vp-\eta(A))\zeta(\alpha,D) +\sum_{C\prec D \preceq B }(\vp-\eta(A))\zeta(\alpha,D)\Big)\\
   &\ge \frac\vp4\Big( \sum_{C\prec D \preceq A }\zeta(\alpha,D) +\sum_{C\prec D \preceq B }\zeta(\alpha,D)\Big)=\frac\vp4 d_{1,\alpha}(A,B)
  \end{align*}
  which proves our claim.
\end{proof}

\section{Refinement argument }\label{S:10}

Before  providing a proof  of Theorem C,  and, thus, the still  missing implication of Theorem A,
we will  introduce  in this and the next section some more notation and  make some preliminary observations. 
The will  consider maps $\Phi:\cS_\alpha\to X$, satisfying weaker conditions  compared to the ones required by  Theorem A and C.  On the one hand it will make an  argument using transfinite induction possible, 
on the other  hand  it is sufficient to arrive to the wanted conclusions.

\begin{defin}\label{D:10.1}
Let $\alpha<\omega_1$. For $r\in(0,1]$ we define 
$$\cS^{(r)}_\alpha= \Big\{ A\in \cS_\alpha:\sum_{D\preceq A} \zeta(\alpha,D)\le r\Big\} .$$
It is not hard to see that $\cS_\alpha^{(r)}$ is a closed subset of $\cS_\alpha$, and hence  compact, and closed under restrictions.
We also put 
$$\cM^{(r)}_\alpha=\MAX(\cS^{(r)}_\alpha)=\{A\in\cS^{(r)}_\alpha: A\;\text{is  maximal in}\;\cS^{(r)}_\alpha \text{ with resp. to ``$\prec$''}\}$$ and for $A\in\cS_\alpha$,
$$\cM^{(r)}_\alpha(A)= \big\{ B\in \cS_\alpha(A): A\cup B\in \cM_\alpha^{(r)}\big\}.$$
\end{defin}

\begin{defin}\label{D:10.2}
Let $X$ be  a Banach space, $\alpha$ be a countable ordinal number, $L$ be an infinite subset of $\N$, and $A_0$ be a set in $\cS_\alpha$ that is either empty or a singleton.
A map $\Phi: \cS_\alpha(A_0)\cap[L]^{<\omega}\rightarrow X$ is called a  {\em semi-embedding of $\cS_\alpha\cap[L]^{<\omega}$ into $X$ starting at $A_0$}, if there is  a number $c>0$ 
so that
\begin{align}\label{E:10.2.1}
& \Big\|\Phi(A)-\Phi(B)\Big\|\le d_{1,\alpha} (A_0\cup A,A_0\cup B) \text{ for all $A,B\in\cS_\alpha(A_0)\cap[L]^{<\omega}$, and }\\
\intertext{for all  $A\in\cS_\alpha(A_0)\cap[L]^{<\omega}$, with $l_1(A_0\cup A)> 0$, for all $r\in(0,1]$,
and for all $B_1$, $B_2$ in $\cM^{(r)}(A_0\cup A)\cap[L]^{<\omega}$ with $B_1 < B_2$:}
\label{E:10.2.2} &\Big\|\Phi(A\cup B_1) \kminus\Phi(A\cup B_2)\Big\|\ge c d_{1,\alpha} (A_0\cup A\cup B_1,A_0\cup A\cup B_2 ).
\end{align}
($l_1(A)$ for $A\in[\N]^{<\omega}$ was introduced  in Definition \ref{D:4.5a}).

We call the  supremum  of all numbers $c>0$ so that \eqref{E:10.2.2}  holds for all $A\in\cS_\alpha(A_0)\cap[L]^{<\omega}$  and $B_1$, $B_2 \in\cM^{(r)}_\alpha(A_0\cup A)\cap[L]^{<\omega}$,
with $B_1<B_2$, the {\em semi-embedding constant of $\Phi$} and denote it by $c(\Phi)$.
\end{defin}

\begin{rem}
If  $\Phi:\cS_\alpha \rightarrow X$ is for some $0<c<C$ a $c$-lower $d_{\infty,\alpha}$ and $C$-upper $d_{1,\alpha}$ embedding,
 we can, after rescaling $\Phi$  if necessary, assume that $C=1$, and from the definition of $d_{1,\alpha}$ and $d_{\infty,\alpha}$
we can easily see that for every $A_0$, that is either empty or a singleton, and $L$ in $[\N]^\omega$,  the restrictions 
$\Phi|_{\cS_\alpha(A_0)\cap[L]^{<\omega}}X$ is a semi-embedding. 

 Assume  $\Phi: \cS_\alpha(A_0)\cap[L]^{<\omega}\rightarrow X$, is a semi-embedding of $\cS_\alpha\cap[L]^{<\omega}$ into $X$ starting at $A_0$. For 
 $A\in\cS_\alpha(A_0)$, with $A\neq\emptyset$,  we put $A'=A\setminus\{\max(A)\}$ and  define
 $$x_{A_0\cup A} =\frac1{\zeta(\alpha, A_0\cup A)}\left( \Phi( A)- \Phi(A')\right).$$
If $A_0 = \emptyset$ put $x_\emptyset  =\Phi(\emptyset)$, whereas if $A_0$ is a singleton, define $x_\emptyset = 0$ and $x_{A_0} = (1/\zeta(\alpha,A_0))\Phi(\emptyset)$.
Note that $\{x_\emptyset\}\cup\{x_{A_0\cup A}:A\in\cS_\alpha(A_0)\cap[L]^{<\omega}\}\subset B_X$. Recall that $\zeta(\alpha,\emptyset) = 0$ and hence, for $A\in \cS_\alpha(A_0)\cap[L]^{<\omega}$, we have 
$$\Phi(A)=x_\emptyset + \sum_{\emptyset\preceq D\preceq A_0\cup A } \zeta(\alpha,D) x_D.$$
We say in that case that 
{\em the family $\{x_\emptyset\}\cup\{x_{A_0\cup A}:A\in\cS_\alpha(A_0)\cap[L]^{<\omega}\}$  generates  $\Phi$}.

In that case the map $\Phi_0: \cS_\alpha(A_0)\cap[L]^{<\omega}\rightarrow X$ with $\Phi_0=\Phi-x_\emptyset$, i.e. for $A\in\cS_\alpha(A_0)\cap[L]^{<\omega}$
\begin{equation}
\label{E:10.1}
\Phi_0(A) = \sum_{\emptyset\preceq D\preceq A_0\cup A} \zeta(\alpha,D) x_D,
\end{equation}
is also a semi-embedding of $\cS_\alpha\cap[L]^{<\omega}$ into $X$ starting at $A_0$, with $c(\Phi_0)=c(\Phi)$.
\end{rem}

\begin{lem}\label{L:10.3}
Let $\gamma,\xi<\omega_1$, with $\gamma < \beta = \omega^{\xi}$,  and let  $B_1<\cdots<B_d$ be in $\MAX(\cS_{\beta\gamma})$,
 so that $\bar B = \{\min(B_j): 1\le j\le d\}$ is a non-maximal $\cS_\beta$ set
with $l_1(\bar B)>0$. Set $D = \cup_{j=1}^dB_j$, let $r\in(0,1]$ and let also $A\in\cM_{\beta\gamma}^{(r)}$ with $D<A$. Then, if $r_0 = \sum_{C\preceq D}\zeta(\beta(\gamma + 1),C) + \zeta(\beta,\bar B)r$,
we have that $A\in\cM_{\beta(\gamma + 1)}^{(r_0)}(D)$.
\end{lem}

\begin{proof}
From Proposition \ref{P:4.4} and Lemma \ref{L:4.6} (iii) we obtain that for $C \preceq A$ we have 
$$\zeta(\beta(\gamma+1),D\cup C) = \zeta(\beta,\bar B\cup\{\min(A)\})\zeta(\beta\gamma, C)= \zeta(\beta,\bar B)\zeta(\beta\gamma, C),$$
which  implies that $D\cup A\in\cS_{\beta(\gamma+1)}^{(r_0)}$.
If we assume that $D\cup A$ is not in $\cM_{\beta(\gamma+1)}^{(r_0)}$, there is $B\in\cS_{\beta(\gamma+1)}^{(r_0)}$ with $D\cup A \prec B$. Possibly after  trimming $B$, we may assume that $B' = D\cup A$. Define $B_0 = B\setminus D$.
Evidently, $A\prec B_0$ and $B_0' = A$. We claim that $B_0\in\cS_{\beta\gamma}$. If we assume that this is not the case, $A$ is a maximal $\cS_{\beta\gamma}$ set. This yields that $\sum_{C\preceq A}\zeta(\beta\gamma,C) = 1$
and hence $r = 1$ and $r_0 = \sum_{C\preceq D}\zeta(\beta(\gamma + 1)) + \zeta(\beta,\bar B) = \sum_{C\preceq D\cup A}\zeta(\beta(\gamma + 1),C)$, i.e. $D\cup A\in\cM_{\beta(\gamma+1)}^{(r_0)}$ which we assumed to be false. Thus,
we conclude that $B_0\in\cS_{\beta\gamma}$ and thus, using Proposition \ref{P:4.4} and the definition of $r_0$, $B_0\in\cS_{\beta\gamma}^{(r)}$, which is  a contradiction, as $A\in\cM_{\beta\gamma}^{(r)}$ and $A\prec B_0$.
\end{proof}

\begin{lem}\label{L:10.4}
Let $\alpha<\omega_1$, $N\in[\N]^{<\omega}$, $A_0$ be a subset of $\N$ that is either empty or a singleton, and $c\in(0,1]$.
Let  $\Psi:\cS_\alpha(A_0)\cap[N]^{<\omega}\to X$ be a semi-embedding of $\cS_\alpha\cap[N]^{<\omega}$ into $X$ starting at $A_0$, generated by a family of vectors 
$\{z_\emptyset\}\cup\{z_{A\cup A_0}: A\in \cS_\alpha(A_0)\cap[N]^{<\omega}\}$, so that $c(\Psi)>c$.
Let  $\vp<c(\Psi)-c$ and  $\{\tilde z_\emptyset\}\cup\{\tilde z_{A\cup A_0}: A\in \cS_\alpha(A_0)\cap[N]^{<\omega}\}\subset B_X$, with  
  $\|\tilde z_{A_0\cup A}- z_{A_0\cup A}\|< \vp$,  for all $A\in \cS_\alpha(A_0)\cap[N]^{<\omega}$ with $A_0\cup A\neq\emptyset$.
Then the  map  $\tilde \Psi:\cS_\alpha(A_0)\cap[N]^{<\omega}\rightarrow X$ defined by 
  $$\tilde \Psi(A)=\sum_{D\preceq A_0 \cup A} \zeta(\alpha,D) \tilde z_D, \text{ for } A\in \cS_\alpha\cap[N]^{<\omega},$$
 is a semi-embedding of $\cS_\alpha\cap[N]^{<\omega}$ into $X$ starting at $A_0$ with $c(\tilde\Psi)> c$.
\end{lem} 

\begin{proof}
For any $r\in(0,1]$, any $A\in \cS_\alpha(A_0)\cap[L]^{<\omega}$ and $B_1$, $B_2\in \cM_\alpha^{(r)}(A_0\cup A)\cap [N]^{<\omega}$, with $B_1< B_2$ we obtain
\begin{align*}
\Big\|\tilde\Psi(A\cup B_1)-\tilde\Psi(A\cup B_2)\Big\| =& \Big\|\sum_{A_0\cup A \prec D\preceq A_0\cup A\cup B_1} \zeta(\alpha, D) \zt_{D}-\sum_{A_0\cup A \prec D\preceq A_0\cup A\cup B_2} \zeta(\alpha, D) \zt_{D}\Big\|\\
\ge& \Big\|\sum_{A_0\cup A \prec D\preceq A_0\cup A\cup B_1} \zeta(\alpha, D) z_{D}-\sum_{A_0\cup A \prec D\preceq A_0\cup A\cup B_2} \zeta(\alpha, D) z_{D}\Big\|\\
&-\vp\left( \sum_{A_0\cup A \prec D\preceq A_0\cup A\cup B_1} \zeta(\alpha, D) + \sum_{A_0\cup A \prec D\preceq A_0\cup A\cup B_2} \zeta(\alpha, D)\right) \\
\ge& (c(\Psi)-\vp)d_{1,\alpha} (A_0\cup A\cup B_1, A_0\cup A\cup B_2 ),
\end{align*} 
which implies \eqref{E:10.2.2}. \eqref{E:10.2.1} follows from the fact that $\tilde z_{A_0\cup A}\kin B_X$ for all $A\kin\cS_\alpha(A_0)\cap[N]^{<\omega}$ with $A_0\cup A\neq\emptyset$.
\end{proof}

For the rest of the  section we will assume  that $X$ has a bimonotone  FDD $(E_n)$.  
For finite or cofinite sets $A\subset \N$, we denote
 the canonical projections from $X$ onto  $\overline{\spa(E_j:j\kin A)}$ by $P_A$, \ie
 $$P_A:X\to X,\,\, \sum_{j=1}^\infty x_j\mapsto  \sum_{j\in A} x_j, \text{ for $x=\sum_{j=1}^\infty x_j\in X$, with $x_j\in E_j$, for $j\in \N$,}$$
 and we write $P_j$ instead of $P_{\{j\}}$, for $j\in\N$.
We denote the linear span of  the $E_j$ by $c_{00}(E_j:j\in\N)$, \ie
$$c_{00}(E_j:j\in\N)=\Big\{ \sum_{j=1}^\infty x_j: x_j\in F_j, \text{ for $j\in\N$ and }\#\{j: x_j\not= 0\}<\infty \Big\}.$$

\begin{defin}\label{D:10.5}  Let $\alpha$ be a countable ordinal number, $M\in[\N]^{\omega}$, and $A_0$ be a subset of $\N$ that is either empty or a singleton.
A semi-embedding $\Phi:\cS_\alpha(A_0)\cap[M]^{<\omega}\rightarrow X$ of $\cS_\alpha\cap[M]^{<\omega}$ into $X$ starting at $A_0$, is said to be $c$-refined, for some $c \le c(\Phi)$, if the following conditions are satisfied:
\begin{itemize}

\item[(i)] the family $\{x_\emptyset\}\cup\{x_{A_0\cup A}:A\in\cS_\alpha(A_0)\cap[M]^{<\omega}\}$  generating  $\Phi$ is contained in $B_X\cap c_{00}(E_j:j\in\N)$,

\item[(ii)] for all $A\in\cS_\alpha(A_0)\cap[M]^{<\omega}$ with $A_0\cup A\neq\emptyset$ we have $$\max(A_0\cup A)\le \max\supp(x_{A_0\cup A}) < \min\{m\in M: m>\max(A_0\cup A)\}.$$

\item[(iii)] for all $r\in(0,1]$, $A\in\cS_\alpha(A_0)\cap[M]^{<\omega}$, with $l_1(A_0\cup A)>0$, and $B_1$, $B_2$ in $\cM_\alpha^{(r)}(A_0\cup A)\cap[M]^{<\omega}$, with $B_1 < B_2$, we have
\begin{equation*}
\begin{split}
\big\|P_{(\max\supp (x_{A_0\cup A}),\infty)}\left(\Phi(A\cup B_1) - \Phi(A\cup B_2)\right)\big\|  \ge cd_{1,\alpha}\left(A_0\cup A\cup B_1, A_0\cup A\cup B_2\right),
\end{split}
\end{equation*}

\item[(iv)] for all $r\in(0,1]$, $A\in\cS_\alpha(A_0)\cap[M]^{<\omega}$ and $B$ in $\cM_\alpha^{(r)}(A_0\cup A)\cap[M]^{<\omega}$ we have
$$\big\|P_{(\max\supp(x_{A_0\cup A}),\infty)}\big(\Phi(A\cup B)\big)\big\| \ge \frac{c}{2}\sum_{A_0\cup A\prec C\preceq A_0\cup A \cup B}\zeta(\alpha,C).$$

\end{itemize}
\end{defin}
\begin{remark}\label{R:10.6}
Let $\xi<\omega_1$,  $\gamma \le \beta = \omega^{\omega^\xi}$  be a limit ordinal, $0\kle c\kleq1$, and  $M\kin[\N]^{\omega}$.
If $a_0\kin\N$,  we note that $\cS_{\beta\gamma}(\{a_0\})\cap[M]^{<\omega}\keq \cS_{\beta\eta(\gamma,a_0)}(\{a_0\})\cap[M]^{<\omega}$,
 and $\zeta(\beta\gamma, \{a_0\}\cup D)\keq\zeta(\beta\eta(\gamma,a_0), \{a_0\}\cup D)$ for $D\kin S_{\beta\gamma}(\{a_0\})\cap[M]^{<\omega}$,
where  $\big(\eta(\gamma,n)\big)_{n\in\N}$ is the sequence provided by Proposition  \ref{P:2.6}.
It follows that a semi-embedding of $\cS_{\beta\gamma}(\{a_0\})\cap[M]^{<\omega}$ into $X$ starting at $\{a_0\}$ that is $c$-refined, 
is a $c$-refined semi-embedding of $S_{\beta\eta(\gamma,n)}(\{a_0\})\cap[M]^{<\omega}$ into $X$. 

Secondly,
if $\Phi:\cS_{\beta\gamma}\cap[M]^{<\omega}\rightarrow X$ is a semi-embedding of $\cS_{\beta\gamma}\cap[M]^{<\omega}$ into $X$ starting at $\emptyset$ that is $c$-refined, then for every $a_0\in M$  and  $N = M\cap[a_0,\infty)$
the map $\Psi=\Phi|_{\cS_{\beta\eta(\gamma,a_0)}(\{a_0\})\cap[N]^{<\omega}}$,  is a semi-embedding of $\cS_{\beta\eta(\gamma,a_0)}\cap[N]^{<\omega}$ into $X$ starting at $\{a_0\}$, that is $c$-refined.
Furthermore, $\Psi$ is generated by the family $\{x_{A_0\cup A}: A\in\cS_{\beta\eta(\gamma,a_0)}(\{a_0\})\cap[N]^{<\omega}\}$, where   $\{x_A: A\in\cS_{\beta
 \gamma}\cap[M]^{<\omega}\}$ is the family generating $\Phi$.
\end{remark}

\begin{lem}\label{L:10.7}
Let  $\xi, \gamma\kle\omega_1$, with $\gamma <\kle\beta \keq \omega^{\omega^\xi}$, $M\kin[\N]^{\omega}$, and $A_0$ be a subset of $\N$ that is either empty or a singleton.
Let also $\Phi:\cS_{\beta(\gamma \kplus 1)}(A_0)\cap[M]^{<\omega}\rightarrow X$ be a semi-embedding of $\cS_{\beta(\gamma+1)}\cap[M]^{\infty}$ into $X$ starting at $A_0$, that is $c$-refined.
The family generating $\Phi$ is denoted by  $\{x_\emptyset\}\cup\{x_{A_0\cup A}:A\in\cS_\alpha(A_0)\cap[M]^{<\omega}\}$.
Extend the set $A_0$ to a set $A_0\cup A_1$,  $A_0<A_1$, which can be written as  $A_0\cup A_1 = \cup_{j=1}^kB_j\kin \cS_{\beta(\gamma+1)}\cap[A_0\cup M]^{<\omega}$, where
$B_1\kle\cdots\kle B_k$ are in 
$\MAX(S_{\beta\gamma})$ and $\bar{B} \keq \{\min(B_j): 1\kleq j\kleq k\}$ is a non-maximal $\cS_\beta$ set with $l_1(\bar{B}) \kgr 0$.
 
Then, for $N = M\cap(\max(A_0\cup A_1),\infty)$ and $n_0 = \max\supp (x_{A_0\cup A_1})$, the map 
$$\Psi:\cS_{\beta\gamma}\cap [N]^{<\omega}\rightarrow X,\quad A\mapsto \frac{1}{\zeta(\beta,\bar B)} P_{\left(n_0,\infty\right)}\left(\Phi(A_1\cup A)\right)$$
is a semi-embedding of $\cS_{\beta\gamma}\cap[N]^{\omega}$ into $X$ starting at $\emptyset$, that is $c$-refined. Moreover, $\Psi$ is generated by
the family $\{z_A: A\in\cS_{\beta\gamma}\cap[N]^{<\omega}\}$, where
\begin{equation}\label{E:10.7.1}
 z_\emptyset = 0\text{  and }z_A = P_{(n_0,\infty)}(x_{A_0\cup A_1\cup A}) \text{ for  } A\kin \cS_{\beta\gamma}\cap[N]^{<\omega}\setminus \{\emptyset\}.
 \end{equation}
\end{lem}

\begin{proof}
By Lemma \ref{L:8.3} we easily obtain that for $A$, $B\in\cS_{\beta\gamma}\cap [N]^{<\omega}$:
\begin{equation}
\label{E:10.7.2}
\Big\|\Psi(A) - \Psi(B)\Big\| \le \frac{1}{\zeta(\beta,\bar B)}d_{1,\beta(\gamma + 1)}(A_0\cup A_1\cup A,A_0\cup A_1\cup B) = d_{1,(\beta\gamma)}(A,B),
\end{equation}
i.e. \eqref{E:10.2.1} from Definition \ref{D:10.2} is satisfied for $\Psi$. We will show that \eqref{E:10.2.1} from Definition \ref{D:10.2} is satisfied for $\Psi$ as well.
Let $r\kin(0,1]$, $A$ be in $\cS_{\beta\gamma}$, with $l_1(A)\kgr0$, and $B_1\kle B_2$ in $\cM_{\beta\gamma}^{(r)}(A)\cap[N]^{<\omega}$ (i.e. $A\cup B_1$, $A\cup B_2\kin\cM_{\beta\gamma}^{(r)}$).
Note that we have $l_{1}(A_0\cup A_1\cup A) \kgr 0$.
If we set $r_0 \keq \sum_{C\preceq A_0\cup A_1}\zeta(\beta(\gamma \kplus 1),C)\kplus \zeta(\beta,\bar B)r$, 
by Lemma \ref{L:10.3} we deduce that $A\cup B_1$ and $A\cup B_2$ are in $\cM_{\beta(\gamma+1)}^{(r_0)}(A_0\cup A_1)\cap[N]^{<\omega}$, \ie $B_1,B_2\kin\cM_{\beta(\gamma+1)}^{(r_0)}(A_0\cup(A_1\cup A))\cap[N]^{<\omega}$.
Definition \ref{D:10.5} (ii) implies $n_0\kleq \max\supp (x_{A_0\cup A_1\cup A})$ and, thus, by  Definition \ref{D:10.5} (iii) for $\Phi$  we deduce
\begin{align}
\label{E:10.7.3}
\Big\|\Psi(A&\cup B_1) - \Psi(A\cup B_2)\Big\|\\
& = \frac{1}{\zeta(\beta,\bar B)}\Big\|P_{(n_0,\infty)}\left(\Phi(A_1\cup A\cup B_1) - \Phi(A_1\cup A\cup B_2)\right)\Big\|\notag\\
&\ge \frac{1}{\zeta(\beta,\bar B)}\Big\|P_{(\max\supp(x_{A_0\cup A_1\cup A}),\infty)}\left(\Phi(A_1\cup A\cup B_1) - \Phi(A_1\cup A\cup B_2)\right)\Big\|\notag\\
&\ge \frac{c}{\zeta(\beta,\bar B)}d_{1,\beta(\gamma + 1)}(A_0\cup A_1\cup A\cup B_1,A_0\cup A_1\cup A\cup B_2)\notag\\
&=cd_{1,\beta\gamma}(A\cup B_1,A \cup B_2),\notag 
\end{align}
where the last equality follows from Lemma \ref{L:8.3}.
In particular, \eqref{E:10.7.2} and \eqref{E:10.7.3}
yield that $\Psi:\cS_{\beta\gamma}\cap[N]^{<\omega}\rightarrow X$ is a semi-embedding of $\cS_{\beta\gamma}\cap[N]^{\infty}$ into $X$ starting at $\emptyset$ with $c(\Psi)\ge c$.

It remains to show that $\Psi$ satisfies Definition \ref{D:10.5} (i) to (iv).
Observe that Definition \ref{D:10.5} (ii) implies that for $C\preceq A_0\cup A_1$, we have $P_{(n_0,\infty)}(x_C) = 0$.
We combine this with \eqref{E:4.4.1} of Proposition \ref{P:4.4} to obtain that for $A\in\cS_{\beta\gamma}\cap [N]^{<\omega}$:
\begin{align}\label{E:10.7.4}
\Psi(A) &= \frac{1}{\zeta(\beta,\bar B)}P_{(n_0,\infty)}(\Phi(A_1\cup A)) \\&
= \frac{1}{\zeta(\beta,\bar B)}\sum_{C\preceq A_0\cup A_1\cup A}\zeta(\beta(\gamma + 1),C)P_{(n_0,\infty)}(x_C)\notag\\
&= \frac{1}{\zeta(\beta,\bar B)}\sum_{A_0\cup A_1\prec C\preceq A_0\cup A_1\cup A}\zeta(\beta(\gamma + 1),C)P_{(n_0,\infty)}(x_C)\notag\\
&= \sum_{C\preceq A}\zeta(\beta\gamma,C)P_{(n_0,\infty)}(x_{A_0\cup A_1\cup C}).\nonumber
\end{align}
For $A\in\cS_{\beta\gamma}\cap[N]^{<\omega}$ define $z_A = P_{(n_0,\infty)}(x_{A_0\cup A_1\cup A})$ and $z_\emptyset = 0$.
It then easily follows by \eqref{E:10.7.4} that $\Psi$ is generated by the family  $\{z_{A}:A\in\cS_\alpha\cap[N]^{<\omega}\}$. Moreover, as $\max\supp(z_A) = \max\supp(x_{A_0\cup A_1\cup A})$,
it is straightforward to check that Definition \ref{D:10.5} (i) and  (ii) are satisfied for $\Psi$. Whereas, observing that for all $A\in\cS_{\beta\gamma}\cap[N]^{<\omega}$, with $A\neq\emptyset$ (which is the case when $l_1(A)>0$),
we have $\max\supp(z_A) = \max\supp(x_{A_0\cup A_1\cup A})\ge n_0$,
an argument similar to the one used to obtain \eqref{E:10.7.3} also yields that $\Psi$
satisfies Definition \ref{D:10.5} (iii) and (iv).
\end{proof}

The main result of this section is the following Refinement Argument.

\begin{lem}\label{L:10.8}
Let   $\alpha<\omega_1$, $M\in[\N]^{\omega}$, and $A_0$ be a subset of $\N$ that is either empty or a singleton.
Let also $\Phi:\cS_\alpha(A_0)\cap[M]^{<\omega}\rightarrow X$ be a semi-embedding of $\cS_\alpha\cap[M]^{<\omega}$ into $X$ starting at $A_0$.
Then, for every $c<c(\Phi)$, there exists $N\in[M]^{\omega}$ and a semi-embedding $\tilde \Phi:\cS_{\alpha}(A_0)\cap[N]^{<\omega}\rightarrow X$ of $\cS_\alpha\cap[N]^{<\omega}$ into $X$ starting at $A_0$,
that is $c$-refined.
\end{lem}

\begin{proof}
Put $\tilde c = (c(\Phi)+c)/2$. Let $\{x_\emptyset\}\cup\{x_{A_0\cup A}: A\in\cS_{\alpha}(A_0)\cap[M]^{<\omega}\}$ be the vectors generating $\Phi$ and choose $\eta>0$ with $\eta < c(\Phi) - \tilde c$. After shifting we can assume without loss of generality that $x_\emptyset=0$.
Set $\tilde x_\emptyset = 0$ and choose for each $A\in\cS_\alpha(A_0)\cap[M]^{<\omega}$, with $A_0\cup A\neq\emptyset$, a vector $\tilde x_{A_0\cup A}\in B_X\cap c_{00}(E_j:j\in\N)$ so that:
\begin{itemize}
 \item[(a)]   $\|\tilde x_{A_0\cup A} - x_{A_0\cup A}\| < \eta/2,$ and $\max(A_0\cup A)\le \max\supp(\tilde x_{A_0\cup A})$.
\end{itemize}
Moreover, recursively choose $\tilde m_1 < \cdots <\tilde m_k <\cdots$ in $M$ so that for all $k$ we have $\tilde m_{k+1} > \max \big\{\max\supp (\tilde x_{A_0\cup A}): A\in\cS_\alpha(A_0)\cap[\{\tilde m_1,\ldots,\tilde m_k\}]\big\}$.
Define $\tilde M = \{\tilde m_k: k\in\N\}$ and $\tilde \Phi:\cS_\alpha(A_0)\cap[\tilde M]^{<\omega}\rightarrow X$ so that for all $A\in\cS_\alpha(A_0)\cap[\tilde M]^{<\omega}$ we have
$$\tilde\Phi(A) = \sum_{C\preceq A_0\cup A}\zeta(\alpha,C)\tilde x_C.$$
By Lemma \ref{L:10.4}, $\tilde \Phi$ is a semi-embedding from $\cS_\alpha\cap[\tilde M]^{<\omega}$ into $X$ starting at $A_0$ for which $c(\tilde \Phi) > \tilde c > c$, and for which 
the conditions (i) and (ii) of Definition \ref{D:10.5} are satisfied.

The goal is to find $N\in[\tilde M]^\omega$ so that, restricting $\tilde\Phi$ to $\cS_\alpha(A_0)\cap[N]^{<\omega}$, Definition \ref{D:10.5} (iii) and (iv) are satisfied as well.
Put $M_0=\tilde M$. Recursively, we will  choose for every $k\in\N$ an infinite set $M_k\subset M_{k-1}$
so that for each $k\in\N$ the following conditions are met:
\begin{itemize}
 
 \item[(b)] $\min (M_{k-1}) < \min(M_k)$, and putting  $m_j  =\min(M_j)$ for $j=1,\ldots,k-1$, then
 
 \item[(c)] for every $A\in\cS_\alpha(A_0)\cap[\{m_1,\ldots,m_{k-1}\}]$ with $l_1(A_0\cup A) > 0$, $r\in(0,1]$,
 and $B_1$, $B_2\in\cM_\alpha^{(r)}(A_0\cup A)\cap[\{m_1,\ldots,m_{k-1}\}\cup M_k]^{<\omega}$ with $B_1 < B_2$, we have
 $$\big\|P_{(\max\supp(\tilde x_{A_0\cup A}),\infty)}\big(\tilde \Phi(A\cup B_1) - \tilde\Phi(A\cup B_2)\big) \big\|\ge \tilde c d_{1,\alpha} (A_0\cup A\cup B_1,A_0\cup A\cup B_2 ),$$
 \item[(d)] for every $A\in\cS_\alpha(A_0)\cap[\{m_1,\ldots,m_{k-1}\}]$ with $l_1(A_0\cup A) > 0$, $r\in(0,1]$, $B\in\cM_\alpha^{(r)}(A_0\cup A)\cap[\{m_1,\ldots,m_{k-1}\}\cup M_k]^{<\omega}$
 we have
 $$\big\|P_{(\max\supp(\tilde x_{A_0\cup A}),\infty)}\big(\tilde \Phi(A\cup B)\big)\big\|\ge  \frac{c}{2}\sum_{A_0\cup A\prec C\preceq A_0\cup A\cup B}\zeta(\alpha,C)$$
 (if $k=1$, $[\{m_1,\ldots,m_{k-1}\}]=\{\emptyset\}$).
\end{itemize}
If we assume that such a sequence $(M_k)_k$ has been chosen, it is straightforward to check that $N=\{m_k: k\in\N\}$ is the desired set. In the case $k=1$, for $A\in\cS_\alpha(A_0)\cap\{\emptyset\}$ we have $A=\emptyset$. Hence, if $A_0 = \emptyset$ then
for all $A\in\cS_\alpha(A_0)\cap\{\emptyset\}$ we have $A_0\cup A=\emptyset$, i.e. $l_1(A_0\cup A) = 0$ and hence, (c) and (d) are always satisfied. Choosing $M_1$ satisfying (b) completes the first inductive step. If, on the other hand,
$A_0$ is a singleton then for all $A\in\cS_\alpha(A_0)\cap\{\emptyset\}$ we have $A_0\cup A=A_0$, i.e. $l_1(A_0\cup A) > 0$. This means that condition (c) and (d) are reduced to the case in which $ A = \emptyset $. The choice of $M_1$ is done
in the same manner as in the general inductive step and we omit it.

Assume that we have chosen, for some $k \ge 1$, infinite sets
$M_k\subset M_{k-1}\subset \cdots \subset M_1\subset  M_0$, so that (b), (c), and (d) are satisfied for all $1\le k'\le k$. Observe that the inductive assumption implies that it is enough to choose $M_{k+1}\in [M_k]^{\omega}$
satisfying (b),  and the conditions  (c) and (d) for sets $A\in\cS_\alpha(A_0)\cap[\{m_1,\ldots,m_{k-1}\}]$ with $l_1(A_0\cup A) > 0$ and $\max(A)=m_{k}=\min(M_{k})$ (or  $A=\emptyset$,  in the case $k=1$, and $A_0$ is a singleton).
We set
\begin{align*}   
\delta &= \min\left\{ \zeta(\alpha,A_0\cup A):  A\in\cS_\alpha(A_0)\cap[\{m_1,m_2,\ldots,m_{k}\}]\text{ and }A_0\cup A\neq\emptyset\right\},\\
\vp &= \frac{\delta}{30} (\tilde c-c),\\
d &= \max\{\max\supp(\tilde x_{A_0\cup A}): A\in\cS_\alpha(A_0)\cap[\{m_1,\ldots,m_k\}]\},
\end{align*}
and choose a finite $\vp$-net $R$ of the interval $(0,1]$, with $1\in R$, which also has the property that for all $A\in  \cS_\alpha(A_0)\cap[\{m_1,m_2,\ldots,m_{k}\}]$ and all $j=0,1,2,\ldots ,l_1(A_0\cup A)$, 
\begin{equation}
\label{E:10.8.1}
j\cdot\zeta(\alpha,A_0\cup A)+\sum_{C\preceq A_0\cup A}\zeta(\alpha,C)\in R. 
\end{equation}
Fix a finite $\frac\vp2$-net $K$ of the unit ball of the finite dimensional space $\spa(E_j:1\kleq j\kleq~\!\!\!d)$.
For every $r\in R$ and $A\in\cS_\alpha(A_0)\cap[\{m_1,m_2,\ldots,m_{k}\}]$ we apply 
 Proposition \ref{P:3.1}  to $\cM_\alpha^{(r)}(A_0\cup A)\cap[ M_k]^{<\omega}$,  and find an infinite subset $\tilde M_{k+1}$ of $M_k$ so that for all $A\kin\cS_\alpha(A_0)\cap[\{m_1,m_2,\ldots,m_{k}\}]$ and $r\kin R$,
there exists $y_A^{(r)}$ in $K$ so that 
$$\big\|y_A^{(r)} - P_{[1,d]}\big(\tilde\Phi(A\cup B) \big)\big\| < \frac{\vp}{2}, \text{for all  $B$ in $\cM_\alpha^{(r)}(A_0\cup A)\cap[\tilde M_{k+1}]^{<\omega}$.}$$
In particular, note that for all $A\in\cS_\alpha(A_0)\cap[\{m_1,m_2,\ldots,m_{k}\}]$ and $r\in R$, for any $B_1$, $B_2$ in $\cM_\alpha^{(r)}(A_0\cup A)\cap[\tilde M_{k+1}]^{<\omega}$ we have
\begin{equation}
\label{E:10.8.2}
 \big\|P_{[1,d]}\big(\tilde\Phi(A\cup B_1) - \tilde\Phi(A\cup B_2)\big) \big\| < \vp.
\end{equation}
Using the property (iv)  of $l_1(\cdot)$  in Lemma \ref{L:4.6}, we can pass to an infinite subset $\widehat M_{k+1}$ of $\tilde M_{k+1}$, so that (b) is satisfied and moreover:
\begin{align}\label{E:10.8.3}\zeta(\alpha,A_0\cup A\cup B)<\vp,
 &\text{ if  $A\in\cS_\alpha(A_0)\cap[\{m_1,\ldots,m_k\}]$, and}\\
 &\text{$B\in\cS_\alpha(A_0\cup A)\cap[M_{k+1}]^{<\omega}$ with $\#B > l_1(A_0\cup A) > 0$}.\notag
 \end{align}

We will show that (c) is satisfied. To that end, fix $0<r\le 1$, $A\in\cS_\alpha(A_0)\cap[\{m_1,\ldots,m_k\}]$ with $\max(A_0\cup A) = m_k$ and $l_1(A_0\cup A) > 0$, and $B_1$, $B_2\in\cM_\alpha^{(r)}(A_0\cup A)\cap[\widehat M_{k+1}]^{<\omega}$ with $B_1 < B_2$.
If both sets $B_1$ and $B_2$ are empty, then (c) trivially holds, as the right-hand side of the inequality has to be zero. 
Otherwise, $B_s \neq\emptyset$, where $s=1$ or $s=2$. Note that $\max(A) = m_k$, i.e. $A_0\cup A\neq\emptyset$ and hence, since $l_1(A_0\cup A)>0$, putting $\tilde C = A_0\cup A\cup\{\min(B_s)\}$, by the definition of $\delta$
we obtain $\zeta(\alpha,\tilde C) = \zeta(\alpha, A_0\cup A)\ge \delta$. This easily yields:
\begin{align*}
d_{1,\alpha}(A_0\cup A\cup B_1, A_0\cup A\cup B_2) &= \sum_{A_0\cup A\prec C\preceq A_0\cup A\cup B_1}\zeta(\alpha,C) + \sum_{A_0\cup A\prec C\preceq A_0\cup A\cup B_2}\zeta(\alpha,C)\\
&\ge \delta = \frac{30 \vp }{\tilde c-c}
\end{align*}
Hence: 
\begin{equation}
\label{E:10.8.4}
 \vp \le \frac{\tilde c-c}{30} d_{1,\alpha}(A_0\cup A\cup B_1, A_0\cup A\cup B_2)
\end{equation}
Arguing similarly, we obtain:
\begin{equation}
\label{E:10.8.5}
r \ge \sum_{C\preceq A_0\cup A\cup B_s}\zeta(\alpha,C) \ge \sum_{C\preceq A_0\cup A}\zeta(\alpha,C) + \zeta(\alpha,\tilde C) \ge \delta > \min(R). 
\end{equation}
Choose $r_0$ to be the maximal element of $R$ with $r_0\le r$. Since $r_0\le r$,
we can find $\tilde B_1$ and $\tilde B_2$ in $\cM_\alpha^{(r_0)}(A_0\cap A)\cap[\widehat M_{k+1}]^{<\omega}$ so that $\tilde B_1\preceq B_1$ and $\tilde B_2\preceq B_2$.
We will show that
\begin{equation}
\label{E:10.8.6}
d_{1,\alpha}(A_0\cup A\cup B_1, A_0\cup A\cup\tilde B_1) < 3\vp\text{ and }d_{1,\alpha}(A_0\cup A\cup B_2, A_0\cup A\cup\tilde B_2) < 3\vp.
\end{equation}
We shall only show that this is the case for $B_1$, for $B_2$ the proof is identical. If $\tilde B_1 = B_1$, then there is nothing to prove and we may therefore assume that $\tilde B_1\prec B_1$.
Define $C_1 = \tilde B_1 \cup \{\min(B_1\setminus \tilde B_1)\}$, $r_1 = \sum_{C\preceq A_0\cup A\cup B_1}\zeta(\alpha,C)$, $\tilde r_1 = \sum_{C\preceq A_0\cup A\cup \tilde B_1}\zeta(\alpha,C)$,
and $r' = \sum_{C\preceq A_0\cup A\cup C_1}\zeta(\alpha, C)$.
The maximality of $\tilde B_1$ in $\cS_{\alpha}^{(r_0)}(A_0\cup A)$ implies
\begin{equation}
\label{E:10.8.7}
\tilde r_1 \le r_0 < r' \le r_1.
\end{equation}
We first assume that $\#C_1 \le l_1(A_0\cup A)$.
In this case however, by \eqref{E:10.8.1}, we obtain that $r' = \sum_{C\preceq A_0\cup A}\zeta(\alpha,C) + (\#C_1)\zeta(\alpha,A_0\cup A)$ is in $R$. This contradicts the maximality of $r_0$.
We conclude that $\#C_1 > l_1(A_0\cup A)$, which by \eqref{E:10.8.3} yields $r' - \tilde r_1 = \zeta(\alpha,A_0\cup A\cup C_1) < \vp$.
Hence, 
\begin{equation}\label{E:10.8.8} d_{1,\alpha}(A_0\cup A\cup B_1, A_0\cup A\cup\tilde B_1) = r_1 - \tilde r_1  = (r_1 - r') + (r' - \tilde r_1) < (r_1 - r_0) + \vp < 3\vp\end{equation}

We verify (c) now as follows
\begin{align*}
 \big\|P_{(d,\infty)}&\big(\tilde \Phi(A\cup B_1) - \tilde\Phi(A\cup B_2)\big)\big\| \\
 \ge&  \big\|P_{(d,\infty)}\big(\tilde \Phi(A\cup \tilde B_1) - \tilde\Phi(A\cup\tilde B_2)\big)\big\| \nonumber\\
 \qquad&-\left(\big\|P_{(d,\infty)}\big(\tilde \Phi(A\cup B_1) - \tilde\Phi(A\cup \tilde B_1)\big)\big\| + \big\|P_{(d,\infty)}\big(\tilde \Phi(A\cup B_2) - \tilde\Phi(A\cup \tilde B_2)\big)\big\|\right)\nonumber\\
 \ge& \big\|P_{(d,\infty)}\big(\tilde \Phi(A\cup \tilde B_1) - \tilde \Phi(A\cup\tilde B_2)\big)\big\|\notag\\
 \qquad&- \left(d_{1,\alpha}(A_0\cup A\cup B_1, A_0\cup A\cup\tilde B_1) + d_{1,\alpha}(A_0\cup A\cup B_2, A_0\cup A\cup\tilde B_2)\right)\phantom{aAAAaAA}\nonumber\\
 >& \big\|P_{(d,\infty)}\big(\tilde \Phi(A\cup \tilde B_1) - \tilde\Phi(A\cup\tilde B_2)\big)\big\| - 6\vp\notag\\
 \ge& \big\|\tilde \Phi(A\cup \tilde B_1) -\tilde \Phi(A\cup\tilde B_2) \big\|- \big\|P_{[1,d]}\big(\tilde \Phi(A\cup \tilde B_1) -\tilde \Phi(A\cup\tilde B_2)\big)\big\| - 6\vp\nonumber\\
 \ge& \big\|\tilde \Phi(A\cup \tilde B_1) -\tilde \Phi(A\cup\tilde B_2) \big\| - 7\vp\notag\\
 \ge& c(\tilde\Phi)d_{1,\alpha}(A_0\cup A\cup\tilde B_1, A_0\cup A\cup \tilde B_2) - 7\vp
 > \tilde cd_{1,\alpha}(A_0\cup A\cup B_1, A_0\cup A\cup B_2).\nonumber
\end{align*}
Here the second inequality follows from \eqref{E:10.2.1}, the third from \eqref{E:10.8.6},  the fifth from \eqref{E:10.8.2}, the sixth one from   \eqref{E:10.2.2}, and the last one from
\eqref{E:10.8.4}.

We will need to pass to a further subset of $\widehat M_{k+1}$ to obtain (d).
An application of the triangle inequality and  (c) (for $k+1$)  yield that 
 for every $A\in \cS_\alpha(A_0)\cap[\{m_1,m_2,\ldots m_k\}]$ with $l_1(A_0\cup A)$, $r\in(0,1]$ and $B_1, B_2\in \cM^{(r)}_\alpha(A_0\cup A)$, $B_1<B_2$, it follows that 
 $$\big\|P_{(\max\supp(\tilde x_{A_0\cup A}),\infty)}\big(\tilde \Phi(A\cup B_1)\big) \big\|\ge\frac{\tilde c}2 d_{1,\alpha} (A_0\cup A\cup B_1,A_0\cup A\cup B_2 ),$$
 or
 $$\big\|P_{(\max\supp(\tilde x_{A_0\cup A}),\infty)}\big( \tilde\Phi(A\cup B_2)\big) \big\|\ge \frac{\tilde c}2 d_{1,\alpha} (A_0\cup A\cup B_1,A_0\cup A\cup B_2 ).$$
We may pass therefore to a further infinite subset $M_{k+1}$ of $\widehat M_{k+1}$, so that
for any $A\in\cS_\alpha(A_0)\cap[\{m_1,m_2,\ldots,m_{k}\}]$, with $l_{1}(A_0\cup A)>0$ and $\max(A_0\cup A) = m_k$,  any $r\in R$, and for any $B$ in $\cM_\alpha^{(r)}(A_0\cup A)\cap[ M_{k+1}]^{<\omega}$ we have
\begin{equation}
\label{E:10.8.9}
\big\|P_{(d,\infty)}\tilde \Phi(A\cup B)\big\|\ge  \frac{\tilde c}{2} \sum_{A_0\cup A\prec C\preceq A_0\cup A\cup B}\zeta(\alpha,C).
\end{equation}
We will show that (d) is satisfied (meaning that \eqref{E:10.8.9} does not only hold if $r\in R$). Fix an arbitrary   $r\in (0,1]$, $A\in\cS_\alpha(A_0)\cap[\{m_1,\ldots,m_k\}]$, with $\max(A_0\cup A) = m_k$ and $l_1(A_0\cup A) > 0$, and $B\in\cM_\alpha^{(r)}(A_0\cup A)\cap[ M_{k+1}]^{<\omega}$.
If $B$ is empty, then the right-hand side of the inequality in (d) is zero and hence the conclusion holds. Arguing identically as in \eqref{E:10.8.4} we obtain
\begin{equation*}
 \vp \le   \frac{\tilde c-c}{30}\sum_{A_0\cup A\prec C\preceq A_0\cup A\cup B}\zeta(\alpha,C)
\end{equation*}
and arguing as in \eqref{E:10.8.5} we obtain $r > \min(R)$, so we may choose the maximal element $r_0$
of $R$ with $r_0\le r$, and $\tilde B$ the maximal element of $\cM^{(r_0)}_\alpha(A_0\cup A)\cap[M_{k+1}]^{<\omega}$ with $\tilde B\preceq B$. The argument yielding \eqref{E:10.8.6}, also yields
$$d_{1,\alpha}(A_0\cup A\cup B,A_0\cup A\cup \tilde B) = \sum_{A_0\cup A\cup \tilde B\prec C\preceq A_0\cup A\cup B}\zeta(\alpha,C) < 3\vp.$$
Arguing in a very similar manner as above:
\begin{align*}
 \big\|P_{(d,\infty)}\tilde \Phi(A\cup B)\big\|&\ge  \big\|P_{(d,\infty)}\tilde \Phi(A\cup \tilde B)\big\| - \big\|P_{(d,\infty)}\tilde \Phi(A\cup B) - P_{(d,\infty)}\tilde\Phi(A\cup\tilde B)\big\|\\
 & \ge\frac{\tilde c}{2} \sum_{A_0\cup A\prec C\preceq A_0\cup A\cup \tilde B}\zeta(\alpha,C) - 3\vp\\
 & \ge\frac{\tilde c}{2} \sum_{A_0\cup A\prec C\preceq A_0\cup A\cup B}\zeta(\alpha,C) - \left(3 + \frac{3\tilde c}{2}\right)\vp\\
 &\ge \frac{\tilde c}{2} \sum_{A_0\cup A\prec C\preceq A_0\cup A\cup B}\zeta(\alpha,C) - \left(\frac{\tilde c - c}{2}\right) \sum_{A_0\cup A\prec C\preceq A_0\cup A\cup B}\zeta(\alpha,C)\\
 &= \frac{c}{2} \sum_{A_0\cup A\prec C\preceq A_0\cup A\cup B}\zeta(\alpha,C).
\end{align*}
\end{proof}

\section{Some further observation on the  Schreier families}\label{S:11}

In this section $\beta$ will be a fixed  ordinal of the form $\beta=\omega^{\omega^\xi}$, with $1\le \xi<\omega_1$.

\subsection{Analysis of a maximal set $B$ in $\cS_{\beta\gamma}$}\label{SS:11.1}  Recall that by Proposition \ref{P:2.6},
for every $\gamma\le \beta$ there exists a sequence $\eta(\gamma,n)$ of ordinal numbers increasing to $\gamma$, so that $\lambda(\beta\gamma,n) = \beta\eta(\gamma,n)$ (recall that $\eta(\gamma,n)$ may also depend on $\beta$).

For every $\gamma\le \beta$ and $B\in\MAX(\cS_{\beta\gamma})$ we define a family of subsets of $B$, which we shall call the {\em  $\beta$-analysis of $B$,} and denote by  $\cA_{\beta,\gamma}(B)$. The definition is done recursively on $\gamma$.

\begin{subequations}
For $\gamma = 1$, set 
\begin{equation}
\label{E:11.1a}
\cA_{\beta,\gamma}(B) = \{B\}.
\end{equation}

Let $\gamma< \beta$ and assume that $\cA_{\beta,\gamma}(B)$ has been defined for all $B\in\MAX(\cS_{\beta\gamma})$.
For $B\in\MAX(\cS_{\beta(\gamma+1)})=\MAX(\cS_{\beta}[\cS_{\beta\gamma}])$ there are (uniquely defined) $B_1<\cdots<B_\ell$ in $\MAX(\cS_{\beta\gamma})$ with
$\{\min B_j:\;j=1,\ldots,\ell\}$ in $\MAX(\cS_\beta)$ so that
$B = \cup_{j=1}^\ell B_j$. Set
\begin{equation}
\label{E:11.1b}
\cA_{\beta,\gamma+1}(B) = \{B\}\cup\left(\cup_{j=1}^\ell\cA_{\beta,\gamma}(B_j)\right).
\end{equation}

Let $\gamma\le \beta$ be a limit ordinal and assume that $\cA_{\beta,\gamma'}(B)$ has been defined for all $\gamma'<\gamma$ and $B\in\MAX(\cS_{\beta\gamma'})$.
If now $B\in\MAX(\cS_{\beta(\gamma)})$, then $B\in\MAX(\cS_{\beta\eta(\gamma,\min (B))})$. Set
\begin{equation}
\label{E:11.1c}
\cA_{\beta,\gamma}(B) = \cA_{\beta,\eta(\gamma,\min (B))}(B).
\end{equation}
\end{subequations}

\begin{remark}\label{R:11.1}
Let $\gamma\le\beta$ and $B\in\MAX(\cS_{\beta\gamma})$. The following properties are straightforward consequences of the definition of $\cA_{\beta,\gamma}(B)$ and a transfinite induction.
\begin{itemize}
 
 \item[(i)] The set $\cA_{\beta,\gamma}(B)$ is a tree when endowed with ``$\supset$''.
 
 \item[(ii)] For $C$, $D$ in $\cA_{\beta,\gamma}(B)$ that are incomparable with respect to inclusion, we have either $C<D$, or $D<C$. 
 
 \item[(iii)] The minimal elements (with respect to inclusion) of $\cA_{\beta,\gamma}(B)$ are in $\cS_\beta$.
 
 \item[(iv)] If $D\in\cA_{\beta,\gamma}(B)$ is a non-minimal element, then its direct successors $(D_j)_{j=1}^\ell$ in $\cA_{\beta,\gamma}(B)$ can be enumerated so that $D_1<\cdots<D_\ell$ and $D = \cup_{j=1}^\ell D_j$. 
 
\end{itemize}

\end{remark}

\subsection{Components of a set $A$ in $\cS_{\beta\gamma}$.}\label{SS:11.2} We recursively define for all non-empty sets $A\in\cS_{\beta\gamma}$ and $\gamma\le \beta$,
a natural number $s(\beta,\gamma,A)$ and subsets $\Cp_{\beta,\gamma}(A,1),\ldots,\Cp_{\beta,\gamma}(A,s(\beta,\gamma,A))$ of $A$. We will call 
$(\Cp_{\beta,\gamma}(A,i))_{i=1}^{s(\beta,\gamma,A)}$ the {\em  components of $A$ in $\cS_{\beta\gamma}$,
with respect to $\cS_\beta$.}
\begin{subequations}

If $\gamma = 1$, i.e. $A$ is a non-empty set in $\cS_\beta$, define
\begin{equation}
\label{E:11.2a}
s(\beta,\gamma,A) = 1 \;\text{and}\; \Cp_{\beta,\gamma}(A,1) = A. 
\end{equation}

Let $\gamma\le \beta$ and assume that $(\Cp_{\beta,\gamma}(A,i))_{i=1}^{s(\beta,\gamma,A)}$ has been defined for all non-empty sets $A$ in $\cS_{\beta\gamma}$. If now $A$ is a non-empty set in $\cS_{\beta(\gamma + 1)} = \cS_{\beta}[\cS_{\beta\gamma}]$,
then there are non-empty sets $A_1<A_2<\cdots<A_d$ in $\cS_{\beta\gamma}$ so that:
\begin{itemize}

\item[(i)] $A = \cup_{i=1}^{d}A_i$,

\item[(ii)] $\{\min A_i: i=1,\ldots,d\}$ is in $\cS_\beta$, and

\item[(iii)] the sets $A_1,\ldots,A_{d-1}$ are in $\MAX(\cS_{\beta\gamma})$.

\end{itemize}
Note that the set $A_d$ may or may not be in $\MAX(\cS_{\beta\gamma})$. It may also be the case that $d = 1$, which in particular happens when $A\in\cS_{\beta\gamma}$. Set $\bar A = A_d$, which is always non-empty, and we define
\begin{equation}
\label{E:11.2b}
\begin{split}
&s(\beta,\gamma+1,A) = s(\beta,\gamma,\bar A) + 1,\\& \;
\Cp_{\beta,\gamma+1}(A,1) = \cup_{i<d}A_i,\;\text{and }\\&\:  \Cp_{\beta,\gamma+1}(A,i) =  \Cp_{\beta,\gamma}(\bar A,i\kminus1),\;\text{if}\; 2\kleq i\kleq s(\beta,\gamma\kplus1,A).
\end{split}
\end{equation}
Note that in the case $d = 1$, $\Cp_{\beta,\gamma+1}(A,1)$ is the empty set.

Let $\gamma\le \beta$ be a limit ordinal and assume that $(\Cp_{\beta,\gamma'}(A,i))_{i=1}^{s(\beta,\gamma',A)}$ has been defined for all $\gamma'<\gamma$ and non-empty sets $A$ in $\cS_{\beta\gamma'}$.
If$A$ is a non-empty set in $\cS_{\beta\gamma}$,
then $A\in\cS_{\beta\eta(\gamma,\min (A))}$ and we define
\begin{equation}
\label{E:11.2c}\begin{split}
&s(\beta,\gamma,A) = s(\beta,\eta(\gamma,\min(A)),A)\;\text{and}\\
&\;\Cp_{\beta,\gamma}(A,i))=   \Cp_{\beta,\eta(\gamma,\min(A))}(A,i)) \text{ for $i=1,2,\ldots, s(\beta,\gamma,A)$}.
\end{split}
\end{equation}

\begin{rem}
Let $\gamma\le\beta$ and $A\in\cS_{\beta\gamma}\setminus\{\emptyset\}$. The following properties follow easily from the definition of $(\Cp_{\beta,\gamma}(A,i))_{i=1}^{s(\beta,\gamma,A)}$ and a transfinite induction on $\gamma$.
\begin{itemize}
 
 \item[(i)] $A = \cup_{i=1}^{s(\beta,\gamma,A)}\Cp_{\beta,\gamma}(A,i)$.
 
 \item[(ii)] For $1\le i < j\le s(\beta,\gamma,A)$ so that both $\Cp_{\beta,\gamma}(A,i)$ and $\Cp_{\beta,\gamma}(A,j)$ are non-empty, we have $\Cp_{\beta,\gamma}(A,i) < \Cp_{\beta,\gamma}(A,j)$.
 \item[(iii)]  $\Cp(A,s(\beta,\gamma(A))\not=\emptyset$.
\end{itemize}

\end{rem}

\end{subequations}

\begin{lem}
\label{L:11.2}
Let $\xi$ and $\gamma$ be countable ordinal numbers with $\gamma\le \beta = \omega^{\omega^\xi}$. Let also $B$ be a set in $\MAX(\cS_{\beta\gamma})$ and $\emptyset\prec A\preceq B$. If  $\cA_{\beta,\gamma}(B)$ is the $\beta$-analysis of $B$
and $(\Cp_{\beta,\gamma}(A,i))_{i=1}^{s(\beta,\gamma,A)}$ are the components of $A$ in $\cS_{\beta\gamma}$ with respect to $\cS_\beta$, then there exists a maximal chain
$B = D_1(A) \supsetneq D_2(A) \supsetneq \cdots \supsetneq D_{s(\beta,\gamma,A)}(A)$ in $\cA_{\beta,\gamma}(B)$ satisfying the following:
\begin{itemize}

\item[(i)] $\Cp_{\beta,\gamma}(A,i)\preceq D_i(A)$ for $1\le i \le s(\beta,\gamma,A)$, and 

\item[(ii)] if $1\le i < s(\beta,\gamma, A)$ then $\Cp_{\beta,\gamma}(A,i)\subsetneq D_i(A)$.
\end{itemize}
\end{lem}
\begin{proof}
We  prove the statement by transfinite induction on $1\le \gamma\le \beta$. If $\gamma = 1$, then $\cA_{\beta,\gamma}(B) = \{B\}$ and $(\Cp_{\beta,\gamma}(A,i))_{i=1}^{s(\beta,\gamma,A)} = \{A\}$, and our claim follows trivially.

Let $\gamma<\beta$ and assume that the statement holds for all non-empty $A\in\cS_{\beta\gamma}$ and $B\in\MAX(\cS_{\beta\gamma})$ with $A\preceq B$. Let now $A$ be a non-empty set in $\cS_{\beta(\gamma +1)}$ and $B\in\MAX(\cS_{\beta(\gamma+1)})$.
If $B = \cup_{i=1}^\ell B_j$, where $B_1<\cdots<B_\ell$ are in $\MAX(\cS_{\beta\gamma})$ with $\{\min(B_j): 1\le j\le \ell\}\in\MAX(\cS_{\beta})$, then by \eqref{E:11.1b} we obtain
$\cA_{\beta,\gamma+1}(B) = \{B\}\cup(\cup_{j=1}^\ell\cA_{\beta,\gamma}(B_j))$. Define 
\begin{equation*}
d= \max\{1\le j\le \ell: B_j\cap A\neq\emptyset\},\,
A_i = B_i\;\text{for}\;1\le i < d,\text{ and }
A_d  = A\cap B_d.
\end{equation*}
Then, by \eqref{E:11.2b}, letting $\bar A = A_d$, we obtain $s(\beta,\gamma+1,A) = s(\beta,\gamma,\bar A) + 1$, $\Cp_{\beta,\gamma+1}(A,1) = \cup_{i<d}A_i$, and $\Cp_{\beta,\gamma+1}(A,i) =  \Cp_{\beta,\gamma}(\bar A,i-1)$,
for $2\le i\le s(\beta,\gamma+1,A)$.
Apply the inductive assumption to $A_d = \bar A\preceq B_d$, to find a maximal chain $B_d = D_1(\bar A) \supsetneq \cdots\supsetneq D_{s(\beta,\gamma,\bar A)}(\bar A)$ in $\cA_{\beta\gamma}(B_d)$ satisfying (i) and (ii)
with respect to $( \Cp_{\beta,\gamma}(\bar A,i))_{i=1}^{s(\beta,\gamma,\bar A)}$. Define
\begin{equation}
D_1(A) = B,\text{ and }
D_i(A) = D_{i-1}(\bar A),\;\text{for}\;2\le i\le s(\beta,\gamma,A).
\end{equation}
Clearly, $(D_i(A))_{i=1}^{s(\beta,\gamma,\bar A)}$ is a maximal chain in $\cA_{\beta,\gamma+1}(B)$. It remains to verify that (i) and (ii) are satisfied, with respect to $(\Cp_{\beta,\gamma+1}(A,i))_{i=1}^{s(\beta,\gamma+1,A)}$.
Assertions (i) and (ii), in the case $i=1$, follow trivially from $\Cp_{\beta,\gamma+1}(A,1) = \cup_{j<d}A_j = \cup_{j<d}B_j \prec \cup_{j<d}B_j \preceq B$.
Assertion (i) and (ii) in the case $i\neq 1$, follow easily  from the inductive assumption and $\Cp_{\beta,\gamma+1}(A,i) = \Cp_{\beta,\gamma}(A',i-1)$, for $2\le i\le s(\beta,\gamma+1,A)$. 
%We now show (iii) in the case $1<i<s(\beta,\gamma+1,A)$.
%Note that $D_i(A)\cap A = D_{i-1}(\bar A)\cap A = (B_d\cap D_{i-1}(\bar A))\cap A = D_{i-1}(\bar A)\cap(B_d\cap A) = D_{i-1}(\bar A)\cap \bar A$. Then $1\le i-1<s(\beta, 
%\gamma+1,A) - 1 = s(\beta,\gamma,\bar A)$ and the result
%follows from the inductive assumption.
%It only remains to observe that (iii) holds for $i=1$. Indeed, the immediate successors $E_1<\cdots<E_d$ of $D_1(A) = B$ are $B_1 <\cdots<B_\ell$ and in the case $A(1)$ is % not empty, we conclude that $t = d-1$ which yields the desired result.

To conclude the proof, the case in which $\gamma\le \beta$ is a limit ordinal number so that the conclusion is satisfied for all $\gamma'<\gamma$, we just observe that the result is an immediate consequence of
\eqref{E:11.1c} and \eqref{E:11.2c}.
\end{proof}

For the next result recall the definition of the doubly indexed fine Schreier families $\cF_{\beta,\gamma}$ introduced in Subsection \ref{SS:2.3}.
\begin{lem}
\label{L:11.3}Let  $\gamma\le\beta$. Then for all  $A\in\cS_{\beta\gamma}$, with $\Cp_{\beta\gamma}(A,i)\neq\emptyset$ for $1\le i\le s(\beta,\gamma,A)$, we have
\begin{equation}\label{E:11.3.1}
\left\{\min(\Cp_{\beta,\gamma}(A,i)): 1\le i\le s(\beta,\gamma,A)\right\}\in\MAX(\cF_{\beta,\gamma}).
\end{equation}
\end{lem}

\begin{proof}
We prove this statement by induction on $\gamma$. If $\gamma = 1$ and $A\in\cS_\beta$ satisfying the assumptions of this Lemma, then $A = \Cp_{\beta,1}(A,1)\neq\emptyset$ and hence the result easily follows from $\MAX(\cF_{\beta,1}) = \{\{n\}: n\in\N\}$.

 Assume  that the result holds for some $\gamma < \beta$ and let $A\in\cS_{\beta(\gamma+1)}$,
with $\Cp_{\beta\gamma}(A,i)\neq\emptyset$ for $1\le i\le s(\beta,\gamma+1,A)$. By the inductive assumption and \eqref{E:11.2b} we obtain that
$B =  \{\min(\Cp_{\beta,\gamma}(A,i)): 2\le i\le s(\beta,\gamma+1,A)\}\in\MAX(\cF_{\beta,\gamma})$. We claim that 
$\tilde A = \big\{\min\big(\Cp_{\beta,\gamma+1}(A,1)\big)\big\}\cup B\in\MAX(\cF_{\beta,\gamma+1})$. Indeed, if this is not the case, 
then by the spreading property of $\cF_{\beta,\gamma+1}$ there is $C\in\cF_{\beta,\gamma+1}$ with $\tilde A\prec C$. Then, 
$B\prec \tilde C = C\setminus\big\{\min\big(\Cp_{\beta,\gamma+1}(A,1)\big)\big\}$. It follows by
$\cF_{\beta,\gamma+1} = \cF_{\beta,1}\sqcup_<\cF_{\beta,\gamma}$, that $\tilde C\in\cF_{\beta,\gamma}$. The maximality of $B$ yields a contradiction.

Assume now that $\gamma\le\beta$ is a limit ordinal number so that the conclusion holds for all $\tilde \gamma<\gamma$. Let $A\in\cS_{\beta\gamma}$ so that $\Cp_{\beta,\gamma}(A,i)\neq\emptyset$ for $1\le i\le s(\beta,\gamma,A)$.
Note that $A\in\cS_{\beta\eta(\gamma,\min(A))}$. and by \eqref{E:11.2c} and the inductive assumption we have $\tilde A = \{\min(\Cp_{\beta,\gamma}(A,i)): 1\le i\le s(\beta,\gamma,A)\}\in\MAX(\cF_{\beta,\eta(\gamma,\min(A))})$. By
\eqref{E:2.12.1} we obtain $\tilde A \in\MAX(\cF_{\beta,\gamma})$.
\end{proof}

For  $\gamma\le \beta$ and a set $B$ in $\MAX(\cS_{\beta\gamma})$ we define
\begin{equation}\label{E:11.3}
\cE_{\beta,\gamma}(B) = \{\emptyset\prec A\preceq B:\;\Cp_{\beta,\gamma}(A,i)\neq\emptyset\;\text{for}\;1\le i\le s(\beta,\gamma,A)\}.
\end{equation} 

\begin{lem}
\label{L:11.4}
 Let  $\gamma < \beta$. If $B$ is in $\MAX(\cS_{\beta(\gamma+1)}) = \cS_\beta[\cS_{\beta\gamma}]$ and $B = \cup_{j=1}^\ell B_j$, where $B_1<\cdots<B_\ell$ are in $\MAX(\cS_{\beta\gamma})$ and $\{\min(B_j): 1\le j\le \ell\}\in\MAX(\cS_{\beta})$, then
\begin{align*}
\big\{A: A\preceq B\text{ and } &A\not\in\cE_{\beta(\gamma+1)}(B)\big\} = \\ 
&\{A: A\preceq B_1\}\cup\left(\bigcup_{m=2}^{\ell}\Big\{\Big(\bigcup_{j=1}^{m-1}B_j\Big)\cup A: A\preceq B_m\text{ and } A\not\in\cE_{\beta,\gamma}(B_{m})\Big\}\right).
\end{align*} 
\end{lem} 

\begin{proof}
Let $\emptyset\not=D\preceq B$. Define $m = \max\{1\le j\le \ell: D\cap B_j\neq\emptyset\}$ and $A = B_{m}\cap D$. Note that $D = (\cup_{j<m}B_j)\cup A$, where $\cup_{j<m}B_j = \emptyset$ if $m = 1$. By \eqref{E:11.2b} we obtain, $s(\beta,\gamma+1,D) = s(\beta,\gamma,A) + 1$, $\Cp_{\beta,\gamma+1}(D,1) = \cup_{j<m}B_j$ and $\Cp_{\beta,\gamma+1}(D,i) = \Cp_{\beta,\gamma}(A,i-1)$ for $2\le i\le s(\beta,\gamma+1,D)$.

Observe that $\Cp_{\beta,\gamma+1}(D,1)=\emptyset$ if and only if $m=1$, i.e. $A = D\preceq B_1$. On the other hand, if $\Cp_{\beta,\gamma+1}(D,1)\neq\emptyset$, then for some $2\le i\le s(\beta,\gamma+1,D)$, we have $\Cp_{\beta,\gamma+1}(D,i) = \emptyset$, if and only if $\Cp_{\beta,\gamma}(A,i-1)=\emptyset$. These observations yield our claim.
\end{proof}

\begin{rem}
Under the assumptions of Lemma \ref{L:11.4},
if for $1\le j\le \ell-1$ we define
 $$\cE_{\beta,\gamma+1}^{(j)}(B) =  \{A\in\cE_{\beta,\gamma+1}(B): \Cp_{\beta,\gamma}(A,1) = \cup_{i\le j}B_i\},$$
then, using a similar argument as the one used in the proof of Lemma \ref{L:11.4}, we obtain
$$\cE_{\beta,\gamma+1}^{(j)}(B) = \left\{(\cup_{i=1}^jB_i)\cup C: C\in\cE_{\beta,\gamma}(B_{j+1})\right\}\text{ and }\cE_{\beta(\gamma+1)}(B) = \cup_{j=1}^{\ell-1}\cE_{\beta,\gamma+1}^{(j)}.$$
\end{rem}

\begin{rem}
Using the fact $\cS_1\subset\cS_\alpha$ that for all countable ordinal numbers $\alpha$, and that $\MAX(\cS_1) = \{F\subset\N: \min(F) = \#F\}$, it is easy to verify that for all $F\in\MAX(\cS_\alpha)$ we have $\max(F) \ge 2\min(F) - 1$.
In particular, if $B_1< B_2$ are both in $\MAX(\cS_{\alpha})$, then
\begin{equation}
\label{E:11.4}
2\min(B_1) \le \min(B_2).
\end{equation}

\end{rem}

\begin{lem}\label{L:11.5}
Let  $\gamma\le \beta $ and let $B$ be a set in $\MAX(\cS_{\beta\gamma})$. Then
\begin{equation}
\label{E:11.5.1}
\sum_{\substack{A\preceq B\\A\not\in\cE_{\beta,\gamma}(B)}}\zeta{(\beta\gamma,A)} < \frac{2}{\min(B)}.
\end{equation}
\end{lem}

\begin{proof}
We  prove the statement by transfinite induction for all $1\le \gamma\le \beta$. If $\gamma = 1$, then the complement of $\cE_{\beta,1}(B)$ only contains the empty set and the result trivially holds.

Let $\gamma<\beta$ and assume that the statement holds for all $B\in\MAX(\cS_{\beta\gamma})$. Let $B\in\MAX(\cS_{\beta(\gamma+1)})$.
Let $B_1<\cdots<B_\ell$ in $\MAX(\cS_{\beta\gamma})$ so that $\{\min(B_j): 1\le j\le \ell\}\in\MAX(\cS_{\beta})$ and $B = \cup_{j=1}^\ell B_j$. For $m=1,\ldots,\ell$, define $D_m = \{\min(B_j): 1\le j\le m\}$.
Proposition \ref{P:4.4} implies the following: if $C\preceq B$, $m = \max\{1\le j\le \ell: C\cap B_j\neq\emptyset\}$ and $A = C\cap B_m$, then $\zeta(\beta(\gamma+1),C)  = \zeta(\beta,D_m)\zeta(\beta\gamma,A)$.
We combine this fact with Lemma \ref{L:11.4}, \eqref{E:11.4}, and \eqref{E:11.5}  to obtain the following:
\begin{eqnarray*}
\sum_{\substack{A\preceq B\\A\not\in\cE_{\beta,\gamma}(B)}}\zeta{(\beta\gamma,A)} &=&
\zeta(\beta, D_1)\sum_{A\preceq B_1}\zeta(\beta\gamma, A) + \sum_{j=2}^\ell\zeta(\beta,D_j)\sum_{\left\{\substack{A\preceq B_j\\A\not\in\cE_{\beta\gamma}(B_j)}\right\}}\zeta(\beta\gamma,A)\\
&< & \zeta(\beta, D_1) + \sum_{j=2}^\ell\zeta(\beta,D_j)\frac{2}{\min(B_j)}\\ 
&\le & \frac{1}{\min(B_1)} + \sum_{j=2}^\ell\zeta(\beta,D_j)\frac{2}{\min(B_j)}\\ 
&\le & \frac{1}{\min(B_1)} + \sum_{j=2}^\ell\zeta(\beta,D_j)\frac{1}{\min(B_1)} \le \frac{2}{\min(B_1)} = \frac{2}{\min(B)}.
\end{eqnarray*}

If $\gamma\le \beta$ is a limit ordinal number so that the conclusion is satisfied for all $\gamma'<\gamma$, we just observe that the result is an immediate consequence of \eqref{E:11.2c} and \eqref{E:4.3}.
\end{proof}

\begin{lem}
\label{being in the same component is backwards closed}
Let  $\gamma\le \beta $ and $B\in\MAX(\cS_{\beta\gamma})$. 
If $A^{(1)}$, $A^{(2)}\in\cE_{\beta,\gamma}(B)$, $(D_{k}(A^{(1)}))_{k=1}^{s(\beta,\gamma,A^{(1)})}$ and
$(D_{k}(A^{(2)}))_{i=1}^{s(\beta,\gamma,A^{(2)})}$ are the maximal chains of $\cA_{\beta,\gamma}(B)$
provided by Lemma \ref{L:11.2}, and $1\le i\le \min\{s(\beta,\gamma,A^{(1)}),s(\beta,\gamma,A^{(2)})\}$ is such that
 $D_i(A^{(1)}) = D_i(A^{(2)})$, then  we have $D_j(A^{(1)}) = D_j(A^{(2)})$ for all $1\le j<i$.
\end{lem}

\begin{proof}
As $(D_{k}(A^{(1)}))_{k=1}^{s(\beta,\gamma,A)}$ and $(D_{k}(A^{(2)}))_{i=1}^{s(\beta,\gamma,A^{(2)})}$ are both maximal chains of $\cA_{\beta,\gamma}(B)$
with $D_i(A^{(1)}) = D_i(A^{(2)})$, the result follows from Remark \ref{R:11.1} (i).
\end{proof}

\begin{lem}
\label{being in the same component is the same as having common minimum}
Let  $\gamma\le \beta$ and $B\in\MAX(\cS_{\beta\gamma})$. If $A^{(1)}$, $A^{(2)}\in\cE_{\beta,\gamma}(B)$, and $1\le i\le \min\{s(\beta,\gamma,A^{(1)}),s(\beta,\gamma,A^{(2)})\}$,
then $D_i(A^{(1)}) = D_i(A^{(2)})$ if and only if $\min(\Cp_{\beta,\gamma}(A^{(1)},i)) = \min(\Cp_{\beta,\gamma}(A^{(2)},i))$, where $(D_{k}(A^{(1)}))_{k=1}^{s(\beta,\gamma,A)}$ and
$(D_{k}(A^{(2)}))_{i=1}^{s(\beta,\gamma,A^{(2)})}$ are the maximal chains of $\cA_{\beta,\gamma}(B)$
provided by Lemma \ref{L:11.2}.
\end{lem}
\begin{proof}
Assume that  $D_i(A^{(1)}) \keq D_i(A^{(2)})$. Lemma \ref{L:11.2} (i) and the assumptions $\Cp(A^{(1)},i)\neq\emptyset$,  and $\Cp(A^{(2)},i)\neq\emptyset$ yield $\min(\Cp_{\beta,\gamma}(A^{(1)},i)) = \min(\Cp_{\beta,\gamma}(A^{(2)},i))$.
For the converse let $A^{(1)}, A^{(2)}\kin\cE_{\beta,\gamma}(B)$ with $\min(\Cp_{\beta,\gamma}(A^{(1)},i))\keq\min(\Cp_{\beta,\gamma}(A^{(2)},i))$. Since all elements of $\cA_{\beta,\gamma}(B)$ either compare with respect to ``$\subset$'', or are disjoint,
Lemma \ref{L:11.2}~(i) and $\min(\Cp_{\beta,\gamma}(A^{(1)},i)) \keq \min(\Cp_{\beta,\gamma}(A^{(2)},i))$ imply that either $D_i(A^{(1)})\ksubset D_i(A^{(2)})$, or $D_i(A^{(2)})\ksubset D_i(A^{(1)})$. We  assume the first, and
towards a contradiction assume that $D_i(A^{(1)})\subsetneq D_i(A^{(2)})$. The maximality of $(D_j(A^{(2)}))_{j=1}^{s(\beta,\gamma,A^{(2)})}$ in $\cA_{\beta,\gamma}(B)$ implies that there is $1\kleq j \kle i$, so that $D_j(A^{(2)}) \keq D_i(A^{(1)})$.
As $A^{(2)}_{\beta,\gamma}(j)\!\neq\!\emptyset$, we obtain by Lemma \ref{L:11.2} (i) $\min(\Cp_{\beta,\gamma}(A^{(1)},i)) \keq \min(D_i(A^{(1)})) \keq \min(D_j(A^{(2)})) \keq \min(\Cp_{\beta,\gamma}(A^{(2)},j))\kle\min(\Cp_{\beta,\gamma}(A^{(2)},i))$ which is a contradiction.
\end{proof}

\subsection{Special families of  convex combinations}\label{SS:11.3}
\begin{defin}
\label{special families of convex combinations}
Let $\gamma\le \beta$ and $B\in\MAX(\cS_{\beta\gamma})$.
A family of non-negative numbers $\{r(A,k): A\in\cE_{\beta,\gamma}(B), 1\le k\le s(\beta,\gamma,A)\}$, is called a {\em $(\beta,\gamma)$-special family of convex combinations for $B$} if the following are satisfied:
\begin{itemize}

\item[(i)] $\sum_{k=1}^{s(\beta,\gamma,A)}r(A,k) = 1$ for all $A\in\cE_{\beta,\gamma}(B)$ and

\item[(ii)] If $A^{(1)}$, $A^{(2)}$ are both in $\cE_{\beta,\gamma}(B)$, $(D_{k}(A^{(1)}))_{k=1}^{s(\beta,\gamma,A)}$ and $(D_{k}(A^{(2)}))_{i=1}^{s(\beta,\gamma,A^{(2)})}$ are the maximal chains in $\cA_{\beta,\gamma}(B)$
provided by Lemma \ref{L:11.2}, and for some $k$ we have $D_k(A^{(1)}) = D_k(A^{(2)})$, then $r(A^{(1)},k) = r(A^{(2)},k)$.

\end{itemize}
\end{defin}

\begin{rem} Let $\gamma\kleq \beta$ and  $B\kin \MAX(\cS_{\beta\gamma})$ and let $\{r(A,k)\!:\!A\kin\cE_{\beta,\gamma}(B), k=1,2,\ldots ,s(\beta,\gamma,A)\}$ be a family of $(\beta,\gamma)$-special convex combinations for $B$.
For $A\in \cE_{\beta,\gamma}(B)$, let $(D_{k}(A)_{i=1}^{s(\beta,\gamma,A)}$ be the maximal chain in $\cA_{\beta,\gamma}(B)$
provided by Lemma \ref{L:11.2}, and let $(A(i))_{i=1}^{s(\beta,\gamma,A)}$ be the components of $A$ in $\cS_{\beta\gamma}$. 

By construction, $D_1(A)=B$ for all $A \in \cE_{\beta,\gamma}(B)$, and thus $r(A,1)$ does not depend on $A$.
Secondly $D_2(A)$, only depends on $A(1)$, thus $r(A^{(1)},2)=r(A^{(2)},2)$ if $\Cp_{\beta,\gamma}(A^{(1)},1)=\Cp_{\beta,\gamma}(A^{(2)},1)$, for 
any $A^{(1)},A^{(2)}\in  \cE_{\beta,\gamma}(B)$. We can continue and inductively we observe  that for all $k\le \min(s(\beta,\gamma,A^{(1)}),s(\beta,\gamma,A^{(2)}))$
if  $\Cp_{\beta,\gamma}(A^{(1)},i)=\Cp_{\beta,\gamma}(A^{(2)},i)$, for all $i=1,2\ldots, k-1$, then 
 $r(A^{(1)},k)=r(A^{(2)},k)$.

Let  $\gamma\kleq \beta$ and $B\kin\MAX(\cS_{\beta\gamma})$. If $\gamma$ is a limit ordinal number, then $B\kin\MAX(\cS_{\beta\eta(\gamma,\min(B))})$
and any $(\beta,\gamma)$-special family of convex combinations $\{r(A,k): A\in\cE_{\beta,\gamma}(B), 1\le k\le s(\beta,\gamma,A)\}$ is also a $(\beta,\eta(\gamma,\min(B)))$-special family of convex combinations,
as $\cE_{\beta,\gamma}(B) = \cE_{\beta,\eta(\gamma,\min(B))}(B)$ and for $A\in\cE_{\beta,\eta(\gamma,\min(B))}(B)$ we have $s(\beta,\eta(\gamma,\min(B)),A) = s(\beta,\gamma,A)$.
\end{rem}

\begin{lem}
\label{special families demoted}
We are given   $\gamma\kle \beta $, $B\in\MAX(\cS_{\beta(\gamma+1)})$, and  a $(\beta,\gamma+1)$-special family of convex combinations $\{r(A,k): A\kin\cE_{\beta,\gamma+1}(B), 1\kleq k\kleq s(\beta,\gamma\kplus1,A)\}$.
  Assume that for some $D\in\cE_{\beta,\gamma+1}(B)$ (and hence for all of them) we have $r(D,1) < 1$. Let $B = \cup_{j=1}^dB_j$, where $B_1<\cdots <B_d$ are the immediate predecessors of
$B$ in $\cA_{\beta\gamma}(B)$. For every $1\le j< d$ consider the family $\{r^{(j)}(C,k): C\in\cE_{\beta,\gamma}(B_{j+1}) \}$, with $$r^{(j)}(C,k) = \frac{1}{1-r(D,1)}r\left(\left(\cup_{i=1}^{j}B_i\right)\cup C,k+1\right)$$
for $k=1,\ldots,s(\beta,\gamma,C) = s(\beta,\gamma,(\cup_{i=1}^{j}B_i)\cup C) - 1$. Then $\{r^{(j)}(C,k): C\in\cE_{\beta,\gamma}(B_{j+1}) \}$ is a $(\beta,\gamma)$-special family of convex combinations.
\end{lem}

\begin{proof}
By \eqref{E:11.2b}, if $C\in\cE_{\beta,\gamma}(B_{j+1}) $ then  $A = (\cup_{i=1}^jB_i)\cup C\in \cE_{\beta,\gamma+1}(B)$  and  $s(\beta,\gamma+1,A) = s(\beta,\gamma,C)+1$, which implies that Definition \ref{special families of convex combinations} (i) is satisfied.
To see that (ii) holds, let $C^{(1)}$, $C^{(2)}$ be in $\cE_{\beta,\gamma}(B_{j+1}) $ so that for some $k$ we have $D_k(C^{(1)}) = D_k(C^{(2)})$. Then by Lemma \ref{being in the same component is the same as having common minimum} we have
$\min(\Cp_{\beta,\gamma} (C^{(1)}, k)=\min( \Cp_{\beta,\gamma} (C^{(2)}, k )$. Setting $A^{(1)} = (\cup_{i=1}^jB_t)\cup C^{(1)}$ and $A^{(2)} = (\cup_{i=1}^jB_t)\cup C^{(2)}$, by \eqref{E:11.2b} we obtain 
 $\Cp_{\beta,\gamma+1}(A^{(1)},k+1) = \Cp_{\beta,\gamma}(C^{(1)},(k)$ and
$\Cp_{\beta,\gamma}(A^{(2)},k+1) =\Cp_{\beta,\gamma}(C^{(2)}, k)$, i.e. $\min(\Cp_{\beta,\gamma}(A^{(1)}, k+1))=\min( \Cp_{\beta,\gamma}(A^{(2)}, k+1))$. By Lemma \ref{being in the same component is the same as having common minimum} we obtain $D_{k+1}(A^{(1)}) = D_{k+1}(A^{(2)})$ and
therefore $r(A^{(1)},k+1) = r(A^{(2)},k+1)$, which yields that $r^{(j)}(C^{(1)},k) = r^{(j)}(C^{(2)},k)$.
\end{proof}

\section{Conclusion of the  proof of Theorems A and C}\label{S:12}

Again, we fix $\xi<\omega_1$ and put $\beta=\omega^{\omega^\xi}$. We secondly assume that $X$ is a Banach space $X$ with a bimonotone FDD $(F_j)$. 
  By the main result in \cite{Sch} every reflexive  Banach space $X$ embeds into a reflexive  Banach space $Z$ with basis, so  that $\Sz(Z)=\Sz(X)$ and  $\Sz(Z^*)=\Sz(X^*)$.
The coordinate projections on finitely or cofinitely many
coordinates are denoted by $P_A$ (see Section \ref{S:10} after Definition \ref{D:10.5}).

\begin{defin}\label{D:12.1}
Let  $\gamma\le \beta $, $M\in[N]^\omega$ and $A_0$ be a subset of $\N$ that is either empty or a singleton.
Let also $\Phi:\cS_{\beta\gamma}(A_0)\cap[M]^{<\omega}\rightarrow X$ be a semi-embedding of $\cS_{\beta\gamma}\cap[M]^{<\omega}$ into $X$, starting after $A_0$, that is $c$-refined, for some $0<c\le 1$. Let
$\{x_\emptyset\}\cup\{x_{A_0\cup A}: A\in\cS_{\beta\gamma}(A_0)\cap[M]^{<\omega}\}$ be the family generating $\Phi$ (recall that notation from the Remark after Definition \ref{D:10.2}).

Let $E\in\MAX(\cS_{\beta\gamma}(A_0)\cap[M]^{<\omega})$. For $A\preceq E$ recall the definition of $s(\beta,\gamma, A)$ and of $(\Cp_{\beta,\gamma}(A,i))_{i=1}^{s(\beta,\gamma, A)}$.
Recall also  from \eqref{E:11.3}
$\cE_{\beta,\gamma}(A_0\cup E)= \big\{\emptyset\prec A\preceq A_0\cup E: \Cp_{\beta,\gamma}(A,i)\not= \emptyset, \text{ for }i=1,2\ldots, s(\beta,\gamma,A)\big\}.$

For each $A\in \cE_{\beta,\gamma}(A_0\cup E)$ we will write $x_A$ as a sum of a block sequence
\begin{equation}\label{E:12.1.1} x_A=\sum_{k=1}^{s(\beta,\gamma,A)} x_{\Phi,A}^{(k)},\end{equation}
with $x_{\Phi,A}^{(k)}=P_{I(A,k)}(x_A)$,  for $k\keq 1,2,\ldots,s(\beta,\gamma,A)$, where
$I(A,1) = [1,\max\supp(x_{\Cp_{\beta,\gamma}(A,1)})]$ and $I(A,k) = (\max\supp(x_{\cup_{i=1}^{k-1}\Cp_{\beta,\gamma}(A,i)}),\max\supp(x_{\cup_{i=1}^{k}\Cp_{\beta,\gamma}(A,i)})]$ for $1\kle k\kle s(\beta,\gamma,A)$.
We call  $((x_{\Phi,A}^{(k)})_{k=1}^{s(\beta,\gamma,A)})_{A\in\cE_{\beta,\gamma}(A_0\cup E)}$ the {\em block step decomposition of $E$ with respect to $\Phi$}.
\end{defin}
\begin{rem}
Let $M\in[\N]^{\omega}$, and  $\gamma\le \beta$ be a limit ordinal number and let  $\eta(\gamma,n)$ be the sequence provided by Proposition  \ref{P:2.6}.
Assume that $A_0$  a singleton or the empty set and  that
$\Phi:\cS_{\beta\gamma}(A_0)\cap[M]^{<\omega}\rightarrow X$ is a semi-embedding of $\cS_{\beta\gamma}\cap[M]^{<\omega}$ into $X$ starting at $A_0$ that is $c$-refined.

If $A_0$ is a singleton, say $A_0=\{a_0\}$,
 let $\Psi:\cS_{\beta\eta(\gamma,a_0)}(A_0)\cap[M]^{<\omega}\rightarrow X$, with $\Psi(A) = \Phi(A)$, be the semi-embedding of $\cS_{\beta\eta(\gamma,a_0)}\cap[M]^{<\omega}$ into $X$ starting at $A_0$, that is $c$-refined,
given by Remark \ref{R:10.6}. Then, for every $E\in\MAX(\cS_{\beta\eta(\gamma,a_0)}(A_0)\cap[N]^{<\omega})$, we have
\begin{equation}
\label{E:12.1}
((x_{\Psi,A}^{(k)})_{k=1}^{s(\beta,\eta(\gamma,a_0),A)})_{A\in\cE_{\beta,\eta(\gamma,a_0)}(A_0\cup E)} = ((x_{\Phi,A}^{(k)})_{k=1}^{s(\beta,\gamma,A)})_{A\in\cE_{\beta,\gamma}(A_0\cup E)}.
\end{equation}

If $A_0=\emptyset$,
 let $a_0\in M$, set $A_0 = \{a_0\}$, $N = M\cap[a_0,\infty)$  and $\Psi=\Phi|_{\cS_{\beta\eta(\gamma,a_0)}(A_0)\cap[N]^{<\omega}}$,  which is, by Remark \ref{R:10.6},  a semi-embedding of $\cS_{\beta\eta(\gamma,a_0)}\cap[N]^{<\omega}$
into $X$ starting at $A_0$, that is $c$-refined. Then, again for every $E\in\MAX(\cS_{\beta\eta(\gamma,a_0)}(A_0)\cap[N]^{<\omega})$, 
\begin{equation}
\label{E:12.2}
((x_{\Psi,A}^{(k)})_{k=1}^{s(\beta,\eta(\gamma,a_0),A)})_{A\in\cE_{\beta,\eta(\gamma,a_0)}(A_0\cup E)} = ((x_{\Phi,A}^{(k)})_{k=1}^{s(\beta,\gamma,A)})_{A\in\cE_{\beta,\gamma}(A_0\cup E)}
\end{equation}
is the step block decomposition of $E$ with respect to $\Psi$.
\end{rem}

Before formulating  and proving the missing parts  our Main Theorems A and B (see upcoming Theorem \ref{T:12.6}) we present the argument which is the 
 main inductive step.
 
 Let $\gamma<\beta$. For $B\in\MAX( \cS_{\beta(\gamma+1)})=\cS_{\beta}[\cS_{\beta\gamma}]$ (Proposition \ref{P:2.6}),
 we let $B_1<B_2<\ldots <B_d$ the (unique) elements of $\MAX(\cS_{\beta\gamma})$ for which $B=\bigcup_{j=1}^d B_j$.
We  also define $\bar B=\{ \min(B_j): j\keq1,2,\ldots,d\}\kin\MAX(\cS_\beta)$ and for $i\keq1,\ldots, d$, $\bar B_j=\{ \min(B_j): j=1,2,\ldots, i\}$.

If $\emptyset \prec A\preceq B$ we can write $A$ as $A=\bigcup_{i=1}^{j-1} B_i \cup C$,  for some $j=1,2,\ldots, d$, and
some $\emptyset\prec C\preceq B_j$, and thus,  by Proposition \ref{P:4.4} $\zeta(\beta(\gamma+1),A)= \zeta(\beta, \bar B_j)\zeta(\beta\gamma, C)$.
We define
\begin{equation}
\label{E:12.4}
\tilde\cE_{\beta,\gamma+1}(B) = \left\{  (\cup_{i=1}^jB_j)\cup C \in\cE_{\beta,\gamma+1}(B):1\kleq j\kle d, l_1(\bar B_j)>0\text{ and }C\preceq B_{j+1}\right\}. 
\end{equation}

Let  $M\kin[N]^\omega$, $A_0$ be a subset of $\N$ that is either empty or a singleton,  and $\Phi:\cS_{\beta\gamma}(A_0)\cap[M]^{<\omega}\rightarrow X$ be a semi-embedding of
$\cS_{\beta\gamma}\cap[M]^{<\omega}$ into $X$, starting after $A_0$, that is $c$-refined, for some $0\kle c\kleq1$.
 Let also $E\kin\MAX(\cS_{\beta(\gamma+1)})(A_0) \cap[M]^{<\omega}$, and put   $B\keq A_0\cup E $. For $1\kleq j\kle d$, with $l_1(\bar B_j)\kgr0$, put $M^{(j)}\keq  M\cap(\max(B_j),\infty)$, and
$\Phi^{(j)}_E: \cS_{\beta\gamma}\cap[M^{(j)}]^{<\omega}\rightarrow X$~with
 \begin{equation} 
 \label{E:12.5}
 \Phi^{(j)}_E(C) = \frac{1}{\zeta(\beta,\bar B_{j+1})}\Phi\left(\left(\left(\cup_{1\le i\le j}B_i\right)\setminus A_0\right)\cup C\right).
 \end{equation}
 Recall that by Lemma \ref{L:10.7} $\Phi^{(j)}_{E}$ is a semi embedding  of $\cS_{\beta\gamma} \cap[M^{(j)}]^{<\omega}$ into $X$ starting at $\emptyset$, that is $c$ refined, 
 recall that for $1\le j< d$
 $$\cE_{\beta,\gamma+1}^{(j)}(A_0\cup E) =  \{A\in\cE_{\beta,\gamma+1}(A_0\cup E): \Cp_{\beta,\gamma}(A,1) = \cup_{i\le j}B_i\},$$
 and, moreover, if $l_1(\bar B_j)>0$ define
 \begin{equation}
 \label{the yjs}
y^{(j)}_{\Phi,E} = \sum_{A\in\cE_{\beta,\gamma+1}^{(j)}(A_0\cup E)}\frac{\zeta(\beta(\gamma+1),A)}{\zeta(\beta,\bar B_j)}\sum_{k=2}^{s(\beta,\gamma+1,A)}x_{\Phi,A}^{(k)}.  
 \end{equation}

\begin{rem}
\begin{enumerate}\item[1)$\qquad\,\,$]
\hskip-1cm By Lemma \ref{L:10.7}, each $\Phi^{(j)}_E$ is a semi-embedding of $\cS_{\beta\gamma}\cap[M^{(j)}]^{<\omega}$ into $X$, starting at $\emptyset$, that is $c$-refined.

\item[2)]  We note for later use that  $(y^{(j)}_{\Phi,E})_{j=1}^d$ is a  sequence in $X$ which satisfies the conditions of  the sequence $(x_j)_{j=1}^d$ in  Theorem \ref{T:7.1} with $\alpha=\beta$ and thus, assuming that $\Sz(X)\le \omega^\beta$, also its conclusion.

\item[3)] By the definition of the components of a set $A$, we conclude  (note that for $j=d$ we have $l_1(\bar B_d) = 0$)
\begin{equation}
\label{tilde script E description}
\tilde\cE_{\beta,\gamma+1}(A_0\cup E) = \bigcup_{\substack{1\le j< d\\ l_1(\bar B_j)>0}}\cE_{\beta,\gamma+1}^{(j)}(A_0\cup E),
\end{equation}
which yields
\begin{equation}
\label{aint it a beauty}
\begin{split}
\sum_{A\in\tilde\cE_{\beta,\gamma+1}(A_0\cup E)}&\zeta(\beta(\gamma+1),A)\sum_{k=1}^{s(\beta,\gamma+1,A)}x_{\Phi,A}^{(k)} =\\
&\sum_{\substack{1\le j < d\\ l_1(\bar B_j)>0}}\zeta(\beta,\bar B_j)y^{(j)}_{\Phi,E} + \sum_{A\in\tilde\cE_{\beta,\gamma+1}(A_0\cup E)}\zeta(\beta(\gamma+1),A)x_{\Phi,A}^{(1)}.
\end{split}
\end{equation}
\end{enumerate}
\end{rem}

\begin{lem}
\label{taking out some lines don't hurt nobody}
Let  $\gamma< \beta $, $M\in[N]^\omega$, $A_0$  be either empty or a singleton in $\N$,
and let  $\Phi:\cS_{\beta(\gamma+1)}(A_0)\cap[M]^{<\omega}\rightarrow X$ be a semi-embedding of $\cS_{\beta(\gamma+1)}\cap[M]^{<\omega}$ into $X$, starting after $A_0$, that is $c$-refined, for some $0<c\le 1$.
Then, for every $E\in\MAX(\cS_{\beta(\gamma+1)}(A_0)\cap[M]^{<\omega})$
\begin{equation}
\label{taking out some lines don't hurt nobody equation}
\Big\|\Phi(E) - \sum_{A\in\tilde\cE_{\beta,\gamma+1}(A_0\cup E)}\zeta(\beta(\gamma+1),A)\sum_{k=1}^{s(\beta,\gamma+1,A)}x_{\Phi,A}^{(k)}\Big\| < \frac{3}{\min(A_0\cup E)}.
\end{equation}
\end{lem}

\begin{proof}
Recall that for $A\in\cE_{\beta,\gamma+1}(A_0\cup E)$, we have $x_A = \sum_{k=1}^{s(\beta,\gamma+1,A)}x_{\Phi,A}^{(k)}$ and that $\Phi(E) = \sum_{A\preceq A_0\cup E}\zeta(\beta(\gamma+1),A)x_A$. Hence,
$$\Big\|\Phi(E) -\!\!\! \sum_{A\in\tilde\cE_{\beta,\gamma+1}(A_0\cup E)}\zeta(\beta(\gamma+1),A)\sum_{k=1}^{s(\beta,\gamma,A)}x_{\Phi,A}^{(k)}\Big\| =
\Big\|\sum_{\substack{A\preceq A_0\cup E:\\A\not\in\tilde\cE_{\beta,\gamma+1}(A_0\cup E)}}\zeta(\beta(\gamma+1),A)x_A\Big\|.$$
We calculate:
\begin{align}
 &\sum_{\substack{A\preceq A_0\cup E\\A\not\in\tilde\cE_{\beta,\gamma+1}(A_0\cup E)}}\zeta(\beta(\gamma+1),A)   \\
 &=\!\!\sum_{\substack{A\preceq A_0\cup E\\A\not\in\cE_{\beta,\gamma+1}(A_0\cup E)}}\zeta(\beta(\gamma+1),A)
  + \sum_{\substack{1\le j< d\\ l_1(\bar B_j) = 0}}\sum_{A\in\cE_{\beta,\gamma+1}^{(j)}(A_0\cup E)}\zeta(\beta(\gamma+1),A)\nonumber\\
 &=\!\! \!\!\sum_{\substack{A\preceq A_0\cup E\\A\not\in\cE_{\beta,\gamma+1}(A_0\cup E)}}\!\!\!\zeta(\beta(\gamma+1),A)
\kplus \sum_{\substack{1\le j< d\\ l_1(\bar B_j) = 0}}\zeta(\beta,\bar B_{j+1})\sum_{A\in\cE_{\beta,\gamma+1}^{(j)}(A_0\cup E)}\!\!\!\zeta(\beta\gamma,A\setminus(\cup_{i\le j}B_j))\nonumber\\
 &<
 \frac{2}{\min(A_0\cup E)} + \frac{1}{\min(\bar B)} = \frac{3}{\min(A_0\cup E)},\nonumber
\end{align}
where the last inequality follows from  Lemmas \ref{L:4.6} (v) and \ref{L:11.5}.
\end{proof}

\begin{lem}\label{what y^(j)'s look like}
Let $\gamma$,  $M\in[N]^\omega$, $A_0$, $\Phi$, and $c$ be as in the statement of Lemma \ref{taking out some lines don't hurt nobody}.
 If $E\in\MAX(\cS_{\beta(\gamma+1)}(A_0)\cap[M]^{<\omega})$ and  $B_1\kle\cdots\kle B_d$ are in $\MAX(\cS_{\beta\gamma})$, with  $A_0\cup E = \cup_{j=1}^dB_j$, and  $\bar B = \{\min B_j: 1\kleq j\kleq d\}\in\MAX(\cS_{\beta})$,
and if  $(y^{(j)}_{E,\Phi})_{j=1}^d$ is defined as in \eqref{the yjs},
then for $j\keq1,\ldots,d-1$ with $l_1(\bar B_j)\kgr0$, where $\bar B_j \keq \{\min B_i: 1\le i\le j\}$, we have $\|y^{(j)}_{\Phi,E}\|\kleq 1$ and $\ran(y^{(j)}_{\Phi,E})\subset(\max(\bar B_j),\max(\bar B_{j+2}))$. Put $\max(\bar B_{d+1}) = \infty$. Then
\begin{align}
\label{the stairs break up like this}
\sum_{A\in\tilde\cE_{\beta(\gamma+1)}(A_0\cup E)}&\zeta(\beta(\gamma+1),A)\sum_{k=1}^{s(\beta,\gamma+1,A)}x_{\Phi,A}^{(k)} \\
&=\sum_{A\in\tilde\cE_{\beta,\gamma+1}(A_0\cup E)}\zeta(\beta(\gamma+1),A)x_{\Phi,A}^{(1)} +  \sum_{\substack{1\le j < d\\ l_1(\bar B_j)>0}}\zeta(\beta,\bar B_{j+1})y_{\Phi,E}^{(j)}.
\notag
\end{align}

\end{lem}

\begin{proof}
Observe that \eqref{the stairs break up like this} immediately follows from \eqref{aint it a beauty} and the fact that for $1\le j<d$ with $l_1(\bar B_j) >0$, we have $\zeta(\beta,\bar B_j) = \zeta(\beta,\bar B_{j+1})$.
For $1\le j< d$ with $l_1(\bar B_j) > 0$ and $A\in\cE_{\beta,\gamma+1}^{(j)}(A_0\cup E)$ note that $\cup_{i=1}^jB_i\prec A\preceq\cup_{i=1}^{j+1}B_i$ and
\begin{align}
\label{just an interval projection}
u_A &=  \sum_{k=2}^{s(\beta,\gamma+1,A)}x_{\Phi,A}^{(k)} \\
&= P_{(\max\supp(x_{\Cp_{\beta,\gamma}(A,1)}), \max\supp(x_A)]}(x_A)= P_{(\max\supp (x_{\cup_{i=1}^jB_j}),\max\supp (x_A)]}x_A\notag
\end{align}
i.e. $\|u_A\|\le 1$ and by Definition \ref{D:10.5} (ii) we obtain
\begin{equation*}
\begin{split}
&\max(\bar B_j) \le \max(\cup_{i=1}^j B_i)\le \max\supp(x_{\cup_{i=1}^jB_i})< \min\supp(u_A)\text{ and}\\
&\max\supp(u_A)\le \max\supp (x_A)< \min\{ m\in M: m\kgr \max(A)\} \le \max (\bar B_{j+2})
\end{split}
\end{equation*}
which yields
\begin{equation}
\label{about those ranges}
\ran(u_A)\subset(\max(\bar B_j),\max(\bar B_{j+2})).
\end{equation}
Furthermore, we have $\zeta(\beta(\gamma + 1),A) = \zeta(\beta,\bar B_{j+1})\zeta(\beta\gamma,A\setminus(\cup_{i=1}^jB_i))$,
and since $\zeta(\beta, \bar B_{j+1}) = \zeta(\beta,\bar B_j)$ (by $l_1(\bar B_j) > 0$), we obtain
\begin{equation}
\label{about those coefficients}
\frac{\zeta(\beta(\gamma+1),A)}{\zeta(\beta,\bar B_j)} = \zeta(\beta\gamma,A\setminus(\cup_{i=1}^jB_i)).
\end{equation}
We deduce
$$y^{(j)}_{\Phi,E} = \sum_{C\in\cE_{\beta,\gamma}(B_{j+1})}\zeta(\beta\gamma,C)u_{(\cup_{i=1}^jB_i)\cup C}.$$
The above in combination with \eqref{about those ranges}, \eqref{about those coefficients}, and the fact that $\|u_A\|\le 1$ yield that $\|y^{(j)}_{\Phi,E}\|\le 1$ and $\ran(y^{(j)}_{\Phi,E})\subset(\max(\bar B_j),\max(\bar B_{j+2}))$.
\end{proof}

\begin{lem}\label{those cut offs are step decompositions of smaller embeddings}
Let $\gamma$, $M\in[N]^\omega$, $A_0$, $\Phi$, and $c$ be as in the statement of Lemma \ref{taking out some lines don't hurt nobody}.
Let $E\in\MAX(\cS_{\beta(\gamma+1)}(A_0)\cap[M]^{<\omega})$, $(B_j)_{j=1}^d$, $\bar B$, $(\bar B_j)_{j=1}^d$ be as in the statement of Lemma \ref{what y^(j)'s look like} and let
 $\Phi^{(j)}_E$, $j=1,\ldots, d-1$ and $M^{(j)}\in [M]^{\omega}$ be defined as in \eqref{E:12.5}.
For $j=1,\ldots,d-1$ with $l_1(\bar B_j)>\emptyset$ denote by
$$\left(\left(z_{\Phi^{(j)}_E,C}^{(k)}\right)_{k=1}^{s(\beta,\gamma',C)}\right)_{C\in\cE_{\beta,\gamma'}(B_{j+1})}$$
the block step decomposition of $B_{j+1}$ with respect to $\Phi^{(j)}_E$. Then for $C\in\cE_{\beta,\gamma}(B_{j+1})$  we have $s(\beta,\gamma+1,(\cup_{i=1}^jB_i)\cup C) = s(\beta,\gamma,C) + 1$ and
 $$z_{\Phi^{(j)}_E,C}^{(k)} = x_{\Phi,(\cup_{i=1}^jB_i)\cup C}^{(k+1)},$$
for $k=1,\ldots, s(\beta,\gamma,C)$.
\end{lem}

\begin{proof}
Fix $C\in\cE_{\beta,\gamma}(B_{j+1})$. By \eqref{E:11.2b}, if we set $A = (\cup_{i=1}^jB_i)\cup C$, then $s(\beta,\gamma+1,A) = s(\beta,\gamma,C) + 1$ and 
$\Cp_{\beta,\gamma+1}(A,i+1) = \Cp_{\beta,\gamma}(A,i)$ for $i=1,\ldots,s(\beta,\gamma,C)$,
and  $\Cp_{\beta,\gamma+1}(A,1) = \cup_{i=1}^jB_i$.
Fix $1\le k\le s(\beta,\gamma,C)$. 
 Let $\{z_{C} :C\kin \cS_{\beta\gamma}\cap[M^{(j)}]^{<\omega}\}$ be the family generating $\Phi_E^{(j)}$ and $n_0 = \max\supp(x_{\cup_{i=1}^jB_i})$, then 
by Definition \ref{D:12.1} and Lemma  \ref{L:10.7} we have
$$x_{\Phi,A}^{(k+1)} = P_{I(A,k+1)}(x_A)\text{ and } z_{\Phi^{(j)}_E,C}^{(k)} = P_{I(C,k)}(z_C)=  P_{I(C,k)}\big(P_{(n_0,\infty)}(x_{\bigcup_{i\le j} B_i\cup C}) \big),$$
with (if $k=1$ replace $\max\supp(z_{\cup_{i=1}^{k-1}\Cp_{\beta,\gamma}(C,i)}$ by $n_0$)
\begin{align*}
I&(C,k) = (\max\supp(z_{\cup_{i=1}^{k-1}\Cp_{\beta,\gamma}(C,i)}),\max\supp(z_{\cup_{i=1}^{k}\Cp_{\beta,\gamma}(C,i)})]\\
&=\! (\max\supp (P_{(n_0,\infty)}x_{(\cup_{i=1}^jB_i)\cup(\cup_{i=1}^{k-1}\Cp_{\beta,\gamma}(C,i))}),\!\max\supp(P_{(n_0,\infty)}x_{(\cup_{i=1}^jB_i)\cup(\cup_{i=1}^{k}\Cp_{\beta,\gamma}(C,i))})]\\
&= (\max\supp (x_{(\cup_{i=1}^jB_i)\cup(\cup_{i=1}^{k-1}\Cp_{\beta,\gamma}(C,i))}),\max\supp(x_{(\cup_{i=1}^jB_i)\cup(\cup_{i=1}^{k}\Cp_{\beta,\gamma}(C,i))})]
 = I(A,k\kplus1),
\end{align*}
where we used $n_0\le \max\supp (x_{\cup_{i=1}^jB_i})$, which  follows from Definition \ref{D:10.5} (ii). Hence,
$$z_{\Phi^{(j)}_E,C}^{(k)} = P_{I(A,k+1)}P_{(n_0,\infty)}(x_A) = P_{I(A,k+1)}(x_A) = x_{\Phi,A}^{(k+1)}.$$
\end{proof}

\begin{prop}
\label{if szlenk is small enough and large Schreier semi-embeds then the steps have large norm}
Assume that  $\Sz(X) \le \omega^\beta$. Then, for every $1\ge c>0$ and $M\in[\N]^\omega$,
there exists $N\in[M]^\omega$ with the following property:
for every $\gamma\le\beta$, $L\in[N]^\omega$, and $A_0\subset N$ that is either empty or a singleton, every   semi-embedding $\Phi:\cS_{\beta\gamma}(A_0)\cap[L]^{<\omega}\rightarrow X$ from $\cS_{\beta\gamma}\cap[L]^{<\omega}$ into $X$
starting at $A_0$, that is $c$-refined, every  $E\in\MAX(\cS_{\beta\gamma}(A_0)\cap[L]^{<\omega})$ and every
$(\beta,\gamma)$-special family of convex combinations $\{r(A,k): A\in\cE_{\beta,\gamma}(A_0\cup E), 1\le k\le s(\beta,\gamma,A)\}$
we have
\begin{equation}
\label{this is too big}
\sum_{A\in\cE_{\beta,\gamma}(A_0\cup E)}\zeta(\beta\gamma,A) \sum_{k=1}^{s(\beta,\gamma,A)} r(A,k)\|x_{\Phi,A}^{(k)}\| \ge \frac{c}{3},
\end{equation}
 where $(x_{\Phi,A}^{(k)})_{A\in\cE_{\beta,\gamma(A_0,E)}}$\!\! is the block step decomposition of $E$ for $\Phi$  (Definition~\ref{D:12.1}).
\end{prop}

\begin{proof}
Fix $M\in[\N]^\omega$ and $1\ge c>0$.
Choose $\vp > 0$ so that $c/2 - \vp > c/3$ and then apply Theorem \ref{T:7.1} to find $N\in [M]^\infty$ so that \eqref{E:7.1.1}
is satisfied for that $\vp$ and $\beta$, and, moreover,
\begin{equation}
\label{growth of subset has to satisfy this}
\Big(\frac{c}{2} - \vp - \frac{4}{\min(N)}\Big)\prod_{m\in N}\Big(1 - \frac{1}{m}\Big) > \frac{c}{3}.
\end{equation}
We claim that this is the desired set. We shall prove by transfinite induction on $\gamma$ the following statement:
if $\gamma\le\beta$, $L\in[N]^\omega$, $A_0\subset L$ that is either empty or a singleton, and $\Phi:\cS_{\beta\gamma}(A_0)\cap[L]^{<\omega}\rightarrow X$ is a semi-embedding of $\cS_{\beta\gamma}\cap[L]^{<\omega}$ into $X$
starting at $A_0$, that is $c$-refined, then for any $E\in\MAX(\cS_{\beta\gamma}(A_0)\cap[L]^{<\omega})$ and any $(\beta,\gamma)$-special family of convex combinations
$\{r(A,k): A\in\cE_{\beta,\gamma}(A_0\cup E), 1\le k\le s(\beta,\gamma,A)\}$
we have
\begin{equation}
\label{this is the inductive assumption of proposition }
\sum_{A\in\cE_{\beta,\gamma}(A_0\cup E)}\zeta(\beta\gamma,A)\sum_{k=1}^{s(\beta,\gamma,A)}r(A,k)\|x_{\Phi,A}^{(k)}\| \ge \left(\frac{c}{2} - \vp - \frac{4}{\min(L)}\right)\prod_{m\in L}^\infty\left(1 - \frac{1}{m}\right).
\end{equation} 
In conjunction with \eqref{growth of subset has to satisfy this}, this will yield the desired result.
Let $\gamma = 1$, $L\subset N$, $A_0$ be a subset of $L$ that is either empty or a singleton, and $\Phi:\cS_\beta(A_0)\cap[L]^{<\omega}\rightarrow X$ be a semi-embedding of $\cS_{\beta}\cap[L]^{<\omega}$ into $X$ starting at $A_0$
that is $c$-refined. Let $E\in\MAX(\cS_\beta(A_0)\cap[N]^{<\omega})$.

By \eqref{E:11.2a}, for each $A\in\cE_{\beta,1}(A_0\cup E)$ we obtain that the block step decomposition of $x_A$ is just $(x_{\Phi,A}^{(1)}) = (x_A)$, and hence if $\{r(A,k):A\in\cE_{\beta,1}(A_0\cup E),  1\le k\le s(\beta,1,A)\}$
is a $(\beta,1)$-special family of convex combinations then $r(A,1) = 1$.
Hence,
\begin{align}
\sum_{A\in\cE_{\beta,1}(A_0\cup E)}\zeta(\beta,A)&\sum_{k=1}^{s(\beta,1,A)}r(A,k)\|x_{\Phi,A}^{(k)}\| = \sum_{A\in\cE_{\beta,1}(A_0\cup E)}\zeta(\beta,A)\|x_{A}\|\nonumber\\
&\ge \Big\|\sum_{A\preceq A_0\cup E}\zeta(\beta,A)x_A\Big\| - \sum_{\substack{A\preceq A_0\cup E\\ A\notin\cE_{\beta,1}(A_0\cup E)}}\zeta(\beta,A)\nonumber\\
&\ge \|\Phi(E)\| - \frac{2}{\min(A_0\cup E)} \text{ (by \eqref{E:11.5.1})}\notag\\
&\ge \frac{c}{2} - \frac{c}{2}\zeta(\beta,A_0) - \frac{2}{\min(A_0\cup B)} \text{ (by  \ref{D:10.5} (iv))}\notag\\
&> \frac{c}{2} - \frac{3}{\min(L)} \text{ (by  \eqref{E:11.5})}.\notag
\end{align}

To verify the induction step, let first $\gamma<\beta$ be an ordinal number for which the conclusion holds. Let $L\subset N$, $A_0$ be a subset of $L$ that is either empty or a singleton,
and $\Phi:\cS_{\beta(\gamma+1)}(A_0)\cap[L]^{<\omega}\rightarrow X$ be a semi-embedding of $\cS_{\beta(\gamma+1)}\cap[L]^{<\omega}$ into $X$ starting at $A_0$
that is $c$-refined. Let $E\in\MAX(\cS_{\beta(\gamma+1)}(A_0)\cap[L]^{<\omega})$ and $\{r(A,k): A\in\cE_{\beta,\gamma+1}(A_0\cup E), 1\le k\le s(\beta,\gamma+1,A)\}$ be a $(\beta,\gamma+1)$-special family of convex combinations.
Let $A_0\cup E = \cup_{j=1}^dB_j$, where $B_1<\cdots<B_d$ are in $\MAX(\cS_{\beta\gamma}\cap[L]^{<\omega})$ and $\bar B = \min\{B_j: 1\le j\le d\}$ is in $\MAX(\cS_{\beta})$.
By   Lemma \ref{what y^(j)'s look like}  and the choice of the set $N$, we obtain
\begin{equation}
\label{yjs are small}
\Big\|\!\!\sum_{\substack{1\le j <d\\ l_1(\bar B_j)>0}}\zeta(\beta,\bar B_{j+1})y_{\Phi,E}^{(j)}\Big\| < \vp.
\end{equation}
Combining first Lemmas \ref{taking out some lines don't hurt nobody} and \ref{what y^(j)'s look like},  then applying  Definition \ref{D:10.5} (iv), \eqref{E:11.5},
and finally using \eqref{yjs are small} we deduce
\begin{align}\label{ARGH}
\sum_{A\in\tilde\cE_{\beta,\gamma\kplus1}(A_0\cup E)}&\!\!\zeta(\beta(\gamma\kplus1),A)\|x_{\Phi,A}^{(1)}\|\\
&\ge
\Big\|\sum_{A\in\tilde\cE_{\beta,\gamma\kplus1}(A_0\cup E)}\!\!\zeta(\beta(\gamma+1),A)x_{\Phi,A}^{(1)}\Big\|\notag\\
&\ge\|\Phi(E)\| - \Big\|\!\!\sum_{\substack{1\le j <d\\ l_1(\bar B_j)>0}}\zeta(\beta,\bar B_{j+1})y_{\Phi,E}^{(j)}\Big\| - \frac{3}{\min(A_0\cup E)}\notag\\
&\ge \frac{c}{2} - 
\frac{c}{2}\zeta(\beta(\gamma\kplus1),A_0) - \vp - \frac{3}{\min(A_0\cup E)}\ge \frac{c}{2} - \vp - \frac{4}{\min(L)}.\notag
\end{align}
We distinguish two cases. If $r(A,1) \keq 1$, for all $A\kin\cE_{\beta,\gamma+1}(A_0\cup E)$, then
\begin{equation*}
\begin{split}
\sum_{A\in\cE_{\beta,\gamma+1}(A_0\cup E)}&\zeta(\beta(\gamma+1),A)\sum_{k=1}^{s(\beta,\gamma+1,A)}r(A,k) \|x_{\Phi,A}^{(k)}\|\\
&= \sum_{A\in\cE_{\beta,\gamma+1}(A_0\cup E)}\zeta(\beta(\gamma+1),A) \|x_{\Phi,A}^{(1)}\|
\ge \sum_{A\in\tilde\cE_{\beta,\gamma+1}(A_0\cup E)}\zeta(\beta(\gamma+1),A)\|x_{\Phi,A}^{(1)}\|.
\end{split}
\end{equation*}
By \eqref{ARGH}, the result follows in that case.

Otherwise we have $r_1=r(A,1)<1$, for  all $A\in\cE_{\beta,\gamma+1}(A_0\cup E)$.
For $1\le j< d$, with $l_1(\bar B_j)>0$, define $L^{(j)} =  L\cap(\max(A_0\cup(\cup_{1\le i \le j}B_i)),\infty)$, and
$\Phi^{(j)}_E: \cS_{\beta\gamma}\cap[L^{(j)}]^{<\omega}\rightarrow X$ as in \eqref{E:12.5}.
By Lemma \ref{L:10.7}, each $\Phi^{(j)}_E$ is a semi-embedding of $\cS_{\beta\gamma}\cap[L^{(j)}]^{<\omega}$ into $X$, starting at $\emptyset$,
that is $c$-refined.

By Lemma \ref{special families demoted}, the family $\{r^{(j)}(C,k): C\in\cE_{\beta,\gamma}(B_{j+1}) \}$, with $$r^{(j)}(C,k) = \frac{1}{1-r(A^{(1)},1)}r\left(\left(\cup_{i=1}^{j}B_i\right)\cup C,k+1\right)$$
for $k=1,\ldots,s(\beta,\gamma,C) = s(\beta,\gamma,(\cup_{i=1}^{j}B_i)\cup C) - 1$ and some $A^{(1)}\in\cE_{\beta,\gamma+1}(A_0\cup E)$, is a $(\beta,\gamma)$-special family of convex combinations. Hence, by the inductive assumption applied to the map $\Phi^{(j)}_E$, and Lemma
\ref{those cut offs are step decompositions of smaller embeddings} we deduce
\begin{align}\label{lotsa work to get here}
 \sum_{A\in\cE^{(j)}_{\beta,\gamma+1}(A_0\cup E)}&\!\!\!\zeta(\beta(\gamma + 1), A)\sum_{k=2}^{s(\beta,\gamma+1,A)}r(A,k)\|x_{\Phi,A}^{(k)}\| \\
=&\zeta(\beta,\bar B_j)\big(1 \kminus r(A^{(1)},1)\big)\sum_{C\in\cE_{\beta,\gamma}(B_{j+1})}\!\!\!\zeta(\beta\gamma,C)\sum_{k=1}^{s(\beta,\gamma,C)}r^{(j)}(C,k)\|z_{\Phi^{(j)}_E,C}^{(k)}\|\nonumber\\
\ge& \zeta(\beta,\bar B_j)\big(1 \kminus r(A^{(1)},1)\big)\Big(\frac{c}{2} - \vp - \frac{4}{\min( L^{(j)})}\Big)\prod_{m\in L^{(j)}}\Big(1 - \frac{3}{m}\Big).\nonumber
\end{align}
We combine \eqref{ARGH} with \eqref{lotsa work to get here}:
\begin{align}
 \sum_{A\in\cE_{\beta,\gamma+1}(A_0\cup E)}&\zeta(\beta(\gamma+1),A)\sum_{k=1}^{s(\beta,\gamma+1,A)}r(A,k)\|x_{\Phi,A}^{(k)}\|\nonumber\\
 &\ge \sum_{A\in\tilde\cE_{\beta,\gamma+1}(A_0\cup E)}\zeta(\beta(\gamma+1),A)\sum_{k=1}^{s(\beta,\gamma+1,A)}r(A,k)\|x_{\Phi,A}^{(k)}\|\nonumber\\
 &=\sum_{A\in\tilde\cE_{\beta,\gamma+1}(A_0\cup E)}\zeta(\beta(\gamma+1),A)r(A,1)\|x_{\Phi,A}^{(1)}\|\nonumber \\
 &\quad+\sum_{\substack{1\le j<d\\ l_1(\bar B_j)>0}}\sum_{A\in\cE^{(j)}_{\beta,\gamma+1}(A_0\cup E)}\zeta(\beta(\gamma + 1), A)\sum_{k=2}^{s(\beta,\gamma+1,A)}r(A,k)
 \|x_{\Phi,A}^{(k)}\|\nonumber\\
 & \ge r(A^{(1)},1)\Big(\frac{c}{2} - \vp - \frac{4}{\min(L)}\Big)\nonumber\\ 
 &\quad+ \sum_{\substack{1\le j<d\\ l_1(\bar B_j)>0}}\zeta(\beta,\bar B_j)(1 - r(A^{(1)},1))\Big(\frac{c}{2} - \vp - \frac{4}{\min(L)}\Big)\prod_{m\in L^{(j)}}\Big(1 - \frac{1}{m}\Big)\nonumber\\
 &\ge \Big(\sum_{\substack{1\le j<d\\ l_1(\bar B_j)>0}}\zeta(\beta,\bar B_j)\Big)\Big(\frac{c}{2} -
  \vp - \frac{4}{\min(L)}\Big)\prod_{m\in L^{(1)}}\Big(1 - \frac{1}{m}\Big)\nonumber\\
 &\ge \Big(1 - \frac{1}{\min(\bar B)}\Big)\Big(\frac{c}{2} - \vp - \frac{4}{\min(L)}\Big)\prod_{m\in L^{(1)}}\Big(1 - \frac{1}{m}\Big)\text{ (by Lemma \ref{L:4.6} (v))} \nonumber\\
 &\ge \Big(\frac{c}{2} - \vp - \frac{4}{\min(L)}\Big)\prod_{m\in L}\Big(1 - \frac{1}{m}\Big).\nonumber
\end{align}
Assume now  that $\gamma\le\beta$ is a limit ordinal number and that the claim holds for all $\gamma'<\gamma$. Let $L\in[ N]^{\omega}$, $A_0$ be a subset of $L$ that is either empty or a singleton,
and $\Phi:\cS_{\beta\gamma}(A_0)\cap[L]^{<\omega}\rightarrow X$ be a semi-embedding of $\cS_{\beta\gamma}\cap[L]^{<\omega}$ into $X$ starting at $A_0$
that is $c$-refined. We distinguish between two cases, namely whether $A_0$ is a singleton or whether it is empty. In the first case, $A_0 = \{a_0\}$, for some $a_0\in L$.
By Remark \ref{R:10.6}, the map $\Psi$ with $\Psi = \Phi$ can be seen as a semi-embedding of $\cS_{\beta\eta(\gamma,a_0)}\cap[L]^{<\infty}$ into $X$ starting at $A_0$,
that is $c$-refined. If $E\in\MAX(\cS_{\beta\gamma}(A_0)\cap[L]^{<\omega})$ then $E\in\MAX(\cS_{\beta\eta(\gamma,a_0)}(A_0)\cap[L]^{<\omega})$ and by \eqref{E:12.1} we have
\begin{equation*}
((x_{\Psi,A}^{(k)})_{k=1}^{s(\beta,\eta(\gamma,a_0),A)})_{A\in\cE_{\beta,\eta(\gamma,a_0)}(A_0\cup E)} = ((x_{\Phi,A}^{(k)})_{k=1}^{s(\beta,\gamma,A)})_{A\in\cE_{\beta,\gamma}(A_0\cup E)},
\end{equation*}
whereas if $\{r(A,k): A\in\cE_{\beta,\gamma}(A_0\cup E), 1\kleq k\kleq s(\beta,\gamma,A)\}$ is a $(\beta,\gamma)$-special family of convex combinations,
by the remark following Definition \ref{being in the same component is the same as having common minimum}, it is a $(\beta,\eta(\gamma,a_0))$-special family of convex combinations as well. Applying the inductive assumption for $\eta(\gamma,a_0)\kle\gamma$
yields
\begin{align}\label{stars and flowers}
\sum_{A\in\cE_{\beta,\gamma}(A_0\cup E)}&\zeta(\beta\gamma,A)\sum_{k=1}^{s(\beta,\gamma,A)}r(A,k)\left\|x_{\Phi,A}^{(k)}\right\|\\
&=\sum_{A\in\cE_{\beta,\eta(\gamma,a_0)}(A_0\cup E)}\zeta(\beta\eta(\gamma,a_0),A)\sum_{k=1}^{s(\beta,\eta(\gamma,a_0),A)}r(A,k)\left\|x_{\Psi,A}^{(k)}\right\|\nonumber\\
&\ge \left(\frac{c}{2} - \vp - \frac{4}{\min(L)}\right)\prod_{m\in L}^\infty\left(1 - \frac{1}{m}\right). \nonumber
\end{align}

In the second case, $A_0$ is empty. Let $B\in\MAX(\cS_{\beta\gamma}\cap[L]^{<\omega})$ and set $a_0 = \min(B)$. By Remark \ref{R:10.6}, if $L' = L\cap[a_0,\infty)$, then
the map $\Psi:\cS_{\beta\eta(\gamma,a_0)}(A_0)\cap[L']\rightarrow X$ with $\Psi(A) = \Phi(A_0\cup A)$, is  a semi-embedding of $\cS_{\beta\eta(\gamma,a_0)}\cap[L']^{<\omega}$ into $X$ starting at $A_0$ that is $c$-refined.
By \eqref{E:12.2} we obtain
$$((x_{\Psi,A}^{(k)})_{k=1}^{s(\beta,\eta(\gamma,a_0),A)})_{A\in\cE_{\beta,\eta(\gamma,a_0)}(B)} = ((x_{\Phi,A}^{(k)})_{k=1}^{s(\beta,\gamma,A)})_{A\in\cE_{\beta,\gamma}(A_0\cup E)}$$
whereas if $\{r(A,k): A\in\cE_{\beta,\gamma}(B), 1\le k\le s(\beta,\gamma,A)\}$ is $(\beta,\gamma)$-special family of convex combinations,
by the Remark following Definition \ref{special families of convex combinations}, it is a $(\beta,\eta(\gamma,a_0))$-special family of convex combinations as well.
The result follows in the same manner as in \eqref{stars and flowers}.
\end{proof}

\begin{thm}\label{T:12.6}
Assume that  $X$ is a   reflexive and separable Banach space,  with the property that $\Sz(X) \kleq \omega^\beta$ and $\Sz(X^*)\kle \beta$.
Then for no $L\in[\N]^\omega$, does there exists a semi-embedding of $\cS_{\beta^2}\cap[L]^{<\omega}$ into $X$, starting at $\emptyset$.
\end{thm}

\begin{proof}
By the main Theorem of \cite{Sch} we can embed $X$ into a  reflexive space $Z$ with basis,  so that $\Sz(Z)=\Sz(X)$,
and $\Sz(Z^*)=\Sz(X^*)$. Thus we
 may assume that $X$ has a basis, which must be shrinking and boundedly complete, since $X$ is reflexive.
By re-norming $X$, we may assume that the bases of $X$ and $X^*$ are bimonotone.
Choose $\alpha_0$ with $\Sz(X^*) \le \omega^{\alpha_0} < \beta$ (this is possible due to the form of $\beta$). Note that
\begin{equation}
\label{the cb indexes equation}
\CB(\cS_{\alpha_0}) = \omega^{\alpha_0} + 1 < \beta + 1 = \CB(\cF_{\beta,\beta}).
\end{equation}
Towards a contradiction, assume that there exists $L\in[\N]^\omega$ and a semi-embedding $\Psi$ of $\cS_{\beta^2}\cap[L]^{<\omega}$ into $X$, starting at $\emptyset$.
By Lemma \ref{L:10.8} there exist $1\kgeq c\kgr0$, $M\in[L]^\omega$, and a semi-embedding $\Phi$ of $\cS_{\beta^2}\cap[M]^{<\omega}$ into $X$,
starting at $\emptyset$, that is $c$-refined. Applying  Proposition \ref{if szlenk is small enough and large Schreier semi-embeds then the steps have large norm}, we may pass to a further
subset of $M$, again denoted by $M$, so that \eqref{this is too big} holds.

Fix $0<\vp < c/3$ and apply Theorem \ref{T:7.1} to the space $X^*$ and the ordinal number $a_0$ to find a further subset of $M$, which we again denote by $M$,
so that \eqref{E:7.1.1} is satisfied.

By  Propositions \ref{P:2.2} and \ref{P:2.14}  we may pass to a subset of $M$, again denoted by $M$, so that $\cS_{\alpha_0}\cap[M]^{<\omega}\subset\cF_{\beta,\beta}$. By Lemma \ref{L:11.2} we obtain that for any $B\in\MAX(\cS_{\beta^2}\cap[M]^{<\omega})$ and $A\in\cE_{\beta,\beta}(B)$, there exists $\tilde A\in\MAX(\cS_{\alpha_0})$ with
\begin{equation}
 \label{this should be useful somewhere}
\tilde A \preceq \left\{\min(\Cp_{\beta,\beta}(A,k)): 1\le k\le s(\beta,\beta,A)\right\}.
\end{equation}

Choose $B\in\MAX(\cS_{\beta^2}\cap[M]^{<\omega})$. We will define a $(\beta,\beta)$-special family of convex combinations
$\{r(A,k): A\in\cE_{\beta,\beta}(B), 1\le k\le s(\beta,\beta,A)\}$.
For $A\in\cE_{\beta,\beta}(B)$ let $\tilde A = \{a_1^A,\ldots,a_{d_A}^A\}\in\MAX(\cS_{\alpha_0})$ be as in \eqref{this should be useful somewhere}. For $1\le k\le s(\beta,\beta,A)$ set 
\begin{equation}\label{this is indeed a special family}
r(A,k) =
\left\{
	\begin{array}{ll}
		\zeta(\alpha_0,\tilde A_k)  & \mbox{if } k\le \#\tilde A \\
		0 & \mbox{otherwise,}
	\end{array}
\right.
\end{equation}
where $\tilde A_k = \{a_1^A,\ldots,a_{k}^A\}$ for $1\le k\le \#\tilde A$.
We will show that this family satisfies Definition \ref{special families of convex combinations} (i) and (ii). The first assertion is straightforward, to
see the second one let $A^{(1)}$, $A^{(2)}\in\cE_{\beta,\beta}(B)$, so that if $(D_{k}(A^{(1)}))_{k=1}^{s(\beta,\beta,A)}$ and
 $(D_{k}(A^{(2)}))_{i=1}^{s(\beta,\beta,A')}$ are the maximal chains of $\cA_{\beta,\beta}(B)$
provided by Lemma \ref{L:11.2}, then for some $k$ we have $D_k(A^{(1)}) = D_k(A^{(2)})$.
By Lemmas \ref{being in the same component is backwards closed} and \ref{being in the same component is the same as having common minimum} we obtain 
$\min(\Cp_{\beta,\beta}(A^{(1)},m)) = \min(\Cp_{\beta,\beta}(A^{(2)},m))$ for $m=1,\ldots,k$,
which implies $\tilde A_m^{(1)} = \tilde A_m^2$ for $m=1,\ldots,\min\{k, \#\tilde A^{(1)}\}$.
By \eqref{this is indeed a special family} it
easily follows that $r(A^{(1)},k) = r(A^{(2)},k)$. Since \eqref{this is too big} is satisfied, we obtain
\begin{equation}
\label{this will ruin it all}
\sum_{A\in\cE_{\beta,\beta}(B)}\zeta(\beta^2,A)\sum_{k=1}^{d_A}\zeta(\alpha_0,\tilde A_k)\left\|x_{\Phi,A}^{(k)}\right\|\ge \frac{c}{3}.
\end{equation}

For each $A\in\cE_{\beta,\beta}(B)$ and $k=1,\ldots,d_A$ choose  $f_A^{(k)}$ in $S_{X^*}$, with $f_A^{(k)}(x^{(k)}_A)=\|x^{(k)}_A\|$ and
\begin{align*}
\ran(f_A^{(k)})&\subset\ran(x_{\Phi,A}^{(k)})\\
&\subset\big(\max\supp(x_{\cup_{i=1}^{k-1}\Cp_{\beta,\beta}(A,i)}),\max\supp(x_{\cup_{i=1}^{k}\Cp_{\beta,\beta}(A,i)})\big]\\
& \subset \big(\min (\Cp_{\beta,\beta}(A,k-1)),\min(\Cp_{\beta,\beta}(A,k+1))\big) = (\max(\tilde A_{k-1}), \max(\tilde A_{k+1})),\nonumber
\end{align*}
where the third  inclusion  follows from Definition \ref{D:10.5} (ii).
As \eqref{E:7.1.1} is satisfied, we obtain that for all $A\in\cE_{\beta,\beta}(B)$
\begin{equation}
\label{these don't go together}
\left\|\sum_{k=1}^{d_A}\zeta(\alpha_0,\tilde A_k)f_A^{(k)}\right\| < \vp.
\end{equation}
We finally calculate 
\begin{align*}
\frac{c}{3} &\le \sum_{A\in\cE_{\beta,\beta}(B)}\zeta(\beta^2,A)\sum_{k=1}^{d_A}\zeta(\alpha_0,\tilde A_k)\big\|x_{\Phi,A}^{(k)}\big\|
\text{ (by \eqref{this will ruin it all})} \\
& = \sum_{A\in\cE_{\beta,\beta}(B)}\zeta(\beta^2,A)\sum_{k=1}^{d_A}\zeta(\alpha_0,\tilde A_k) f_A^{(k)}(x_{\Phi,A}^{(k)})
\text{ (by choice of $f^{(k)}_A$)}\\
& = \sum_{A\in\cE_{\beta,\beta}(B)}\zeta(\beta^2,A)
\sum_{k=1}^{d_A}\zeta(\alpha_0,\tilde A_k) f_A^{(k)}\left(\sum_{m=1}^{s(\beta,\beta,A)}x_{\Phi,A}^{(m)}\right)
\text{ (since $\ran(f_A^{(k)})\subset\ran(x_{\Phi,A}^{(k)}$)}
\\
& = \sum_{A\in\cE_{\beta,\beta}(B)}\zeta(\beta^2,A)\sum_{k=1}^{d_A}\zeta(\alpha_0,\tilde A_k) f_A^{(k)}\left(x_A\right)
\text{ (by \eqref{E:12.1.1})}
\\
&\le \sum_{A\in\cE_{\beta,\beta}(B)}\zeta(\beta^2,A)\left\|\sum_{k=1}^{d_A}\zeta(\alpha_0,\tilde A_k)f_A^{(k)}\right\| <\vp
\text{ (by \eqref{these don't go together})}.
\end{align*}
This contradiction completes the proof.
\end{proof}
Before proving Corollary \ref{C:1.1} we will need the following observation.

\begin{prop}\label{P:13.5}
Let $X$ be a Banach space, $\alpha<\omega_1$ and $L$ be an infinite subset of the natural numbers so that there exist numbers $0<c<C$ and a map $\Phi:\cS_\alpha\cap[L]^{<\omega}\to X$
that is a $c$-lower-$d_{\alpha,\infty}$ and $C$-upper-$d_{\alpha,1}$ embedding. Then for every $\beta<\alpha$ there exists $n\in\N$ and a map $\Phi_\beta:\cS_{\beta}\cap[L\cap(n,\infty)]^{<\omega}\to X$ that is a
$c$-lower-$d_{\beta,\infty}$ and $C$-upper-$d_{\beta,1}$ embedding.
\end{prop}
Before proving Proposition \ref{P:13.5} we need some preliminary observation.
 \begin{lem}\label{L:13.6}
Let $\alpha$ be an ordinal number and $A\subset[0,\alpha]$ satisfying:
\begin{itemize}
 \item[(i)] $\alpha\in A$ and
 \item[(ii)] if $\beta\in A$ and $\gamma<\beta$, then there is $\gamma\le \eta<\beta$ with $\eta\in A$.
\end{itemize}
Then $A = [0,\alpha]$.
\end{lem}

\begin{proof}
If the conclusion is false, let $\alpha_0$ be the minimum ordinal for which $\text{ there exists }A_0\subsetneq[0,\alpha_0]\text{ satisfying (i) and (ii)}$.
Fix an arbitrary $\gamma<\alpha_0$. By the properties of $A_0$, there exists $\gamma\le\alpha<\alpha_0$ with $\alpha\in A_0$.
If we set $A = A_0\cap[0,\alpha]$, then it easily follows that $A$ and $\alpha$ satisfy (i) and (ii). By the minimality of $\alpha_0$, we conclude $A = [0,\alpha]$ and hence, since $\gamma \in[0,\alpha] = A$ and $A \subset A_0$, we have $\gamma\in A_0$.
Since $\gamma$ was chosen arbitrarily we conclude $[0,\alpha_0] = A_0$, a contradiction that completes the proof.
\end{proof}

\begin{lem}\label{L:13.7}
Let $\alpha<\omega_1$ be a limit ordinal number. Then there exists a sequence of successor ordinal numbers $(\mu(\alpha,n))_n$ satisfying:
\begin{itemize}
 \item[(i)] $\mu(\alpha,n)<\alpha$ for all $n\in\N$ and $\lim_n\mu(\alpha,n) = \alpha$,
 \item[(ii)] $\cS_\alpha = \{A\in[\N]^{<\omega}: A\in\cS_{\mu(\alpha,\min(A))}\}\cup\{\emptyset\}$ and $\cS_{\mu(\alpha,n)}\cap\big[[n,\infty)\big]^{<\omega}\subset\cS_\alpha$ for all $n\in\N$, and
 \item[(iii)] for $A\in\cS_{\alpha}\setminus\{\emptyset\}$, $z_{(\alpha,A)} = z_{(\mu(\alpha,\min(A)),A)}$.
\end{itemize}
\end{lem}

\begin{proof}
We define $(\mu(\alpha,n))_n$ by transfinite recursion on the set of countable limit ordinal numbers. For $\alpha = \omega$ we set $(\mu(\omega,n))_n = (\lambda(\omega,n))_n$. If $\alpha$ is a limit ordinal so that for all $\alpha'<\alpha$ the corresponding
sequence has been defined, set for each $n\in\N$
$$
\mu(\alpha,n)=
\left\{
	\begin{array}{ll}
		\lambda(\alpha,n)  & \mbox{if } \lambda(\alpha,n)\text{ is a successor ordinal number, or}\\
		\mu(\lambda(\alpha,n),n) & \mbox{ otherwise.}
	\end{array}
\right.
$$
The fact that (ii), (iii) and the first part of (i) hold is proved easily by transfinite induction using \eqref{E:2.5.1}  in Corollary \ref{C:2.5} and the definition of repeated averages. To show that $\lim_n\mu(\alpha,n) = \alpha$, we will show
that for arbitrary $L\in[\N]^{\omega}$ we have $\sup_{n\in L}\mu(\alpha,n) = \alpha$. Fix $L\in[\N]^{\omega}$ and $\beta<\alpha$. Then, since $\CB(\cS_\alpha\cap[L]^{<\omega}) = \omega^\alpha+1>\omega^\beta+1$, we have
$\emptyset\in(\cS_\alpha\cap[L]^{<\omega})^{(\omega^\beta+1)}$ and hence there exists $n\in L$ with $\{n\}\in(\cS_\alpha\cap[L]^{<\omega})^{(\omega^\beta)}$. Using \eqref{E:2.2} we obtain 
\begin{align*}
\emptyset&\in(\cS_\alpha\cap[L]^{<\omega})^{(\omega^\beta)}(\{n\}) \ksubset \cS_\alpha^{(\omega^\beta)}(\{n\}) \keq (\cS_\alpha(\{n\}))^{(\omega^\beta)} \keq (\cS_{\mu(\alpha,n)}(\{n\}))^{(\omega^{\beta})} 
\keq \cS_{\mu(\alpha,n)}^{(\omega^\beta)}(\{n\}),
\end{align*}
which implies $\{n\}\in\cS_{\mu(\alpha,n)}^{(\omega^\beta)}$,     i.e. $\CB(\cS_{\mu(\alpha,n)})\geqslant \omega^\beta + 1$ which yields $\mu(\alpha,n)\geqslant \beta$.
\end{proof}

\begin{rem}
It is uncertain whether $\cS_{\mu(\alpha,n)}\subset\cS_{\mu(\alpha,n+1)}$, it is even uncertain whether $\mu(\alpha,n) \le \mu(\alpha,n+1)$. It is true however that for $2\le n\le m$ we have $\mu(\alpha,n)\le \mu(\alpha,m) + 1$.
\end{rem}

\begin{proof}[Proof of Proposition \ref{P:13.5}]
We shall first treat two very specific cases. In the first case, $\alpha = \beta +1$. Fix $n_0\geqslant 2$ with $n_0\in L$ and $B_0\in\MAX(\cS_\beta\cap[L]^{<\omega})$ with $\min(B_0) = n_0$. Define $n = \max(B_0)$
and $\Phi_\beta:\cS_{\beta}\cap[L\cap(n,\infty)]^{<\omega}\to X$ with $\Phi_\beta(A) = n_0\Phi(B_0\cup A)$. Then $\Phi_\beta$ is the desired embedding.

In the second case, $\alpha$ is a limit ordinal and for some $n_0\geqslant 2$ with $n_0\in L$ we have $\beta + 1 = \mu(\alpha,n_0)$.  Fix $B_0\in\MAX(\cS_\beta\cap[L]^{<\omega})$ with $\min(B_0) = n_0$, define
$n = \max(B_0)$ and $\Phi_\beta:\cS_{\beta}\cap[L\cap(n,\infty)]^{<\omega}\to X$ with $\Phi_\beta(A) =  n_0\Phi(B_0\cup A)$. Then, using the properties of $(\mu(\alpha,k))_k$,
it can be seen that $\Phi_\beta$ is well defined and it is the desired embedding.

In the general case, define $A$ to be the set of all $\beta\le \alpha$ for which such an $n$ and $\Phi_\beta$ exist. Since $\alpha\in A$, it remains to show that $A$ satisfies (ii) of Lemma \ref{L:13.6}.
Indeed, fix $\beta\in A$ and $\gamma <\beta$. If $\beta = \eta+1$, then by the first case we can deduce that $\eta\in A$ and $\gamma\le \eta < \beta$. Otherwise, $\beta$ is a limit ordinal.
Let $\Phi_\beta$ and $n_\beta$ witness the fact that $\beta\in A$ and by Lemma \ref{L:13.7} (i) we may choose $n\in L$ with $n > n_\beta$ so that $\mu(\beta,n)>\gamma+1$.
If $\eta$ is the predecessor of $\mu(\beta,n)$, then by the second case we deduce that $\eta\in A$ and $\gamma\le \eta<\beta$.
\end{proof}

\begin{proof}[Proof of Corollary \ref{C:1.1}] We first  recall a result by Causey \cite[Theorem 6.2]{Ca} which says that for  a countable ordinal $\xi$ it follows that   $\gamma=\omega^\xi$ is the Szlenk index of some separable Banach space $X$,
 if and only if  $\xi$ is not of the form $\xi=\omega^\eta$, with $\eta$ being a limit ordinal. Since $\alpha=\omega^{\omega^{\omega^\alpha}}$, $\alpha$ cannot be the Szlenk index of some separable Banach space.
 
 \noindent``(a)$\Rightarrow$(b)'' From (a) and  Causey's result we have $\Sz(X)\kle\omega^\alpha$ and  $\Sz(X^*)\kle \omega^\alpha$, and thus there exists a $\xi<\omega_1$
with  $\beta=\omega^{\omega^{\xi}}\le \omega^{\omega^{\xi+1}}<\alpha$, so that $\Sz(X)<\omega^{\beta}$ and $\Sz(X)<\beta$.
  It follows therefore  from Theorem \ref{T:12.6} that for no
  $L\in[\N]^{\omega}$  there exist numbers $0<c<C$ and a map $\Phi:\cS_{\beta^2}\cap[L]^{<\omega}\to X$
that is a $c$-lower-$d_{\beta^2,\infty}$ and $C$-upper-$d_{\beta^2,1}$ embedding. 
 Since $\beta^2\le  \omega^{\omega^{\xi+1}}<\alpha$ we  conclude our claim from Proposition \ref{P:13.5}

\noindent
 ``(b)$\Rightarrow $(a)''  follows from Theorems \ref{T:9.1} and \ref{T:9.3}.
\end{proof}

In order to proof Corollary \ref{C:1.2}
recall that every separable Banach space is isometrical equivalent to a subspace of $C[0,1]$,
the space of continuous functions on $[0,1]$. The set $ \cSB$ of all
closed subspaces of $C[0,1]$ is given the {\em Effros-Borel structure},
which is the $\sigma$-algebra generated by the sets $\{
F\kin\cSB:\, F\cap U\kneq\emptyset\}$, where $U$ ranges over
all open subsets of $C[0,1]$.
\begin{proof}[Proof of Corollary \ref{C:1.1}] By \cite[Theorem D]{OSZ1} the set
$$\cC_\alpha=\{ X\in  \cSB: X\text{ reflexive, and }\max(\Sz(X),\Sz(X^*))\le \alpha\}$$
is analytic. So it is left to show that its complement is also analytic.
Since by Corollary \ref{C:1.2}
\begin{align*}
\cSB\setminus \cC_\alpha&=\{ X\kin\cSB: X\text{ not reflexive}\}\cup \{ X\kin  \cSB:  X\text{reflexive and}\max(\Sz(X),\Sz(X^*))> \alpha\}\\
&=\{ X\kin\cSB: X\text{ not reflexive}\}\cup \{ X\kin  \cSB:  (\cS_{\alpha},d_{1,\alpha}) \text{ bi-Lipschitzly embeds into } X\}
\end{align*}
and since by \cite[Corollary 3.3]{Bos} the set of reflexive spaces in $\cSB$ is co-analytic, we deduce our claim from the following  Lemma.
This Lemma seems to be wellknown, but for completeness we include a short  proof.
\end{proof}
\begin{lem}\label{P13.6} Let $(M,d)$ be a separable metric space. Then
$$\cA(M)=\{X\in\cSB: (M,d) \text{ bi-Lipschitzly embeds into $X$}\} $$ 
is analytic in $\cSB$.
\end{lem}
\begin{proof} After passing to a countable dense subset of $M$ we can assume that $M$ is countable and we  enumerate $M$ into $\{m_n:n\kin\N\}$.
The set $\big(\N^\N\big)^\N$, \ie the set of  all sequences in $\N^\N$ is a Polish space with respect to the product topology.
We write an element $s\in \big(\N^\N\big)^\N$ as  $s=(s_n:n\in\N)$ with $s_n=(s_{n,j})_{j=1}^\infty\in \N^\N$.
By the Kuratowski-Ryll-Nardzewski Selection Theorem \cite[Theorem 12.13]{Ke} there are Borel maps $d_n:\cSB\to C[0,1]$, $n\in\N$, so that for every $X\in\cSB,$ the sequence
$\big(d_n(X)\big)$ is in $X$ and  dense in $X$. Put $D(X)=\{d_n(X): n\in\N\}$, for $X\in \cSB$.
We first note that there is a bi-Lipschitz map $\Psi:M\to X$ if and only if 
\begin{align*}
\exists C\kgr1\,&\forall n\kin\N\,\exists \Psi_n:M\to D(X), \text{ so that } (\Psi_n(m)) \text{ is Cauchy in $X$ for all $m\kin M$ and}\\
&\frac1C d(m,m')<\lim_{n\to\infty}\| \Psi_n(m)-\Psi_n(m')\|<Cd(m,m') \text{ for all $m,m'\in M$.}
\end{align*}
Since the set 
$$\left \{ (s,X)\kin \big(\N^\N\big)^\N\times\cSB: \begin{matrix}
\forall k\kin\N\,\forall l\kin \N\, \exists n_0\kin\N\,\forall n_1,n_2\kge n_0\quad \|d_{s_{n_1,k}}(X)-d_{s_{n_2,k}}(X)\|<\frac1l\\
\exists C\kin\N \, \forall k,l\kin\N \, \exists n_0\kin\N\,\forall n\kge n_0\qquad\qquad\qquad\qquad\qquad\qquad\qquad\\
 \qquad \frac1C d(m_k,m_l)< \|  d_{s_{n,k}}(X)-d_{s_{n,l}}(X) \|<Cd(m_k,m_l)
\end{matrix} \right\}  $$
is Borel in the product space $\cSB\times \big(\N^\N\big)^\N$, its projection on the second coordinate is analytic. On the other hand the image of  this projection
consists exactly of the spaces $X\in\cSB$ for which there is a bi-Lipschitz map $\Psi: M\to X$ (consider for each $s\in \big(\N^\N\big)^\N$, $n\in\N$ and $X\in\cSB$, the map
$\Psi_n:M\ni m_k\mapsto d_{s_{n,k}}(X)\in X$).
\end{proof}

\section{Final comments and open questions}\label{S:13}

The proof of Theorem A, yields the following equivalences.  ``(a)$\Rightarrow$(b) '' follows from Proposition \ref{P:8.2},
``(d)$\Rightarrow$(a)'' from Theorem \ref{T:12.6}, while ``(b)$\Rightarrow$(c)$\Rightarrow$(d)'' is  trivial.
\begin{cor}\label{C:13.1} For a separable  Banach space $X$  the following statements
are equivalent
\begin{enumerate}
\item[a)] $X$ is not reflexive,
\item[b)] For all $\alpha<\omega_1$ there exists  for some numbers $0<c<C$, a  $c$-lower $d_{\infty,\alpha}$, $C$-upper $d_{1,\alpha}$ embedding of $\cS_\alpha$ into $X$. 
\item[c)]  For all $\alpha<\omega_1$  there is a map $\Psi_\alpha:\cS_\alpha\to X$ and some $0<c\le 1$, so that 
for all  $A,B,C\in\cS_\alpha$,  with the property that  $A\succeq C$, $B\succeq C$ and $A\setminus C< B\setminus C$
 $$c d_{1,\alpha}(A,B)\le \|\Psi(A)-\Psi(B)\| \le d_{1,\alpha}(A,B).$$
\item[d)] For all $\alpha<\omega_1$, there is an $L\in[\N]^{<\omega}$ and a semi embedding  $\Psi_\alpha:\cS_\alpha\cap[L]^{<\omega}\to X$.
\end{enumerate}  
\end{cor}

As mentioned before we can consider  for $\alpha<\omega_1$ and $A\in\cS_\alpha$ the vector
$z_A$ to be an element in   $\ell_1^+$, with $\|x\|_{\ell_1}\le 1$.
 We define 
 $$T_\alpha=\{ (A,B)\in \cS_\alpha\times\cS_\alpha : \exists C\preceq A \text{ and } C\preceq A, \text{ with } A\setminus C< B\setminus C\}.$$
 We note that  $d_{1,\alpha}(A,B)= \|z_A-z_B\|_1$, for $(A,B)\in T_\alpha$.
Using this notation we deduce the following sharpening of \cite[Theorem 3.1]{Os3} from Corollary \ref{C:13.1}

\begin{cor}\label{C:13.2}
Let $X$ be a separable Banach space. Then the following are equivalent:
\begin{enumerate}
\item[a)] $X$ is not reflexive.
\item[b)] For all $\alpha<\omega_1$  there is a map $\Psi_\alpha:\cS_\alpha\to X$ and some $0<c\le 1$, so that 
$$c d_{1,\alpha}(A,B)\le \|\Psi_\alpha(A)-\Psi_\alpha(B)\|\le d_{1,\alpha}(A,B), \text{ whenever $(A,B)\in \cT_\alpha$.}$$
\end{enumerate}
\end{cor}

Next we would like to show that  the James space $J$ (c.f. \cite[Definition 4.43]{FHHMZ}) has the property that $\Sz(J)=\Sz(J^*)=\omega$, and thus,  since $J$ is not reflexive, that Theorem C becomes false without the assumption  that $X$ is reflexive.

If $X$ is a Banach space with an FDD and $1<p \le \infty$. The space $X$ is said to satisfy an upper $\ell_p$-estimate, if there
is a constant $C$ so that for every block sequence $(x_k)_{k=1}^n$ in $X$ we have $C\|\sum_{k=1}^nx_k\|\le (\sum_{k=1}^n\|x_k\|^p)^{1/p}$. A lower $\ell_p$-estimate is defined in the obvious way.
If $X$ is a Banach space with an FDD that satisfies an upper $\ell_p$-estimate, note that by a standard argument, the FDD of $X$ has to be shrinking and hence $X^*$ is separable.
Therefore the assumption that $X$ has an FDD with an upper $\ell_p$-estimate implies by \cite[Theorem 3.4]{OS2} that $\Sz(X)=\omega$.

\begin{ex}\label{Ex:13.3}
Let  $J$ denote James space, the above implies that $\Sz(J) = \Sz(J^*) = \omega$. To see that, we recall that $J$ has two bases $(d_i)_i$ and $(e_i)_i$, the first one of which is shrinking and satisfies an upper $\ell_2$-estimate, whereas the second one is boundedly complete and satisfies a lower $\ell_2$-estimate.
The first fact yields that $\Sz(J)\le \omega$, whereas the second fact and a duality argument yields that the space $X$ spanned by the basis $(e_i^*)_i$ satisfies an upper $\ell_2$-estimate, i.e. $\Sz(X) \le \omega$. It is well known that $X$ is of co-dimension one in $J^*$ and that $J^*$ is isomorphic to its hyperplanes and therefore $\Sz(J^*)\le \omega$.

\end{ex}
This Example shows that Theorem C cannot be true if we omit the requirement that $X$ is reflexive. 
Nevertheless the following variation of Corollary \ref{C:1.1} holds.
\begin{cor}\label{C:13.4}
Assume that $\alpha<\omega_1$ is an ordinal, for which $\alpha=\omega^\alpha$. Then the following statements are equivalent for a
separable Banach space $X$.
\begin{enumerate}
\item[a)] $X$ is reflexive and $\max(\Sz(X),\Sz(X^*))\le \alpha$,
\item[b)] There is no map $\Psi:\cS_\alpha \to X$ and $0<c\le 1$ so that for all $A,B,C\in\cS_\alpha$ with the property that $C\preceq A$, $C\preceq B$, and $A\setminus C< B\setminus C$, we have
$$cd_{1,\alpha}(A,B)\le \|\Psi(A)-\Psi(B)\|\le d_{1,\alpha}(A,B).$$
\end{enumerate}
\end{cor}
\begin{proof} Let $\Psi:\cS_\alpha \to X$ satisfying the condition stated in (b) for some $c>0$. Then $\tilde\Psi=\big(\Psi-\Psi(\emptyset)\big)/2$  also has this property for $c/2$,  maps $\cS_\alpha$ into $B_X$ and $\tilde\Psi(\emptyset)=0$.
\noindent  ``(a)$\Rightarrow$(b)'' follows from Theorem \ref{T:12.6},  the in the proof of Corollary \ref{C:1.1} cited result from \cite{Ca},
and the same argument involving Proposition \ref{P:13.5} in the proof of Corollary \ref{C:1.1}.
``(b)$\Rightarrow$ (a)'' follows from Proposition \ref{P:8.2}, Theorems \ref{T:9.1} and \ref{T:9.3}.
\end{proof} 
\begin{rem} The statement of Corollary \ref{C:13.4}  also holds for $\alpha=\omega$. This can be seen from the proof of the main result in \cite{BKL}.
\end{rem}

 We finish with stating two open problems.
 \begin{probs}\label{Pr:13.8} \begin{enumerate}
 \item[a)\qquad\,\,] \hskip-.9cm Does there exists a family of metric spaces $(M_i,d_i)$ which is  a family of test spaces for reflexivity in the sense of \cite{OS2}, \ie for which it is true that  a  separable Banach 
 space $X$  is reflexive if and only if not all of the $M_i$ uniformly bi-Lipschitzly embed into $X$?
 \item[b)] Does there exists a countable family of metric spaces $(M_i,d_i)$ which is  a family of test spaces for reflexivity?
 \item[c)] It follows from Theorem B, that if   $X$ is a separable Banach space with non separable bidual  then 
 $(\cS_{\alpha},d_{1,\alpha})$ bi-Lipschitzly embeds into $X$, for all $\alpha<\omega_1$.
 Is the converse true, or in Ostrovskii's language, are the spaces   $(\cS_{\alpha},d_{1,\alpha})$, $\alpha<\omega_1$ test spaces for spaces with  separable bi-duals?
 \end{enumerate}
 \end{probs}

\end{document}